\documentclass[11pt]{article}
\usepackage{graphics,amsmath,amssymb}
\usepackage{latexsym}
\usepackage{epsfig}

\def\a{\alpha}
\def\be{\beta}

\def\epsilon{\varepsilon}
\def\ga{\gamma}

\def\la{\lambda}
\def\La{\Lambda}
\def\phi{\varphi}
\def\si{\sigma}

\def\om{\omega}

\newtheorem{theorem}{Theorem}[section]
\newtheorem{lemma}[theorem]{Lemma}

\newtheorem{corollary}[theorem]{Corollary}

\newtheorem{proposition}[theorem]{Proposition}
\newtheorem{remark}[theorem]{Remark}

\def\Z{{\mathbb Z}}
\def\N{{\mathbb N}}
\def\C{{\mathbb C}}
\def\R{{\mathbb R}}
\def\D{{\cal D}}
\def\I{{\mathbb I}}

\newenvironment{Proof}{\removelastskip\par\medskip
\noindent{\em Proof.} \rm}{\penalty-20\null$\square$\par\medbreak}

\newenvironment{Proofy}{\removelastskip\par\medskip
\noindent{\em Proof of Theorem~\ref{th:inv.ingham}.} \rm}{\penalty-20\null$\square$\par\medbreak}

\newenvironment{Proofv}{\removelastskip\par\medskip
\noindent{\em Proof of Theorem~\ref{th:dir.ingham}.} \rm}{\penalty-20\null$\square$\par\medbreak}

\newenvironment{Proof1}{\removelastskip\par\medskip
\noindent{\em Proof of Theorem~\ref{th:diringham}.} \rm}{\penalty-20\null$\square$\par\medbreak}

\newenvironment{Proof2}{\removelastskip\par\medskip
\noindent{\em Proof of Proposition~\ref{pr:haraux-inv}.} \rm}{\penalty-20\null$\square$\par\medbreak}

\newenvironment{Proof3}{\removelastskip\par\medskip
\noindent{\em Proof of Theorem~\ref{th:inv.ingham1}.} \rm}{\penalty-20\null$\square$\par\medbreak}     

\newenvironment{ProofL1}{\removelastskip\par\medskip
\noindent{\em Proof of Lemma~{\rm\ref{le:stimaK}}.} \rm}{\penalty-20\null$\square$\par\medbreak}

\title{\bf Control problems for weakly coupled \\
systems with memory}

\author{Paola Loreti
\thanks{Dipartimento di Scienze di Base e Applicate per l'Ingegneria
Sezione di Matematica,
Sapienza Universit\`a di Roma,
Via Antonio Scarpa 16, 00161 Roma (Italy); e-mail: $<$paola.loreti@sbai.uniroma1.it$>$ }
\and Daniela Sforza
\thanks{Dipartimento di Scienze di Base e Applicate per l'Ingegneria
Sezione di Matematica, Sapienza Universit\`a di Roma,
Via Antonio Scarpa 16, 00161 Roma (Italy); e-mail: $<$daniela.sforza@sbai.uniroma1.it$>$ }}

\begin{document}
\date{}

\maketitle

\begin{abstract}
We investigate control problems for {\it wave-Petrovsky} coupled systems in the presence of memory terms. By writing the solutions as Fourier series, we are able to prove Ingham type estimates, and hence reachability results. Our findings  have applications in viscoelasticity theory and linear acoustic theory.
\end{abstract}

\bigskip
\noindent
{\bf Keywords:} coupled systems; convolution kernels; Fourier series; Ingham estimates; reachability.

\section{Introduction}

We will analyze control problems for {\it wave-Petrovsky} weakly coupled systems in the presence of memory terms. In particular, we will solve the reachability  to a given target in a finite time, by using a harmonic approach based
on Ingham type estimates.  

In the papers \cite {LoretiSforza, LoretiSforza1} we studied reachability problems for a class of integro-differential equations
\begin{equation*}
\displaystyle 
u_{tt}(t,x) -u_{xx}(t,x)+\beta\int_0^t\ e^{-\eta(t-s)} u_{xx}(s,x)ds= 0\,,
\qquad
t\in (0,T)\,,\ \  x\in(0,\pi)\,,
\end{equation*}
then generalized to spherical domains in \cite{LoretiSforza2} and to more general kernels in \cite{LoretiSforzaP}.


The interest for researching this type of control problems comes from the theory
of viscoelasticity. Exponential kernels naturally arise in linear viscoelasticity theory, such as in the analysis of Maxwell fluids  or 
Poynting-Thomson solids,
see  e.g.
\cite {Pruss,Re1}.  For other references in viscoelasticity theory see the seminal papers of Dafermos \cite{D1,D2} and  \cite{RHN,LPC}. 
For other type of kernels, see \cite{LV}. 

As it is well known,  viscoelastic relaxation kernels have to be completely monotone functions, that is,
continuously differentiable to every order functions $K(t)$ satisfying 
$$
(-1)^n K^{(n)}(t)\geq 0\,\qquad   \forall  n\in\N\,, \, \, \forall t\geq 0\,.
$$ 
This class of  relaxation kernels includes, 
 as a significant case, the Prony sum
$$K(t)=\sum_{ i=1}^N \beta_i e^{ -\eta_i t}$$
with $ \beta_i>0$ and $ \eta_i\geq 0$,  $i=1, \dots,N$. 
Prony-sum kernels have many implications for the dispersion and the attenuation phenomena in acoustic  theory  \cite{Han, Han1, HaS1}. Moreover,  the analysis of the 1-$d$ wave equation of a 
vibrating string has  analogies with  seismic wave propagation  \cite {StWy}.
It could be interesting  to consider in the model  the effect of  viscosity as an attenuation phenomenon
for seismic events.

Continuing along the lines traced by the research papers  \cite {LoretiSforza, LoretiSforza1, LoretiSforza2}, we have done  further investigations,  which split  into   the following  three directions $a)$, $b)$ and $c)$.
\begin{itemize}
\item[$a)$]
The study of a more general relaxation kernel of Prony type in a single wave equation. This problem presents some difficulties with
respect to the case of kernels consisting in a single exponential function, because we have to handle  a more
complicated spectral analysis,  to compare the coefficients of the materials
and to find conditions under which the reachability control problem may have a positive
solution (\cite  {LoretiSforza3}, in preparation).
\item[$b)$]
The analysis of weakly coupled systems  of  {\it wave-wave} type, with a memory term having a single-exponential kernel
as in  \cite { LoretiSforza1}. 
 To find the eigenvalues, one has to study  a fifth-degree equation: it turns out  that the two couples of complex conjugate roots have the same asymptotic behavior (\cite { LoretiSforza4}, preprint).
See \cite{KR}  for one of the first papers on wave-wave coupled PDE's without memory.

\item[$c)$]
The study of weakly coupled systems   
 of {\it wave-Petrovsky} type, again with  memory terms consisting in a single-exponential kernel.
The analysis of weakly   coupled PDE's of wave-Petrovsky type without memory began in \cite{KL1},
where the harmonic analysis approach was successfully applied to get osservability results.
\end{itemize}

\noindent 
All these research lines need a deep analysis and extensive computations,  with  significant differences. In this paper  we consider
the third  research problem $c)$.
We add to a wave equation an integral relaxation term and couple it with a Petrovsky type equation in 
the following way
\begin{equation}\label{eq:problem-uI}
\begin{cases}
\displaystyle 
u_{1tt}(t,x) -u_{1xx}(t,x)+\beta\int_0^t\ e^{-\eta(t-s)} u_{1xx}(s,x)ds+Au_2(t,x)= 0\,,
\\
\phantom{u_{1tt}(t,x) -u_{1xx}(t,x)+\int_0^t\ k(t-s) u_{1xx}(s,x)ds+Au_2(t,x)= 0}
t\in (0,T)\,,\quad x\in(0,\pi)
\\
\displaystyle
u_{2tt}(t,x) +u_{2xxxx}(t,x)+Bu_1(t,x)= 0
\,,
\end{cases}
\end{equation}
$0<\beta<\eta$, $A\,,B\in\R$,
with null initial conditions 
\begin{equation}
u_1(0,x)=u_{1t}(0,x)=u_2(0,x)=u_{2t}(0,x)=0\qquad  x\in(0,\pi)\,,
\end{equation} 
and boundary conditions
\begin{equation}\label{eq:bound-u1i}
u_1(t,0)=0\,,\qquad u_1(t,\pi)=g_1(t)\qquad t\in (0,T) \,,
\end{equation}
\begin{equation}\label{eq:bound-u2i}
u_2(t,0)=u_{2xx}(t,0)=u_2(t,\pi)=0\,,\qquad u_{2xx}(t,\pi)=g_2(t)\qquad t\in (0,T) \,.
\end{equation}
We can consider $g_i$, $i=1,2$, as control functions.  The reachability problem consists in proving the existence of  $g_i\in L^2(0, T)$  that steer a weak solution of system \eqref{eq:problem-uI}, subject to boundary conditions \eqref{eq:bound-u1i}--\eqref{eq:bound-u2i}, from the null state to a given one in finite time. To better explain,
we embrace the definition of reachability problem for systems with memory given by several authors in the literature, see for example \cite{Lio3,K1,K2,Las,L0,L, LV, MN}.  

Indeed,
we mean the following:  given  $T>0$ and 
$$(u_{10},u_{11},u_{20},u_{21})\in  L^{2}(0,\pi)\times H^{-1}(0,\pi)\times H_0^{1}(0,\pi)\times H^{-1}(0,\pi)\,,$$
to find   $g_i\in L^2(0,T)$ such that the weak solution  $u$ of problem \eqref{eq:problem-uI}--\eqref{eq:bound-u2i} 
verifies the final conditions
\begin{equation}\label{eq:problem-u1i}
u_1(T,x)=u_{10}(x)\,,\qquad u_{1t}(T,x)=u_{11}(x)\,,
\qquad x\in(0,\pi)\,,
\end{equation}
\begin{equation}\label{eq:problem-u1ii}
u_2(T,x)=u_{20}(x)\,,\qquad u_{2t}(T,x)=u_{21}(x)\,,
\qquad x\in(0,\pi)\,.
\end{equation}

We are able to bring about reachability results 
without any smallness assumption
on the convolution kernels, as suggested by J.-L. Lions in \cite[p. 258]{Lio3}.
A common way to study exact controllability
problems is the so-called Hilbert Uniqueness Method,
introduced by Lagnese -- Lions, see \cite{LaLi,Lio1,Lio2,Lio3}.
We  will apply this method to  system \eqref{eq:problem-uI}.
The HUM method is based on  a ``uniqueness theorem" for the adjoint problem.
To prove such uniqueness theorem we will employ some typical techniques  of harmonic analysis, see \cite{Ru,KL2}. This approach
relies on Fourier series development for the solution $(u_1,u_2)$ of the adjoint
problem, which can be written as follows
\begin{equation}\label{eq:u1I}
u_1(t,x)=\sum_{n=1}^{\infty}\Big(R_ne^{r_nt}+C_ne^{i\om_nt}+\overline{C_n}e^{-i\overline{\omega_n}t}
+D_{n}e^{i p_{n}t}+\overline{D_{n}}e^{-i \overline{p_n}t}\Big)\sin(nx)
\,,
\end{equation}
\begin{equation}\label{eq:u2I}
u_2(t,x)=\sum_{n=1}^{\infty}
\Big(d_nD_{n}e^{i p_{n}t}
+\overline{d_n}\overline{D_{n}}e^{-i \overline{p_n}t}\Big)\sin(nx)
+ e^{-\eta t}\sum_{n=1}^{\infty}\Re p_{n}
\Big(\frac{D_{n}}{\eta+ip_n}+\frac{\overline{D_{n}}}{\eta-ip_n}
\Big)\sin(nx)
\,,
\end{equation}
where 
\begin{equation}\label{eq:rnI}
r_n=\be-\eta-{\be(\be-\eta)^2\over\la_{n}}
+O\Big({1\over{\la_{n}^{3/2}}}\Big)\,,
\end{equation}
\begin{equation}
\omega_n=
\sqrt{\lambda_{n}}+{\beta\over2}\Big({3\over4}\beta-\eta
\Big){1\over\sqrt{\lambda_{n}}}
+i\Big({\be\over 2}-{\be(\be-\eta)^2\over2\la_{n}}\Big)
+O\Big({1\over{\lambda_{n}^{3/2}}}\Big)
\,,
\end{equation}
\begin{equation}\label{eq:pnI}
p_n=\la_{n}+{AB\over 2\la_{n}^3}+O\Big({1\over{\la_{n}^{4}}}\Big)
\,.
\end{equation}
In this framework  Ingham type estimates \cite{Ing} play an important role. 
We need to establish inverse and direct inequalities for functions \eqref{eq:u1I}-\eqref{eq:u2I} evaluated at $x=\pi$, see
\eqref{eq:inv.ingham} and \eqref{eq:diringham} later on, obtaining them in the same sharp time of the nonintegral case.

In this approach the main difficulties are the following:
\begin{enumerate}
\item
The study of the distribution of the eigenvalues on the complex
plane. Indeed,
the spectral analysis of the coupled system leads to a full fifth-degree 
equation governing  the eigenvalues behavior.  
A method due to Haraux \cite{Ha}, subsequent to
the seminal work of Ingham \cite{Ing}, enables us to consider only the
asymptotic behavior of the eigenvalues related to the spatial operator.
In order to get the asymptotic behavior  of the eigenvalues, see \eqref{eq:rnI}-\eqref{eq:pnI}, we need to develop an accurate  asymptotic computation (see Section \ref{se:specana}).
\item
The generalization of  Ingham's approach and the proof of the inverse inequality.
This means that  we are able to  generalize  the results contained in \cite{Ing} and \cite{Ha},
see Theorems  \ref{th:inv.ingham1}, \ref{th:diringham} and 
Proposition \ref{pr:haraux-inv} later on.
In particular, a difficulty is
the presence in $u_2$ of a series constant in time, but depending on the coefficients $D_n$, see \eqref{eq:u2I}.
Due to its form, this series is difficult to handle. To overcome this impasse, as a first step we can neglect the dependence on $D_n$ and treat the whole series simply as a constant. Following this approach, we have to use Haraux's method: we introduce the usual operator which annihilates the constant, so that we can apply an inverse estimate holding in the case the constant is null, and then  recover the constant itself, see Theorem \ref{th:Dn+D} later on.

\item
Due to the finite speed of propagation, we expect the controllability time $T$ to be sufficiently large. Indeed, we will find that   $T>2\pi/\gamma$, where $\gamma$ is the gap of a branch of eigenvalues related to the integro-differential operator, see 
Theorem \ref{th:reachres}. The achievement of the time estimate    $T>2\pi/\gamma$ will require an accurate compensation in the analysis of  the  terms appearing in formulas \eqref{eq:u1I} and \eqref{eq:u2I}, see Theorem \ref{th:Inversa} later on. 
 
\end{enumerate}

The plan of our paper is the following. In Section 2 we  give some preliminary results. In Section 3 we describe the Hilbert Uniqueness Method in an abstract setting. In Section 4 we give a detailed spectral analysis for a coupled system with memory.
In Section 5 we prove our main results:
Theorems \ref{th:inv.ingham1}, \ref{th:diringham} and Proposition \ref{pr:haraux-inv}.
Finally, in Section 6 we give a reachability result for a coupled system with memory.

\section{Preliminaries}
Let $X$ be a real Hilbert space with scalar product
$\langle \cdot \, ,\, \cdot \rangle$ and norm $\| \cdot \|$. 
For any
$T\in\, (0,
\infty]$   we denote by
$L^1(0,T;X)$ the usual spaces of measurable
functions
$v:(0,T)\to X$ such that one has

$$
\|v\|_{1,T}:=\int_0^T \|v(t)\|\,dt<\infty\,.
$$
 We shall use the shorter notation $\|v\|_1$ for
$\|v\|_{1,\infty}$.
We denote by $L_{loc}^1 (0,\infty;X)$ the space of functions
belonging to
$L^1(0,T;X)$ for any $T\in (0,\infty)$.
In the case of $X=\R$, we will use the abbreviations
$L^1(0,T)$ and
$L_{loc}^1(0,\infty)$ to denote the spaces $L^1(0,T;\R)$
and
$L_{loc}^1(0,\infty;\R)$,
respectively.

Classical results for integral equations
(see, e.g., \cite[Theorem 2.3.5]{GLS})
ensure that, for any  kernel $k\in L_{loc}^1(0,\infty)$ and 
$\psi\in L_{loc}^1 (0,\infty;X)$,  the problem
\begin{equation}\label{integral}
\varphi(t)-k*\varphi(t)=\psi(t),\qquad t\ge 0\,,
\end{equation}
admits a unique solution $\varphi\in L_{loc}^1(0,\infty;X)$. In particular, if we take $\psi=k$ in \eqref{integral}, 
we can consider the unique
solution $\varrho_k\in L_{loc}^1(0,\infty)$ of 
\begin{equation*}
\varrho_k (t)-k*\varrho_k (t)=k (t),\qquad t\ge 0\,.
\end{equation*}
Such a solution is called the {\em resolvent kernel} of $k$.
Furthermore, for any $\psi$ the solution $\varphi$ of (\ref{integral}) is given by the
variation of constants formula
\begin{equation*}
\varphi(t)=\psi(t)+\varrho_k *\psi(t),\qquad t\ge 0\,,
\end{equation*}
where $\varrho_k$ is the resolvent kernel of $k$.

We recall some results concerning integral equations in case of decreasing exponential kernels, see for example \cite[Corollary 2.2]{LoretiSforza1}.
\begin{proposition}\label{pr:unicita}
For $0<\beta<\eta$ and $T>0$ the following properties hold true.
\begin{description}
\item [(i)] The resolvent kernel of $k(t)=\beta e^{-\eta t}$ is $\varrho_k(t)=\beta e^{(\beta-\eta)t}$.

\item [(ii)]Given $\psi\in
L_{loc}^1 (-\infty,T;X)$, a function $\varphi\in L_{loc}^1 (-\infty,T;X)$ is a solution of
\begin{eqnarray*}
\varphi(t)-\beta\int_t^T\ e^{-\eta(s-t)}\varphi(s)ds=\psi(t)
\qquad t\le T\,,
\end{eqnarray*}
if and only if
\begin{eqnarray*}
\varphi(t)=\psi(t)+\beta\int_{t}^Te^{(\beta-\eta)(s-t)}\psi(s)\ ds
\qquad t\le T\,.
\end{eqnarray*}
Moreover, there exist two positive constants $c_1\,,c_2$ depending on $\beta,\eta,T$ such that
\begin{equation}\label{eq:unicita}
c_1\int_0^T|\varphi(t)|^2\ dt
\le\int_0^T|\psi(t)|^2\ dt
\le c_2\int_0^T|\varphi(t)|^2\ dt\,.
\end{equation}

\end{description}

\end{proposition}
\begin{lemma}\label{le:fifth}
Given $\beta\,,\eta\in\R$ and $\la\,,A\,,B\in\R\setminus\{0\}$, a couple $(f, g)$ of functions belonging to $ C^2([0,\infty))$ is a solution of the system
\begin{equation}\label{eq:system0}
\begin{cases}\displaystyle
f{''}(t)+\la f(t)-\la\be \int_0^t e^{-\eta(t-s)}f(s) ds+Ag(t)=0\,,\qquad t\ge 0\,,
\\\displaystyle
g{''}(t)+\la^2 g(t)+Bf(t)=0
\,,\qquad t\ge 0\,,
\end{cases}
\end{equation}
if and only if $f\in C^5([0,\infty))$ is  a solution of the problem
\begin{equation}\label{eq:fifth}
\begin{cases}
\displaystyle 
f^{(5)}(t)+\eta f^{(4)}(t)+(\la+\la^2) f{'''}(t)+(\la (\eta- \be)+\la^2\eta) f''(t)+(\la^3-AB)f'(t)
\\\displaystyle\hskip5cm
+(\la^3 (\eta- \be)-\eta AB)f(t)=0
\,,\qquad t\ge 0\,,
\\\\\displaystyle
 f^{(4)}(0)=-(\la+\la^2) f{''}(0)+\la \be f'(0)-(\eta\la \be +\la^3 -AB)f(0)\,,
\end{cases}
\end{equation}
and $g\in C^2([0,\infty))$ is given by
\begin{equation}\label{eq:Ag}
g(t)=-\frac1A\big[f{''}(t)+\la f(t)-\la\be \int_0^t e^{-\eta(t-s)}f(s) ds\big]\,.
\end{equation}
\end{lemma}
\begin{Proof}
Let $(f, g)$ be a solution of (\ref{eq:system0}). 
Differentiating the first equation in (\ref{eq:system0}), we get
\begin{equation}\label{eq:0gprime}
f{'''}(t)+\la f{'}(t)
+\eta\la\be \int_0^t e^{-\eta(t-s)}f(s) ds-\la \be f(t)+Ag'(t)=0\,,
\end{equation}
whence
\begin{equation}\label{eq:gprime0}
Ag'(0)=-f{'''}(0)-\la f{'}(0)+\la \be f(0)\,.
\end{equation}
Substituting in \eqref{eq:0gprime} the  identity
\begin{equation*}
\la\be \int_0^t e^{-\eta(t-s)}f(s) ds=f{''}(t)+\la f(t)+Ag(t)\,,
\end{equation*}
we obtain
\begin{equation}\label{eq:gprime}
f{'''}(t)+\eta f{''}(t)+\la f{'}(t)+\la (\eta- \be) f(t)+Ag'(t)+\eta Ag(t)=0\,.
\end{equation}
Differentiating yet again, we have
\begin{equation*}
f^{(4)}(t)+\eta f{'''}(t)+\la f{''}(t)+\la (\eta- \be) f'(t)+Ag''(t)+\eta Ag'(t)=0\,,
\end{equation*}
whence, by using the second equation in (\ref{eq:system0}), 
that is
$Ag{''}(t)=-ABf(t)-\la^2 Ag(t)$,
we get
\begin{equation}\label{eq:fourth}
f^{(4)}(t)+\eta f{'''}(t)+\la f{''}(t)+\la (\eta- \be) f'(t)-ABf(t)+\eta Ag'(t)-\la^2 Ag(t)=0\,.
\end{equation}
Thanks to \eqref{eq:gprime0} and $Ag(0)=-f{''}(0)-\la f(0)$, we have
\begin{multline*}
f^{(4)}(0)=-\eta f{'''}(0)-\la f{''}(0)-\la (\eta- \be) f'(0)+ABf(0)-\eta Ag'(0)+\la^2 Ag(0)\\
=-\eta f{'''}(0)-\la f{''}(0)-\la (\eta- \be) f'(0)+ABf(0)+\eta f{'''}(0)\\
+\eta\la f{'}(0)-\eta\la \be f(0)-\la^2 f{''}(0)-\la^2\la f(0)
\\
=-(\la+\la^2) f{''}(0)+\la \be f'(0)-(\eta\la \be +\la^3 -AB)f(0)
\,,
\end{multline*}
so the equation for $f^{(4)}(0)$ in  (\ref{eq:fifth}) holds true.
By differentiating \eqref{eq:fourth} we obtain: $f\in C^5([0,\infty))$ and
\begin{equation*}
f^{(5)}(t)+\eta f^{(4)}(t)+\la f{'''}(t)+\la (\eta- \be) f''(t)-ABf'(t)+\eta Ag''(t)-\la^2 Ag'(t)=0\,.
\end{equation*}
Therefore, by using again
$g{''}(t)=-Bf(t)-\la^2 g(t)$
we get 
\begin{equation*}
f^{(5)}(t)+\eta f^{(4)}(t)+\la f{'''}(t)+\la (\eta- \be) f''(t)-ABf'(t)-\eta ABf(t)-\la^2 Ag'(t)-\eta\la^2  Ag(t)=0\,.
\end{equation*}
It follows from \eqref{eq:gprime} 
\begin{equation*}
-Ag'(t)-\eta Ag(t)=f{'''}(t)+\eta f{''}(t)+\la f{'}(t)+\la (\eta- \be) f(t)\,,
\end{equation*}
and hence we have
\begin{multline*}
f^{(5)}(t)+\eta f^{(4)}(t)+(\la+\la^2) f{'''}(t)+(\la (\eta- \be)+\la^2\eta) f''(t)
\\
+(\la^3-AB)f'(t)+(\la^3 (\eta- \be)-\eta AB)f(t)
=0\,,
\end{multline*}
that is
$f$ verifies the differential equation in (\ref{eq:fifth}). 
Finally,
from the first equation in (\ref{eq:system0}) we deduce that
$g$ is given by \eqref{eq:Ag}.

Conversely, if $f$ is a solution of $(\ref{eq:fifth})$,  
multiplying the differential equation
 by $e^{\eta t}$ and integrating from $0$ to $t$, we obtain 
\begin{multline*}
\int_0^t e^{\eta s}f^{(5)}(s)\ ds+
\eta \int_0^t e^{\eta s}f^{(4)}(s)\ ds
+(\la+\la^2)\int_0^t e^{\eta s}f{'''}(s)\ ds
 +\eta(\la+\la^2)\int_0^t e^{\eta s}f{''}(s)\ ds 
 \\
-\la\be\int_0^t e^{\eta s}f{''}(s)\ ds
+(\la^3-AB)\int_0^t e^{\eta s}f{'}(s)\ ds
+(\la^3(\eta-\be)-\eta AB)\int_0^t e^{\eta s}f(s)\ ds=0\,.
\end{multline*}
Integrating by parts the first, the third, the fifth and the sixth integral, we have
\begin{multline*}
e^{\eta t}f^{(4)}(t)-f^{(4)}(0)
+(\la+\la^2) e^{\eta t}f{''}(t) -(\la+\la^2) f{''}(0)
-\la\be e^{\eta t}f{'}(t)+\la\be f{'}(0)+\eta\la\be e^{\eta t}f{}(t) \\
-\eta\la\be  f{}(0)-\eta^2\la\be\int_0^t e^{\eta s}f{}(s)\ ds
+(\la^3-AB)e^{\eta t}f(t)-(\la^3-AB)f(0)
-\la^3\be\int_0^t e^{\eta s}f(s)\ ds=0\,.
\end{multline*}
Using the identity for $f^{(4)}(0)$ in  (\ref{eq:fifth}) and multiplying by $e^{-\eta t}$, we obtain
\begin{multline}\label{eq:fourthbis}
f^{(4)}(t)
+(\la+\la^2) f{''}(t) 
-\la\be f{'}(t)+\eta\la\be f{}(t) 
-\eta^2\la\be\int_0^t e^{-\eta(t- s)}f{}(s)\ ds\\
+(\la^3-AB)f(t)
-\la^3\be\int_0^t e^{-\eta(t- s)}f(s)\ ds=0\,.
\end{multline}
Moreover, by \eqref{eq:Ag} it follows
\begin{equation*}
Ag'(t)=-f{'''}(t)-\la f'(t)+\la\beta f(t)-\eta\la\be \int_0^t e^{-\eta(t-s)}f(s) ds\,,
\end{equation*}
and hence
\begin{equation*}
Ag''(t)=-f^{(4)}(t)-\la f''(t)+\la\beta f'(t)-\eta\la\be f(t) +\eta^2\la\be\int_0^t e^{-\eta(t-s)}f(s) ds\,.
\end{equation*}
Therefore, thanks to the previous identity and \eqref{eq:fourthbis} we have
\begin{equation*}
Ag''(t)=
\la^2 f{''}(t) 
+(\la^3-AB)f(t)
-\la^3\be\int_0^t e^{-\eta(t- s)}f(s)\ ds\,,
\end{equation*}
whence, in view of \eqref{eq:Ag} we get
\begin{equation*}
Ag''(t)=-\la^2 Ag(t)-ABf(t)\,.
\end{equation*}
Finally, by \eqref{eq:Ag} and the above equation, it follows that the couple $(f, g)$ is a solution of the system \eqref{eq:system0}.
\end{Proof}

\section{The Hilbert Uniqueness Method}\label{se:HUM}
For reader's convenience,  in this section we will describe the Hilbert Uniqueness Method for coupled systems.
For another  approach based
on the ontoness  of the solution operator, see e.g. \cite{LasT, T1}.

Given $k\in L_{loc}^1(0,\infty)$ and $A\,,B\in\R$,
we consider the following coupled system:
\begin{equation}\label{eq:problem-u}
\begin{cases}
\displaystyle 
u_{1tt}(t,x) -u_{1xx}(t,x)+\int_0^t\ k(t-s) u_{1xx}(s,x)ds+Au_2(t,x)= 0\,,
\\
\phantom{u_{1tt}(t,x) -u_{1xx}(t,x)+\int_0^t\ k(t-s) u_{1xx}(s,x)ds+Au_2(t,x)= 0\,,\qquad}
t\in (0,T)\,,\quad x\in(0,\pi)
\\
\displaystyle
u_{2tt}(t,x) +u_{2xxxx}(t,x)+Bu_1(t,x)= 0
\,,
\end{cases}
\end{equation}
with null initial conditions 
\begin{equation}
u_1(0,x)=u_{1t}(0,x)=u_2(0,x)=u_{2t}(0,x)=0\qquad  x\in(0,\pi)\,,
\end{equation} 
and boundary conditions
\begin{equation}\label{eq:bound-u1}
u_1(t,0)=0\,,\qquad u_1(t,\pi)=g_1(t)\qquad t\in (0,T) \,,
\end{equation}
\begin{equation}\label{eq:bound-u2}
u_2(t,0)=u_{2xx}(t,0)=0\,,\qquad u_2(t,\pi)=g_2(t)\,,\qquad u_{2xx}(t,\pi)=g_3(t)\qquad t\in (0,T) \,.
\end{equation}
For a reachability problem we mean the following: given $T>0$ and taking
$
(u_{10},u_{11},u_{20},u_{21})
$
in a suitable space to define later,
find $g_i\in L^2(0,T)$, $i=1,2,3$ such that the weak
solution $u$ of problem  \eqref{eq:problem-u}-\eqref{eq:bound-u2}
verifies the final conditions
\begin{equation}\label{eq:problem-u1}
u_1(T,x)=u_{10}(x)\,,\quad u_{1t}(T,x)=u_{11}(x)\,,\quad
u_2(T,x)=u_{20}(x)\,,\quad u_{2t}(T,x)=u_{21}(x)\,,
\quad x\in(0,\pi)\,.
\end{equation}
One can solve such reachability problems by the HUM method. To see that, we proceed as follows.

Given
$(z_{10},z_{11},z_{20},z_{21})\in (C^\infty_c(0,\pi))^4$,
we introduce the {\it adjoint} system of (\ref{eq:problem-u}), that is 
\begin{equation}\label{eq:adjoint}
\begin{cases}
\displaystyle 
z_{1tt}(t,x) -z_{1xx}(t,x)+\int_t^T\ k(s-t) z_{1xx}(s,x)ds+Bz_2(t,x)= 0\,,\\
\phantom{u_{1tt}(t,x) -u_{1xx}(t,x)+\int_0^t\ k(t-s) u_{1xx}(s,x)ds+Au_2(t,x)= 0\,,\qquad}
t\in (0,T)\,,\quad x\in(0,\pi)
\\
\displaystyle
z_{2tt}(t,x) +z_{2xxxx}(t,x)+Az_1(t,x)= 0\,,
\\
z_1(t,0)=z_1(t,\pi)=z_2(t,0)=z_2(t,\pi)=z_{2xx}(t,0)=z_{2xx}(t,\pi)=0\qquad t\in [0,T]\,,
\end{cases}
\end{equation}
with  final data 
\begin{equation} \label{eq:final}
z_1(T,\cdot)=z_{10}\,,\quad z_{1t}(T,\cdot)=z_{11}\,,\quad 
z_2(T,\cdot)=z_{20}\,,\quad z_{2t}(T,\cdot)=z_{21}\,.
\end{equation} 
The above problem is well-posed, see e.g. \cite{Pruss}. 
Thanks to the regularity of the final data, the solution  $(z_1,z_2)$ of \eqref{eq:adjoint}--\eqref{eq:final}
is regular enough to 
consider  the  nonhomogeneous problem
\begin{equation}\label{eq:phi}
\left \{\begin{array}{l}\displaystyle
\phi_{1tt}(t,x) -\phi_{1xx}(t,x)+\int_0^t\ k(t-s) \phi_{1xx}(s,x)ds+A \phi_ 2(t,x)= 0
\\
\phantom{\phi_{1tt}(t,x) -\phi_{1xx}(t,x)+\int_0^t\ k(t-s) \phi_{1xx}(s,x)ds+A \phi_ 2(t,x)= 0\,,\qquad}
t\in (0,T)\,,\quad x\in(0,\pi)\,,
\\
\displaystyle
\phi_{2tt}(t,x) + \phi_{2xxxx}(t,x)+B \phi_ 1(t,x)= 0
\\
\\
\phi_1(0,x)= \phi_{1t}(0,x)= \phi_ 2(0,x)= \phi_{2t}(0,x)=0\qquad  x\in(0,\pi)\,,
\\
\\
\displaystyle
\phi_1(t,0)=0\,,\quad
\phi_1(t,\pi)=z_{1x}(t,\pi)-\int_t^T\ k(s-t)z_{1x}(s,\pi)ds\qquad t\in [0,T]\,,
\\
\\
\phi_2(t,0)=\phi_{2xx}(t,0)=0\,,\quad 
\phi_2(t,\pi)=-z_{2xxx}(t,\pi)\,,\quad
\phi_{2xx}(t,\pi)=-z_{2x}(t,\pi)\qquad t\in [0,T]\,.
\end{array}\right .
\end{equation} 

As in the non integral case, it can be proved  that  problem \eqref{eq:phi} admits a unique solution
$\phi$.
So, we can introduce the following linear operator: for any  $(z_{10},z_{11},z_{20},z_{21})\in \big(C^\infty_c(0,\pi)\big)^4$
we define
\begin{equation}\label{eq:psi0}
\Psi(z_{10},z_{11},z_{20},z_{21})=(-\phi_{1t}(T,\cdot),\phi_{1}(T,\cdot),-\phi_{2t}(T,\cdot),\phi_{2}(T,\cdot))
\,.
\end{equation}
For any $(\xi_{10},\xi_{11},\xi_{20},\xi_{21})\in \big(C^\infty_c(0,\pi)\big)^4$, let $(\xi_1,\xi_2)$ be the solution of 
\begin{equation}\label{eq:adjoint10}
\left \{\begin{array}{l}\displaystyle
\xi_{1tt}(t,x) -\xi_{1xx}(t,x)+\int_t^T\ k(s-t) \xi_{1xx}(s,x)ds+B\xi_2(t,x)= 0\,,\\
\phantom{\xi_{1tt}(t,x) -\xi_{1xx}(t,x)+\int_0^t\ k(t-s) \xi_{1xx}(s,x)ds+A\xi_2(t,x)= 0\,,\qquad}
t\in (0,T)\,,\quad x\in(0,\pi)
\\
\displaystyle
\xi_{2tt}(t,x) +\xi_{2xxxx}(t,x)+A\xi_1(t,x)= 0\,,
\\
\\
\xi_1(t,0)=\xi_1(t,\pi)=\xi_2(t,0)=\xi_2(t,\pi)=\xi_{2xx}(t,0)=\xi_{2xx}(t,\pi)=0\qquad t\in [0,T]\,,
\\
\\
\xi_1(T,\cdot)=\xi_{10}\,,\quad \xi_{1t}(T,\cdot)=\xi_{11}\,,\quad 
\xi_2(T,\cdot)=\xi_{20}\,,\quad \xi_{2t}(T,\cdot)=\xi_{21}\,.
\end{array}\right .
\end{equation} 
We will prove that
\begin{multline}\label{eq:psi}
\langle\Psi(z_{10},z_{11},z_{20},z_{21}),(\xi_{10},\xi_{11},\xi_{20},\xi_{21})\rangle_{L^2(0,\pi)}
\\
=\int_0^T\phi_1(t,\pi)\Big(\xi_{1x}(t,\pi)-\int_t^T\ k(s-t)\ \xi_{1x}(s,\pi)\ ds\Big) \ dt
\\
 -\int_0^T\phi_{2xx}(t,\pi)\xi_{2x}(t,\pi)\ dt-\int_0^T\phi_{2}(t,\pi)\xi_{2xxx}(t,\pi)\ dt
 \,. 
\end{multline}
To this end, we multiply the first equation in (\ref{eq:phi}) by $\xi_1$ and integrate 
on $[0,T]\times[0,\pi]$, so we have
\begin{multline*}
\int_0^\pi \int_0^T\phi_{1tt}(t,x)\xi_1(t,x)\ dt \ dx
-\int_0^T\int_0^\pi\phi_{1xx}(t,x)\xi_1(t,x)\ dx\ dt
\\
+\int_0^\pi\int_0^T\int_0^t\ k(t-s)\phi_{1xx}(s,x)\ ds\ \xi_1(t,x)\ dt\ dx
+A\int_0^T\int_0^\pi \phi_{2}(t,x)\xi_1(t,x)\ dx \ dt=0\,. 
\end{multline*}
If we take into account that 
\begin{equation*}
\int_0^T\int_0^t\ k(t-s)\phi_{1xx}(s,x)\ ds\ \xi_1(t,x)\ dt  =
\int_0^T\phi_{1xx}(s,x)\int_s^T\ k(t-s)\ \xi_1(t,x)\ dt \ ds 
\end{equation*}
and
integrate by parts,  then we have
\begin{multline*}
 \int_0^\pi\big(\phi_{1t}(T,x)\xi_{10}(x)- \phi_1(T,x)\xi_{11}(x)\big)\ dx
 +\int_0^\pi \int_0^T\phi_1(t,x)\xi_{1tt}(t,x)\ dt \ dx
 \\
+\int_0^T\phi_1(t,\pi)\xi_{1x}(t,\pi)\ dt -\int_0^T\int_0^\pi\phi_1(t,x)\xi_{1xx}(t,x)\ dx\ dt
 \\
-\int_0^T\phi_{1}(s,\pi)\int_s^T\ k(t-s)\ \xi_{1x}(t,\pi)\ dt \ ds 
+\int_0^\pi\int_0^T\phi_1(s,x)\int_s^T\ k(t-s)\ \xi_{1xx}(t,x)\ dt \ ds \ dx
\\
+A\int_0^T\int_0^\pi \phi_{2}(t,x)\xi_1(t,x)\ dx \ dt=0\,. 
\end{multline*}
As a consequence of the above equation and
\begin{equation*}
\xi_{1tt} -\xi_{1xx}+\int_t^T\ k(s-t) \xi_{1xx}(s,\cdot)ds=-B\xi_2\,,
\end{equation*}
we obtain
\begin{multline}\label{eq:xi1}
 \int_0^\pi\big(\phi_{1t}(T,x)\xi_{10}(x)- \phi_1(T,x)\xi_{11}(x)\big)\ dx
+\int_0^T\phi_1(t,\pi)\Big(\xi_{1x}(t,\pi)-\int_t^T\ k(s-t)\ \xi_{1x}(s,\pi)\ ds\Big) \ dt 
\\
+\int_0^T\int_0^\pi \big(A\phi_{2}(t,x)\xi_1(t,x)-B\phi_1(t,x)\xi_{2}(t,x)\big)\ dx \ dt=0\,. 
\end{multline}
In a similar way, we multiply the second equation in (\ref{eq:phi}) by $\xi_2$ and integrate by parts
on $[0,T]\times[0,\pi]$ to get
\begin{multline*}
 \int_0^\pi\big(\phi_{2t}(T,x)\xi_{20}(x)- \phi_2(T,x)\xi_{21}(x)\big)\ dx
 +\int_0^\pi \int_0^T\phi_2(t,x)\xi_{2tt}(t,x)\ dt \ dx
 \\
-\int_0^T\big(\phi_{2xx}(t,\pi)\xi_{2x}(t,\pi)+\phi_{2}(t,\pi)\xi_{2xxx}(t,\pi)\big)\ dt 
+\int_0^T\int_0^\pi\phi_2(t,x)\xi_{2xxxx}(t,x)\ dx\ dt
\\
+B\int_0^T\int_0^\pi \phi_{1}(t,x)\xi_2(t,x)\ dx \ dt=0\,, 
\end{multline*}
whence, in virtue of
\begin{equation*}
\xi_{2tt} +\xi_{2xxxx}=-A\xi_1\,,
\end{equation*}
we get
\begin{multline}\label{eq:xi2}
 \int_0^\pi\big(\phi_{2t}(T,x)\xi_{20}(x)- \phi_2(T,x)\xi_{21}(x)\big)\ dx
-\int_0^T\big(\phi_{2xx}(t,\pi)\xi_{2x}(t,\pi)+\phi_{2}(t,\pi)\xi_{2xxx}(t,\pi)\big)\ dt 
\\
+\int_0^T\int_0^\pi\big(B \phi_{1}(t,x)\xi_2(t,x)-A\phi_{2}(t,x)\xi_1(t,x)\big)\ dx \ dt=0\,. 
\end{multline}
If we sum equations \eqref{eq:xi1} and \eqref{eq:xi2}, then we have
\begin{multline}\label{eq:prenorm}
\langle\Psi(z_{10},z_{11},z_{20},z_{21}),(\xi_{10},\xi_{11},\xi_{20},\xi_{21})\rangle_{L^2(0,\pi)}
\\
= \int_0^\pi\big(-\phi_{1t}(T,x)\xi_{10}(x)+ \phi_1(T,x)\xi_{11}(x)-\phi_{2t}(T,x)\xi_{10}(x)+ \phi_2(T,x)\xi_{11}(x)\big)\ dx
 \\
 =\int_0^T\phi_1(t,\pi)\Big(\xi_{1x}(t,\pi)-\int_t^T\ k(s-t)\ \xi_{1x}(s,\pi)\ ds\Big) \ dt
 \\
 -\int_0^T\big(\phi_{2xx}(t,\pi)\xi_{2x}(t,\pi)+\phi_{2}(t,\pi)\xi_{2xxx}(t,\pi)\big)\ dt \,,
 \end{multline}
that is, \eqref{eq:psi} holds true.

Now, taking $(\xi_{10},\xi_{11},\xi_{20},\xi_{21})=(z_{10},z_{11},z_{20},z_{21})$ in  (\ref{eq:psi}), we have
\begin{multline}\label{eq:psi1}
\langle\Psi(z_{10},z_{11},z_{20},z_{21}),(z_{10},z_{11},z_{20},z_{21})\rangle_{L^2(0,\pi)}
\\
=
\int_0^T\Big|z_{1x}(t,\pi)-\int_t^T\ k(s-t)\ z_{1x}(s,\pi)\ ds\Big|^2 \ dt
 +\int_0^T\big|z_{2x}(t,\pi)\big|^2\ dt+\int_0^T\big|z_{2xxx}(t,\pi)\big|^2\ dt\,. 
\end{multline}
As a consequence, we can introduce a semi-norm on the space $\big(C^\infty_c(\Omega)\big)^4$. Precisely,
we define, for $(z_{10},z_{11},z_{20},z_{21})\in \big(C^\infty_c(\Omega)\big)^4$,
\begin{multline}\label{eq:normF}
\|(z_{10},z_{11},z_{20},z_{21})\|_{F}:=
\\
\displaystyle\Big(
\int_0^T\Big|z_{1x}(t,\pi)-\int_t^T\ k(s-t)\ z_{1x}(s,\pi)\ ds\Big|^2 \ dt
 +\int_0^T\big|z_{2x}(t,\pi)\big|^2\ dt+\int_0^T\big|z_{2xxx}(t,\pi)\big|^2\ dt
\Big)^{1/2}\,.
\end{multline}
If $k(t)=\beta e^{-\eta t}$,
thanks to
\eqref{eq:unicita}, 
$\|\cdot\|_{F}$ is a norm if and only if the following uniqueness theorem holds.
\begin{theorem}\label{th:uniqueness}
If $(z_1,z_2)$ is the solution of problem {\rm (\ref{eq:adjoint})--(\ref{eq:final})} such that
$$
z_{1x}(t,\pi)=z_{2x}(t,\pi)=z_{2xxx}(t,\pi)=0\,,\qquad \forall t\in [0,T]\,,
$$
then 
$$
z_1(t,x)=z_2(t,x)= 0 \qquad\forall (t,x)\in [0,T]\times[0,\pi]\,.
$$
\end{theorem}
If theorem \ref{th:uniqueness} holds true, then we can define the Hilbert space $F$ as the completion of $ \big(C^\infty_c(\Omega)\big)^4$ for
the norm (\ref{eq:normF}). Moreover, the operator $\Psi$ extends uniquely to a continuous operator, denoted again by $\Psi$, from $F$ to the dual
space $F'$ in such a way that   $\Psi:F\to F'$ is an isomorphism.

In conclusion,
if we prove  the uniqueness result given by theorem \ref{th:uniqueness} and 
$$F=H^1_0(0,\pi)\times L^2(0,\pi)\times H^3_0(0,\pi)\times H^{1}(0,\pi)$$
with the equivalence of the respective norms, then we can solve the reachability problem
\eqref{eq:problem-u}--\eqref{eq:problem-u1}
taking $(u_{10},u_{11},u_{20},u_{21})\in  L^{2}(0,\pi)\times H^{-1}(0,\pi)\times H^{-1}(0,\pi)\times H^{-3}(0,\pi)$.

In addition, if $g_2(t)\equiv0$ in the reachability problem
\eqref{eq:problem-u}--\eqref{eq:problem-u1}, we must take $\phi_2(t,\pi)\equiv0$ in problem \eqref{eq:phi}. So, in view of \eqref{eq:prenorm}
formula \eqref{eq:normF} becomes
\begin{multline}\label{eq:normF1}
\|(z_{10},z_{11},z_{20},z_{21})\|_{F}:=
\displaystyle\Big(
\int_0^T\Big|z_{1x}(t,\pi)-\int_t^T\ k(s-t)\ z_{1x}(s,\pi)\ ds\Big|^2 \ dt
 +\int_0^T\big|z_{2x}(t,\pi)\big|^2\ dt\,,
\end{multline}
and the uniqueness result can be written in this way
\begin{theorem}\label{th:uniqueness1}
If $(z_1,z_2)$ is the solution of problem {\rm (\ref{eq:adjoint})--(\ref{eq:final})} such that
$$
z_{1x}(t,\pi)=z_{2x}(t,\pi)=0\,,\qquad \forall t\in [0,T]\,,
$$
then 
$$
z_1(t,x)=z_2(t,x)= 0 \qquad\forall (t,x)\in [0,T]\times[0,\pi]\,.
$$
\end{theorem}
Finally, by proving  the uniqueness result given by Theorem \ref{th:uniqueness1} and 
$$F=H^1_0(0,\pi)\times L^2(0,\pi)\times H^1_0(0,\pi)\times H^{-1}(0,\pi)$$
with the equivalence of the respective norms, we can solve the reachability problem
\eqref{eq:problem-u}--\eqref{eq:problem-u1} for
$(u_{10},u_{11},u_{20},u_{21})\in  L^{2}(0,\pi)\times H^{-1}(0,\pi)\times H_0^{1}(0,\pi)\times H^{-1}(0,\pi)$.

\section{Spectral analysis 
}\label{se:specana}
In this section we will elaborate a detailed spectral analysis for the adjoint problem.

Let 
$L:D(L)\subset X\to X$ be a self-adjoint positive linear
 operator on $X$ with dense domain $D(L)$ and let $\{\la_j\}_{j\ge1}$ be a strictly increasing sequence  of  eigenvalues for the operator $L$ with
$\la_j>0$ and $\la_j\to\infty$ such that the sequence of the corresponding eigenvectors $\{w_j\}_{j\ge1}$ constitutes a Hilbert basis for $X$.

Fix two real numbers $A$\,,$B$
and consider the following weakly coupled system:
\begin{equation}\label{eq:system}
\begin{cases}
\displaystyle 
u_1''(t) +Lu_1(t)-\beta\int_0^t\ e^{-\eta(t-s)}L u_1(s)ds+Au_2(t)= 0\,,\qquad
t\ge 0
\\
\displaystyle
u_2''(t) +L^2u_2(t)+Bu_1(t)= 0
\,,\qquad t\ge 0\,,
\\
u_1(0)=u_{10}\,,\quad u'_1(0)=u_{11}\,,
\\
u_2(0)=u_{20}\,,\quad u'_2(0)=u_{21}\,.
\end{cases}
\end{equation}

We have
\begin{equation*}
u_{10}=\sum_{j=1}^{\infty}\alpha_{1j}w_{j}\,,\qquad\quad\alpha_{1j}=\langle u_{10},w_j\rangle \,,
\quad\sum_{j=1}^{\infty}\alpha_{1j}^2\lambda_j<\infty\,,
\end{equation*} 
\begin{equation*}
u_{11}=\sum_{j=1}^{\infty}\rho_{1j}w_{j}\,,\qquad\quad\rho_{1j}=\langle u_{11},w_j\rangle \,,\quad\sum_{j=1}^{\infty}\rho_{1j}^2<\infty\,,
\end{equation*} 
\begin{equation*}
u_{20}=\sum_{j=1}^{\infty}\alpha_{2j}w_{j}\,,\qquad\quad\alpha_{2j}=\langle u_{20},w_j\rangle \,,
\quad\sum_{j=1}^{\infty}\alpha_{2j}^2\lambda_j<\infty\,,
\end{equation*} 
\begin{equation*}
u_{21}=\sum_{j=1}^{\infty}\rho_{2j}w_{j}\,,\qquad\quad\rho_{2j}=\langle u_{21},w_j\rangle \,,
\quad\sum_{j=1}^{\infty}\frac{\rho_{2j}^2}{\lambda_j}<\infty\,.
\end{equation*} 

We will seek the solution $(u_{1}(t),u_{2}(t))$  of system \eqref{eq:system} with components written as sums of series, that is
\begin{equation*}
u_{1}(t)=\sum_{j=1}^{\infty}f_{1j}(t)w_{j}\,,\quad u_{2}(t)=\sum_{j=1}^{\infty}f_{2j}(t)w_{j}\,,
\qquad
f_{ij}(t)=\langle u_{i}(t),w_j\rangle\,,\quad i=1,2
\,.
\end{equation*}
If we put the above  expressions for $u_{1}(t)$ and $u_{2}(t)$
into \eqref{eq:system} and multiply  by $w_j$, $j\in\N$, then  we have that $(f_{1j}(t),f_{2j}(t))$ is the solution of system
\begin{equation}\label{eq:secondsys}
\begin{cases}\displaystyle
f_{1j}^{''}(t)
+\la_{j}f_{1j}(t)-\la_{j}\be \int_0^t e^{-\eta(t-s)}f_{1j}(s) ds
+Af_{2j}(t)=0
\\\displaystyle
f_{2j}^{''}(t)+\la_j^2 f_{2j}(t)+Bf_{1j}(t)=0
\,,
\\
f_{1j}(0)=\a_{1j}\,, \quad f_{1j}^{'}(0)=\rho_{1j}\,,
\\
f_{2j}(0)=\a_{2j}\,, \quad f_{2j}^{'}(0)=\rho_{2j}\,.
\end{cases}
\end{equation}
Thanks to lemma \ref{le:fifth}, $(f_{1j}(t),f_{2j}(t))$ is the solution of problem \eqref{eq:secondsys} if and only if  $f_{1j}(t)$ is the solution of the Cauchy problem
\begin{equation}\label{eq:third}
\begin{cases}
\displaystyle 
f_{1j}^{(5)}(t)+\eta f_{1j}^{(4)}(t)+(\la_{j}^2+\la_{j}) f_{1j}'''(t)+(\eta\la_{j}^2+\la_{j} (\eta-\be))f_{1j}''(t)+(\la_{j}^3-AB)f_{1j}'(t)
\\\displaystyle\hskip5cm
+(\la_{j}^3(\eta-\be)-\eta AB)f_{1j}(t)=0
\,,\qquad t\ge 0\,,
\\
 f_{1j}(0)=\alpha_{1j},
\quad
 f_{1j}'(0)=\rho_{1j},
 \\
 f_{1j}''(0)=-\la_{j}\alpha_{1j}-A\alpha_{2j},
 \quad
 f_{1j}'''(0)=-A\rho_{2j}-\la_{j} \rho_{1j}+\la_{j} \beta\alpha_{1j},
 \\
 f_{1j}^{(4)}(0)=(\la_{j}^2+\la_{j})(\la_{j}\alpha_{1j}+A\alpha_{2j})+\la_{j}\be \rho_{1j}
 -\la_{j}\eta\be \alpha_{1j}-(\la_{j}^3-AB)\alpha_{1j}\,,
\end{cases}
\end{equation} 
and $f_{2j}(t)$ is given by
\begin{equation}\label{eq:f2j}
f_{2j}(t)=-\frac1A\Big[f_{1j}^{''}(t)
+\la_{j}f_{1j}(t)-\la_{j}\be \int_0^t e^{-\eta(t-s)}f_{1j}(s) ds\Big]
\,.
\end{equation}
We proceed to solve  $(\ref{eq:third})$.
To this end,
we have to evaluate the solutions of the characteristic equation of the fifth degree
\begin{equation}\label{eq:fchar}
\La^{5}+\eta\La^{4}+(\la_{j}^2+\la_{j})\La^{3}+(\eta\la_{j}^2+\la_{j} (\eta-\be))\La^{2}+(\la_{j}^3-AB)\La+\la_{j}^3(\eta-\be)-\eta AB=0\,.
\end{equation}
The asymptotic behavior  of the solutions of  equation (\ref{eq:fchar}) as $j\to\infty$ is the following
\begin{equation}\label{eq:lambda1}
\La_{1j}=\be-\eta
-{\beta\big(\be-\eta\big)^2\over\la_{j}}+O\Big({1\over{\la_{j}^{2}}}\Big)
=\be-\eta
+O\Big({1\over{\la_{j}}}\Big)
\,,
\end{equation}
\begin{multline}\label{eq:lambda2}
\La_{2j}=- {\be\over 2}
+{\beta\big(\be-\eta\big)^2\over2}{1\over\la_{j}}
+
i\Big[\sqrt{\la_{j}}+{\be\over2}\Big({3\over4}\beta-\eta\Big){1\over\sqrt\lambda_{j}}\Big] 
+O\Big({1\over{\la_{j}^{3/2}}}\Big)
\\
=- {\be\over 2}
+i\Big[\sqrt{\la_{j}}+{\be\over2}\Big({3\over4}\beta-\eta
\Big){1\over\sqrt\lambda_{j}}\Big]
+O\Big({1\over{\la_{j}}}\Big)
\,,
\end{multline}
\begin{equation}\label{eq:lambda3}
\La_{3j}=\overline{\La_{2j}}
=- {\be\over 2}
-i\Big[\sqrt{\la_{j}}+{\be\over2}\Big({3\over4}\beta-\eta
\Big){1\over\sqrt\lambda_{j}}\Big]
+O\Big({1\over{\la_{j}}}\Big)
\,,
\end{equation}
\begin{equation}\label{eq:lambda4}
\La_{4j}
=-{\be AB\over 2\la_{j}^5}
+i\Big(\la_{j}+{AB\over 2\la_{j}^3}+{AB\over 2\la_{j}^4}+{AB\over 2\la_{j}^5}\Big)
+O\Big({1\over{\la_{j}^{6}}}\Big)
=i\la_{j}+O\Big({1\over{\la_{j}^{3}}}\Big)
\,,
\end{equation}
\begin{equation}\label{eq:lambda5}
\La_{5j}=\overline{\La_{4j}}
=-i\la_{j}+O\Big({1\over{\la_{j}^{3}}}\Big)
\,.
\end{equation}
Therefore, we can write the solution of (\ref{eq:third}) in
the form
\begin{equation}\label{eq:f1j}
f_{1j}(t)=C_{1j}e^{t\La_{1j}}+C_{2j}e^{t\La_{2j}}+C_{3j}e^{t\La_{3j}}+C_{4j}e^{t\La_{4j}}+C_{5j}e^{t\La_{5j}}
=\sum_{k=1}^{5}C_{kj}e^{t\La_{kj}}\,,
\end{equation}
where $C_{kj}$ are complex numbers.
To determine the coefficients $C_{kj}$, we have to impose the initial conditions in $(\ref{eq:third})$, that is we must solve 
 the system
\begin{equation}\label{vandermonde}
\left \{\begin{array}{l}
\phantom{\La_{1j}}C_{1j}+\phantom{\La_{1j}}C_{2j}+\phantom{\La_{1j}}C_{3j}+\phantom{\La_{1j}}C_{4j}
+\phantom{\La_{1j}}C_{5j}=f_{1j}(0)\\ 
\\
\La_{1j}C_{1j}+\La_{2j}C_{2j}+\La_{3j}C_{3j}+\La_{4j}C_{4j}+\La_{5j}C_{5j}
=f_{1j}^{'}(0)\\
\\
\La_{1j}^2C_{1j}+\La_{2j}^2C_{2j}+\La_{3j}^2C_{3j}+\La_{4j}^2C_{4j}+\La_{5j}^2C_{5j}=f_{1j}^{''}(0)\\
\\
\La_{1j}^3C_{1j}+\La_{2j}^3C_{2j}+\La_{3j}^3C_{3j}+\La_{4j}^3C_{4j}+\La_{5j}^3C_{5j}=f_{1j}^{'''}(0)\\
\\
\La_{1j}^4C_{1j}+\La_{2j}^4C_{2j}+\La_{3j}^4C_{3j}+\La_{4j}^4C_{4j}+\La_{5j}^4C_{5j}=f_{1j}^{(4)}(0)\,.
\end{array}\right .
\end{equation}
Therefore, we have the following asymptotic behavior as $j\to\infty$ of the coefficients $C_{kj}$:
\begin{equation}\label{eq:Cj}
\begin{cases}
C_{1j}={\beta\over\la_j}(\rho_{1j}+\alpha_{1j}(\beta-\eta))
+( \alpha_{1j}+\rho_{1j})O\Big({1\over{\la_{j}^{2}}}\Big)
\\
C_{2j}={\alpha_{1j}\over 2}
-\frac i{4\la_{j}^{1/2}}(\beta\alpha_{1j}+2\rho_{1j})
+( \alpha_{1j}+\rho_{1j})O\Big({1\over{\la_{j}}}\Big)
\\
C_{3j}=\overline{C_{2j}}
\\
C_{4j}=\frac{A\alpha_{2j}}{2\la^2_{j}}
+(\alpha_{2j}-i\rho_{2j})\frac {A}{2\la_{j}^{3}}
+( \alpha_{2j}+\rho_{2j})O\Big({1\over{\la^{7/2}_{j}}}\Big)
\\
C_{5j}=\overline{C_{4j}}\,.
\end{cases}
\end{equation}
Thanks to the expressions of $C_{kj}$ and $\La_{kj}$, $k=1,2,3$, we note that the function
\begin{equation}\label{eq:f1j3}
f^{*}_{1j}(t)=C_{1j}e^{t\La_{1j}}+C_{2j}e^{t\La_{2j}}+\overline{C_{2j}}e^{t\overline{\La_{2j}}}\,,
\end{equation}
verifies the problem
\begin{equation*}
\left \{\begin{array}{l}\displaystyle
(f^{*}_{1j})^{''}(t)
+\la_{j}f^{*}_{1j}(t)-\la_{j}\be \int_0^t e^{-\eta(t-s)}f^{*}_{1j}(s) ds=0\,,\\
\\
f^{*}_{1j}(0)=\a_{1j}\,, \qquad (f^{*}_{1j})^{'}(0)=\rho_{1j}\,,
\end{array}\right .
\end{equation*} 
see \cite[Section 6]{LoretiSforza1}.
Therefore, in view of
 \eqref{eq:f2j} the coefficients $f_{2j}$ are given by
\begin{multline}\label{eq:f2jbis}
f_{2j}(t)=-\frac1A
\Big(C_{4j}\Big(\La_{4j}^2+\la_{j}-\frac{\beta\la_{j}}{\eta+\La_{4j}}\Big)e^{t\La_{4j}}
+\beta e^{-\eta t}\la_{j}
\frac{C_{4j}}{\eta+\La_{4j}}
\Big)
\\
-\frac1A
\Big(\overline{C_{4j}}\Big(\overline{\La_{4j}}^2+\la_{j}-\frac{\beta\la_{j}}{\eta+\overline{\La_{4j}}}\Big)e^{t\overline{\La_{4j}}}
+\beta e^{-\eta t}\la_{j}
\frac{\overline{C_{4j}}}{\eta+\overline{\La_{4j}}}
\Big)
\qquad
t\ge 0\,.
\end{multline}

The proof of the following lemma is based on considerations similar to those used for analogous results in \cite{LoretiSforza}, but, for the sake of completeness, we prefer to give it.
\begin{lemma}
The following estimates hold true:
\begin{itemize}
\item[(i)] there exist some constants $c_1\,,c_2>0$ such that we have, for any $j\in\N$,
\begin{equation}\label{eq:|Cj2|}
\frac{c_1}{\la_{j}}\big(\alpha^2_{1j}\la_{j}+\rho^2_{1j}\big)
\le
|C_{2j}|^2
\le
\frac{c_2}{\la_{j}}\big(\alpha^2_{1j}\la_{j}+\rho^2_{1j}\big)\,;
\end{equation}
\item[(ii)] there exists a constant $c>0$ such that  we have, for any $j\in\N$,
\begin{equation}\label{eq:C1overC2}
{{\vert C_{1j}\vert}\over{\vert C_{2j}\vert}}\le {c \over {{\lambda^{1/2}_j}}}\,;
\end{equation}
\item[(iii)] there exist some constants $c_1\,,c_2>0$ such that we have, for any $j\in\N$,
\begin{equation}\label{eq:|Cj4|}
\frac{c_1}{\la^5_{j}}\Big(\alpha^2_{2j}\la_{j}+\frac{\rho^2_{2j}}{\la_{j}}\Big)
\le
|C_{4j}|^2
\le
\frac{c_2}{\la^5_{j}}\Big(\alpha^2_{2j}\la_{j}+\frac{\rho^2_{2j}}{\la_{j}}\Big)\,.
\end{equation}


\end{itemize}
\end{lemma}
\begin{Proof}
(i) First, we observe that
\begin{equation*}
|C_{2j}|^2={1\over 4}\Big(\alpha^2_{1j}+\frac{\rho^2_{1j}}{\la_{j}}\Big)
+\alpha^2_{1j}O\Big({1\over{\la_{j}}}\Big)
+\alpha_{1j}\rho_{1j}O\Big({1\over{\la_{j}}}\Big)
+\rho^2_{1j}O\Big({1\over{\la^2_{j}}}\Big)\,.
\end{equation*}
We can assume that for any $j\in\N$ $\alpha_{1j}\not=0$ or $\rho_{1j}\not=0$, and hence by the previous formula we obtain
\begin{equation*}
\frac{|C_{2j}|^2}{\alpha^2_{1j}+\frac{\rho^2_{1j}}{\la_{j}}}
={1\over 4}
+\frac{\big(\alpha^2_{1j}+\frac{\rho^2_{1j}}{\la_{j}}\big)O\Big({1\over{\la_{j}}}\Big)
+\alpha_{1j}\frac{\rho_{1j}}{\la_{j}^{1/2}}O\Big({1\over{\la_{j}^{1/2}}}\Big)
}
{\alpha^2_{1j}+\frac{\rho^2_{1j}}{\la_{j}}}
\to\frac14\,,\qquad\mbox{as}\qquad j\to\infty\,,
\end{equation*}
so, \eqref{eq:|Cj2|} follows.

\noindent (ii)
Since
\begin{equation*}
|C_{1j}|\le{|\alpha_{1j}|(\eta-\beta)+|\rho_{1j}|\over\la_j}
\left(\beta+
\frac{ |\alpha_{1j}|O\Big({1\over{\la_{j}}}\Big)+|\rho_{1j}|O\Big({1\over{\la_{j}}}\Big)}
{|\alpha_{1j}|(\eta-\beta)+|\rho_{1j}|}\right)\,,
\end{equation*}
we have, for any $j\in\N$,
\begin{equation}
|C_{1j}|\le c^*{|\alpha_{1j}|(\eta-\beta)+|\rho_{1j}|\over\la_j}
\,,
\end{equation}
for some $c^*>0$. Therefore, by using also \eqref{eq:|Cj2|} we get, for any $j\in\N$,
\begin{equation*}
{{\vert C_{1j}\vert}\over{\vert C_{2j}\vert}}\le {c^*\over {\sqrt{c_1\lambda_j}}}
\frac{|\alpha_{1j}|(\eta-\beta)+|\rho_{1j}|}{\sqrt{\alpha^2_{1j}\la_{j}+\rho^2_{1j}}}
\le {c\over {\lambda^{1/2}_j}}
\,,
\end{equation*}
so, we obtain \eqref{eq:C1overC2}.

\noindent (iii)
Notice that
\begin{equation*}
|C_{4j}|^2=\frac{A^2}{4}\Big(\frac{\alpha^2_{2j}}{\la^4_{j}}+\frac{\rho^2_{2j}}{\la^6_{j}}\Big)
+\alpha^2_{2j}o\Big({1\over{\la^4_{j}}}\Big)+\rho^2_{2j}o\Big({1\over{\la^6_{j}}}\Big)
\,,
\end{equation*}
and hence it follows
\begin{equation*}
\frac{|C_{4j}|^2}{\frac{\alpha^2_{2j}}{\la^4_{j}}+\frac{\rho^2_{2j}}{\la^6_{j}}}
={A^2\over 4}
+\frac{\alpha^2_{2j}o\Big({1\over{\la^4_{j}}}\Big)+\rho^2_{2j}o\Big({1\over{\la^6_{j}}}\Big)}
{\frac{\alpha^2_{2j}}{\la^4_{j}}+\frac{\rho^2_{2j}}{\la^6_{j}}}
\to{A^2\over 4}\,,\qquad\mbox{as}\qquad j\to\infty\,,
\end{equation*}
that is, \eqref{eq:|Cj4|} holds true.
\end{Proof}
In conclusion, keeping in mind \eqref{eq:f1j} and \eqref{eq:f2jbis},  the components $u_1(t)$ and  $u_2(t)$ of the
solution for the Cauchy problem \eqref{eq:system} are given by 
\begin{equation*}
u_1(t)=\sum_{j=1}^{\infty}\big(C_{1j}e^{t\La_{1j}}+C_{2j}e^{t\La_{2j}}+\overline{C_{2j}}e^{t\overline{\La_{2j}}}
+C_{4j}e^{t\La_{4j}}+\overline{C_{4j}}e^{t\overline{\La_{4j}}}\big)w_{j}
\,,
\end{equation*}
\begin{multline*}
u_2(t)=-\frac1A\sum_{j=1}^{\infty}
\Big(C_{4j}\Big(\La_{4j}^2+\la_{j}-\frac{\beta\la_{j}}{\eta+\La_{4j}}\Big)e^{t\La_{4j}}
+\overline{C_{4j}}\Big(\overline{\La_{4j}}^2+\la_{j}-\frac{\beta\la_{j}}{\eta+\overline{\La_{4j}}}\Big)e^{t\overline{\La_{4j}}}\Big)w_{j}
\\
-\frac{\beta }Ae^{-\eta t}\sum_{j=1}^{\infty}
\la_{j}\Big(\frac{C_{4j}}{\eta+\La_{4j}}
+\frac{\overline{C_{4j}}}{\eta+\overline{\La_{4j}}}
\Big)w_{j}
\,,
\end{multline*}
for any $t\ge 0$, where $\La_{kj}$ and $C_{kj}$ are defined by formulas \eqref{eq:lambda1}--\eqref{eq:lambda5} and \eqref{eq:Cj} respectively. 
We introduce, for any $n\ge1$, the following numbers
$r_{n}\,,R_{n}\in\R$ and
$\om_{n}\,,C_{n}\,,p_{n}\,,D_{n}\in\C$:
\begin{equation*}
r_n=\Lambda_{1n}=\be-\eta
+O\Big({1\over{\la_{n}}}\Big)
\,,
\end{equation*}
\begin{equation*}
\Re\omega_n=\Im\Lambda_{2n}=
\sqrt{\lambda_{n}}+{\beta\over2}\Big({3\over4}\beta-\eta
\Big){1\over\sqrt{\lambda_{n}}}+O\Big({1\over{\lambda_{n}}}\Big)
\,,
\end{equation*}
\begin{equation*}
\Im\omega_n=-\Re\Lambda_{2n}=
{\be\over 2}
+O\Big({1\over{\la_{n}}}\Big)
\,,
\end{equation*}
\begin{equation*}
\Re p_n=\Im\Lambda_{4n}=\la_{n}+O\Big({1\over{\la_{n}^{3}}}\Big)\,,
\end{equation*}
\begin{equation*}
\Im p_n=-\Re\Lambda_{4n}=O\Big({1\over{\la_{n}^{5}}}\Big)\,,
\end{equation*}
$$
R_n=C_{1n},
\qquad
C_n=C_{2n},
\qquad
D_n=C_{4n}\,.
$$
Thanks to these notations, the functions $u_1$ and $u_2$ can be written in the form  

%
\begin{equation}\label{eq:vsum}
u_1(t)=\sum_{n=1}^{\infty}\Big(R_ne^{r_nt}+C_ne^{i\om_nt}+\overline{C_n}e^{-i\overline{\omega_n}t}
+D_{n}e^{i p_{n}t}+\overline{D_{n}}e^{-i \overline{p_n}t}\Big)w_{n}\qquad
t\ge 0\,,
\end{equation}
\begin{equation}\label{eq:vsumbis}
u_2(t)=\sum_{n=1}^{\infty}
\Big(d_nD_{n}e^{i p_{n}t}
+\overline{d_n}\overline{D_{n}}e^{-i \overline{p_n}t}\Big)w_{n}
+\D e^{-\eta t}w_{n}
\quad
t\ge 0\,,
\end{equation}
where
\begin{equation}
d_n=\frac1A \Big(p^2_n-\Re p_n+\frac{\be\Re p_n}{\eta+ip_{n}}\Big)\,,
\end{equation}
\begin{equation}\label{eq:esprD}
\D=-\frac\beta A\sum_{n=1}^{\infty}\Re p_{n}
\Big(\frac{D_{n}}{\eta+ip_n}+\frac{\overline{D_{n}}}{\eta-ip_n}
\Big)\,.
\end{equation}
\begin{lemma}
There exist  constants $m_1,m_2>0$ such that
\begin{equation}\label{eq:alldn0}
m_1|p_n|^2\le|d_n|\le m_2|p_n|^2\,\qquad\forall n\in\N\,.
\end{equation}

\end{lemma}
\begin{Proof}
We note that
for $n_0$ sufficiently large we have, for any $n\ge n_0$,  
\begin{equation*}
|d_n|^2=
\Big|p^2_n-\Re p_n+\frac{\be\Re p_n}{\eta+ip_{n}}\Big|^2
=|p_n|^4\Big|1-\frac{\Re p_n}{p^2_n}+\frac{\be\Re p_n}{p^2_n(\eta+ip_{n})}\Big|^2
\le\frac{3}{2}|p_n|^4
\,,
\end{equation*}
and
\begin{equation*}
|d_n|^2=\Big|p^2_n-\Re p_n+\frac{\be\Re p_n}{\eta+ip_{n}}\Big|^2
=|p_n|^4\Big|1-\frac{\Re p_n}{p^2_n}+\frac{\be\Re p_n}{p^2_n(\eta+ip_{n})}\Big|^2
\ge\frac{|p_n|^4}{2}\,.
\end{equation*}
Since $ip_n$ is not a solution of the cubic equation
\begin{equation*}
\La^{3}+\eta\La^{2}+\Re p_n\La+\Re p_n(\eta-\be)=0\,,
\end{equation*}
we have for any $n\in\N$
\begin{equation*}
1-\frac{\Re p_n}{p^2_n}+\frac{\be\Re p_n}{p^2_n(\eta+ip_{n})}\not=0\,,
\end{equation*}
whence
\begin{equation*}
\min_{n\le n_0}\Big|1-\frac{\Re p_n}{p^2_n}+\frac{\be\Re p_n}{p^2_n(\eta+ip_{n})}\Big|>0\,,
\quad
\max_{n\le n_0}\Big|1-\frac{\Re p_n}{p^2_n}+\frac{\be\Re p_n}{p^2_n(\eta+ip_{n})}\Big|>0\,.
\end{equation*}
Therefore, there exist  constants $m_1,m_2>0$ such that \eqref{eq:alldn0} holds true.
\end{Proof}
\begin{remark} {\rm In the following section, we will skip the dependence on $w_{n}$ in \eqref{eq:vsum} and \eqref{eq:vsumbis}, because  that is not restricting,
as we will see in Theorem \ref{th:reachres}.}
\end{remark}

\section{Ingham type inequalities}\label{se:invdir}
In this section we will establish the inverse and direct inequalities for $(u_1,u_2)$, where
\begin{equation}\label{eq:vsum1}
u_1(t)=\sum_{n=1}^{\infty}\Big(R_ne^{r_nt}+C_ne^{i\om_nt}+\overline{C_n}e^{-i\overline{\omega_n}t}
+D_{n}e^{i p_{n}t}+\overline{D_{n}}e^{-i \overline{p_n}t}\Big)
\qquad
t\in\R 
\,,
\end{equation}
\begin{equation}
u_2(t)=\sum_{n=1}^{\infty}
\Big(d_nD_{n}e^{i p_{n}t}
+\overline{d_n}\overline{D_{n}}e^{-i \overline{p_n}t}\Big)
+\D e^{-\eta t}
\qquad
t\in\R 
\,,
\end{equation}
$r_{n}\,,R_{n}\,,\D\in\R$ and
$\om_{n}\,,C_{n}\,,p_{n}\,,D_{n}\in\C$, $p_n\not=0$,
by assuming that 
\begin{equation}\label{eq:hp1}
\lim_{n\to\infty}(\Re p_{n+1}-\Re p_{n})=+\infty\,,
\end{equation}
\begin{equation}\label{eq:hp2}
\lim_{n\to\infty}\Im p_{n}=0\,,
\end{equation}
and for some $\gamma>0$, $\alpha\in\R$, $n'\in\N$, $\mu>0$, $\nu> 1/2$, $m_1\,,m_2>0$
\begin{equation}\label{eq:hom1}
\liminf_{n\to\infty}({\Re}\om_{n+1}-{\Re}\om_{n})=\gamma\,,
\end{equation}
\begin{equation}\label{eq:hom2}
\lim_{n\to\infty}{\Im}\om_n=\alpha
\,,
\qquad
r_n\le -{\Im}\om_n\,\qquad\forall\ n\ge n'\,,
\end{equation}
\begin{equation}\label{eq:hom3}
|R_n|\le \frac{\mu}{n^{\nu}}|C_n|\,\quad\forall\ n\ge n'\,,
\qquad
|R_n|\le \mu|C_n|\,\quad\forall\ n\le n'\,,
\end{equation}
\begin{equation}\label{eq:alldn}
m_1|p_n|^2\le|d_n|\le m_2|p_n|^2\,\qquad\forall n\in\N\,.
\end{equation}
We note that from \eqref{eq:hp1} it follows
\begin{equation}\label{eq:cesaro}
\lim_{n\to\infty}\frac{\Re p_n}n=+\infty\,,
\end{equation}
see \cite[p. 54]{Ce} and  from $\lim_{n\to\infty}|p_n|=+\infty$ and $p_n\not=0$ it follows
that
there exists $a_0>0$ such that
\begin{equation}\label{eq:a0}
|p_n|\ge a_0\qquad\forall n\in\N\,.
\end{equation}

\subsection{Preliminary results}
First, to prove inverse type estimates we need to introduce an auxiliary function, see \cite{Ed}.
Indeed, we  define
\begin{equation}\label{eq:k}
k(t):=\left \{\begin{array}{l}
\displaystyle\sin \frac{\pi t}{T}\,\qquad\qquad \mbox{if}\,\, t\in\ [0,T]\,,\\
\\
0\,\qquad\qquad\quad\  \ \ \  \mbox{otherwise}\,.
\end{array}\right .
\end{equation}
For the reader's convenience, we list some easy to check properties of $k$ in the following lemma.

\begin{lemma} \label{th:k}
Set
\begin{equation}\label{eqn:K}
K(u):=\frac{\pi T}{\pi^2-T^2u^2}\,,\qquad u\in \C\,,
\end{equation}
the following properties hold  for any $u\in \C$
\begin{equation}\label{eqn:sinek1}
\int_{0}^{\infty} k(t)e^{iu t}dt
= (1+e^{iu T})K(u)
\,,
\end{equation}
\begin{equation}\label{eqn:sinek2bis}
\overline{K(u)}=K(\overline{u})\,,
\end{equation}
\begin{equation}\label{eqn:sinek2}
\big|K(u)\big|=\big|K(\overline{u})\big|\,,
\end{equation}
\begin{equation}\label{eq:sinek3}
\big|K(u)\big|\le \frac{\pi T}{|T^2(\Re u)^2-T^2(\Im u)^2-\pi^2|}\,.
\end{equation}
\end{lemma}
The following result is a crucial tool in the proof of Ingham type inverse estimate.
\begin{proposition}\label{pr:En}
Under assumptions {\rm (\ref{eq:hp1})--(\ref{eq:hp2})},
for $T>0$, $\varepsilon\in (0,1)$, $M>\frac{2\pi}{T(1-\varepsilon)}$ and
for any complex number sequence $\{E_n\}$ with $\sum_{n=1}^\infty\ |E_{n}|^2<+\infty$  there exists 
$n_0=n_0(\varepsilon)\in\N$ independent of coefficients $E_n$ such that if $E_n=0$ for $n<n_0$,
then we have
\begin{multline}\label{eq:En2}
\left|\int_{0}^{\infty} k(t)\Big| \sum_{n=n_0}^\infty  E_{n}e^{i p_{n}t}+\overline{E_{n}}e^{-i \overline{p_n}t}\Big|^2\ d t
-2\pi T\sum_{n=n_0}^\infty\ \frac{1+e^{-2\Im p_{n} T}}{\pi^2+4T^2(\Im p_{n})^2}|E_{n}|^2\right|
\\
\le
\frac{4\pi}{TM^2}\sum_{n=n_0}^\infty\ (1+e^{-2\Im p_{n} T})|E_{n}|^2
\,.
\end{multline}

\end{proposition}

\begin{Proof}
First of all, we note that by using (\ref{eqn:sinek1}) we have
\begin{multline*}
\int_{0}^{\infty} k(t)\Big| \sum_{n=1}^\infty E_{n}e^{i p_{n}t}
+\overline{E_{n}}e^{-i \overline{p_{n}}t}\Big|^2\ d t
\\
=\int_{0}^{\infty} k(t)\sum_{n=1}^\infty\big(E_{n}e^{i p_{n}t}+\overline{E_{n}}e^{-i \overline{p_{n}}t}\big)
\sum_{m=1}^\infty\big(\overline{E_{m}}e^{-i\overline{p_{m}}t}+E_{m}e^{i p_{m}t}\big)\ dt
\\
=\sum_{n, m=1}^\infty E_{n}\overline{E_{m}} (1+e^{i( p_{n}-\overline{p_{m}}) T})K( p_{n}-\overline{p_{m}})
+\sum_{n, m=1}^\infty E_{n} E_{m} (1+e^{i( p_{n}+p_{m}) T})K(p_{n}+p_{m})
\\
+\sum_{n, m=1}^\infty \overline{E_{n}E_{m}} (1+e^{-i(\overline{ p_n+p_m}) T})K( \overline{p_n+p_m})
+\sum_{n, m=1}^\infty\overline{E_{n}}E_{m} (1+e^{-i(\overline{ p_n}-p_m) T})K( \overline{p_n}-p_{m})
\,.
\end{multline*}
In view of \eqref{eqn:sinek2bis} we have
\begin{multline*}
\int_{0}^{\infty} k(t)\Big| \sum_{n=1}^\infty E_{n}e^{i p_{n}t}
+\overline{E_{n}}e^{-i \overline{p_{n}}t}\Big|^2\ d t
\\
=2\sum_{n, m=1}^\infty \Re \big[E_{n}\overline{E_{m}} (1+e^{i( p_{n}-\overline{p_m}) T})K( p_{n}-\overline{p_m})\big]
\\
+2\sum_{n, m=1}^\infty \Re \big[E_{n}E_{m} (1+e^{i( p_{n}+p_{m}) T})K( p_{n}+p_{m})\big]
\,.
\end{multline*}
We note that for $m=n$
\begin{equation*}
K( p_{n}-\overline{p_n})=K(2i\Im p_{n})=\frac{\pi T}{\pi^2+4T^2(\Im p_{n})^2}\,.
\end{equation*}
Therefore, we deduce
\begin{multline*}
\int_{0}^{\infty} k(t)\Big| \sum_{n=1}^\infty E_{n}e^{i p_{n}t}+\overline{E_{n}}e^{-i \overline{p_n}t}\Big|^2\ d t
-2\pi T\sum_{n=1}^\infty\ \frac{1+e^{-2\Im p_{n} T}}{\pi^2+4T^2(\Im p_{n})^2}|E_{n}|^2
\\
=2\sum_{n, m=1,n\not=m}^\infty \Re \big[E_{n}\overline{E_{m}} (1+e^{i( p_{n}-\overline{p_m}) T})
K( p_{n}-\overline{p_m})\big]
\\
+2\sum_{n, m=1}^\infty \Re \big[E_{n}E_{m} (1+e^{i( p_{n}+p_{m}) T})K( p_{n}+p_{m})\big]
\,,
\end{multline*}
whence
\begin{multline}\label{eq:F2up}
\left|\int_{0}^{\infty} k(t)\Big| \sum_{n=1}^\infty E_{n}e^{i p_{n}t}+\overline{E_{n}}e^{-i\overline{ p_n}t}\Big|^2\ d t
-2\pi T\sum_{n=1}^\infty\ \frac{1+e^{-2\Im p_{n} T}}{\pi^2+4T^2(\Im p_{n})^2}|E_{n}|^2\right|
\\
\le
2\sum_{n, m=1,n\not=m}^\infty |E_{n}| |E_{m}|(1+e^{-(\Im p_{n}+\Im p_{m}) T})
|K( p_{n}-\overline{p_m})|
\\
+2\sum_{n, m=1}^\infty |E_{n}| |E_{m}|(1+e^{-(\Im p_{n}+\Im p_{m}) T})
|K( p_{n}+p_{m})|
\,.
\end{multline}
We observe that,
in virtue of \eqref{eqn:sinek2}, we have
$$
|K( p_{n}-\overline{p_m})|=|K( p_{m}-\overline{p_n})|\,,
$$
whence
\begin{multline*}
\sum_{n, m=1,n\not=m}^\infty |E_{n}||E_{m}| |K( p_{n}-\overline{p_m})|
\le
\frac{1}{2}\sum_{n, m=1,n\not=m}^\infty\ \big(|E_{n}|^2+ |E_{m}|^2\big) |K( p_{n}-\overline{p_m})|
\\
=
\frac{1}{2}\sum_{n=1}^\infty\ |E_{n}|^2\sum_{m=1,m\not=n}^\infty\ |K( p_{n}-\overline{p_m})|
+\frac{1}{2}\sum_{m=1}^\infty\ |E_{m}|^2\sum_{n=1,n\not=m}^\infty\ |K( p_{m}-\overline{p_n})|
\\
= \sum_{n=1}^\infty\ |E_{n}|^2\sum_{m=1,m\not=n}^\infty\ |K( p_{n}-\overline{p_m})|
\,.
\end{multline*}
Similarly, we have
\begin{equation*}
\sum_{n, m=1}^\infty |E_{n}||E_{m}|e^{-(\Im p_{n}+\Im p_{m}) T} |K( p_{n}-\overline{p_m})|
\le
 \sum_{n=1}^\infty\ e^{-2\Im p_{n}T}|E_{n}|^2 \sum_{m=1}^\infty\ |K( p_{n}-\overline{p_m})|
\,,
\end{equation*}
\begin{equation*}
\sum_{n, m=1}^\infty |E_{n}||E_{m}|(1+e^{-(\Im p_{n}+\Im p_{m}) T}) |K( p_{n}+p_{m})|
\le
 \sum_{n=1}^\infty\ (1+e^{-2\Im p_{n}T})|E_{n}|^2\sum_{m=1}^\infty\ |K( p_{n}+p_{m})|
\,,
\end{equation*}
so plugging the above inequalities into formula \eqref{eq:F2up}, we obtain 
\begin{multline}\label{eq:En}
\left|\int_{0}^{\infty} k(t)\Big| \sum_{n=1}^\infty E_{n}e^{i p_{n}t}+\overline{E_{n}}e^{-i \overline{p_n}t}\Big|^2\ d t
-2\pi T\sum_{n=1}^\infty\ \frac{1+e^{-2\Im p_{n} T}}{\pi^2+4T^2(\Im p_{n})^2}|E_{n}|^2\right|
\\
\le
2\sum_{n=1}^\infty\ (1+e^{-2\Im p_{n}T})|E_{n}|^2\sum_{m=1,m\not=n}^\infty\ |K( p_{n}-\overline{p_m})|
\\
+2 \sum_{n=1}^\infty\ (1+e^{-2\Im p_{n}T})|E_{n}|^2\sum_{m=1}^\infty\ |K( p_{n}+p_{m})|
\,.
\end{multline}
In the following lemma we single out the estimates concerning the sums depending on $K$ in the right-hand side of the above formula, because we will also use them in the proof of the direct estimate.
\begin{lemma}\label{le:stimaK}
For any $\varepsilon\in (0,1)$ and $M>\frac{2\pi}{T(1-\varepsilon)}$ there exists 
$n_0=n_0(\varepsilon)\in\N$ such that for any $n\ge n_0$ we have
\begin{equation}\label{eq:minus}
\sum_{m=n_0,m\not=n}^\infty\ |K( p_{n}-\overline{p_m})|
\le
\frac{2\pi}{TM^2}\,,
\end{equation}
\begin{eqnarray}\label{eq:plus}
\sum_{m=n_0}^\infty|K( p_n+p_m)|
\le 
\frac{{4\pi}}{TM^2}
\sum_{m=n_0}^\infty\frac{1}{4m^{2}-1}
\,.
\end{eqnarray}
\end{lemma}
\begin{ProofL1}
To prove the first inequality, we observe that, thanks to (\ref{eq:sinek3}), we get  
\begin{eqnarray}\label{eq:modulo0}
\sum_{m=1,m\not=n}^\infty\ |K( p_{n}-\overline{p_m})|
\le
\pi T\sum_{m=1,m\not=n}^\infty\frac{1}
{\big|T^2( \Re p_n-\Re p_m)^2-T^2( \Im p_n+\Im p_m)^2-\pi^2\big|}.
\end{eqnarray}
From assumption (\ref{eq:hp1}) it follows that for any $M>0$ there exists $n_0\in\N$ such that
\begin{equation*}
\Re p_{n+1}-\Re p_n\ge M\qquad\quad\forall n\ge n_0\,,
\end{equation*}
whence
\begin{equation*}
|\Re p_n-\Re p_m|\ge M |n-m|\,,\qquad\forall n\,,m\ge n_0\,.
\end{equation*}
Thanks to the previous estimate, we have
\begin{equation*}
T^2( \Re p_n-\Re p_m)^2-T^2( \Im p_n+\Im p_m)^2-\pi^2
\ge
T^2M^2 (n-m)^2-T^2( \Im p_n+\Im p_m)^2-\pi^2
\,.
\end{equation*}
Moreover,
since $\lim_{n\to\infty}\Im p_n=0$, fix $0<\varepsilon<1$, for $n_0\in\N$ sufficiently large we have 
\begin{equation*}
|\Im p_n|<\frac{M}{4}\varepsilon
\qquad
\forall n\ge n_0
\,,
\end{equation*}
so, for any $n\,,m\in\N$, $n,m\ge n_0$\,, we have
\begin{equation*}\label{eq:}
T^2(\Im p_n+\Im p_m)^2+\pi^2
< \frac{1}{4}\big(T^2M^2\varepsilon^2+4\pi^2\big)\,.
\end{equation*}
Now, as $M>\frac{2\pi}{T(1-\varepsilon)}$ we have $T^2M^2\varepsilon^2+4\pi^2<T^2M^2$, so from the above inequality it follows
\begin{eqnarray}\label{eq:imm0}
T^2(\Im p_n+\Im p_m)^2+\pi^2
<
\frac14T^2M^2
\,,
\end{eqnarray}
and hence 
for $m\not=n$,  
\begin{eqnarray*}\label{eq:}
T^2(\Re p_n-\Re p_m)^2-T^2(\Im p_n+\Im p_m)^2-\pi^2
\ge
T^2M^2(n-m)^2-\frac14T^2M^2>0.
\end{eqnarray*}
Putting the previous formula into (\ref{eq:modulo0}), for any $n\ge n_0$ we obtain 
\begin{multline*}
\sum_{m=n_0,m\not=n}^\infty\ |K( p_{n}-\overline{p_m})|
\\
\le
4\pi T\sum_{m=n_0,m\not=n}^\infty\
\frac{1}{4T^2M^2(m-n)^2-T^2M^2}
=\frac{4\pi}{TM^2}\sum_{m=n_0,m\not=n}^\infty\
\frac{1}{4(m-n)^2-1}\\
\le
\frac{4\pi}{TM^2}\sum_{j=1}^{\infty}\ \frac{1}{4j^2-1}
=\frac{2\pi}{TM^2}\sum_{j=1}^{\infty}\
\Big(\frac{1}{2j-1}-\frac{1}{2j+1}\Big)=
\frac{2\pi}{TM^2}\,,
\end{multline*}
that is \eqref{eq:minus}.

As regards the second estimate, again by (\ref{eq:sinek3}) we have
\begin{eqnarray}\label{eqn:om_n+p_m}
\sum_{m=n_0}^\infty\
|K(p_n+p_m)|
\le 
\pi T\sum_{m=n_0}^\infty\frac{1}
{\big|T^2(\Re p_n+\Re p_m)^2-T^2(\Im p_n+\Im p_m)^2-\pi^2\big|}\,.
\end{eqnarray}
From \eqref{eq:cesaro}, we have for any $M>0$
\begin{eqnarray*}
\Re p_n\ge M n\,,\qquad\forall n\ge n_0\,.
\end{eqnarray*}
By using the previous inequality and \eqref{eq:imm0}, we get for $M>\frac{2\pi}{T(1-\varepsilon)}$
\begin{equation*}
T^2(\Re p_n+\Re p_m)^2-T^2(\Im p_n+\Im p_m)^2-\pi^2
\ge T^2M^2m^2-\frac14T^2M^2
=\frac{T^2M^2}4(4m^2-1)
\,.
\end{equation*}
Therefore from (\ref{eqn:om_n+p_m}), by using the above estimate, we get
\begin{eqnarray*}
\sum_{m=n_0}^\infty|K( p_n+p_m)|
\le 
\frac{{4\pi}}{TM^2}
\sum_{m=n_0}^\infty\frac{1}{4m^{2}-1}
\,,
\end{eqnarray*}
that is \eqref{eq:plus}.
\end{ProofL1}
{\it Proof of Proposition~{\rm \ref{pr:En}} (continued)}.
If we assume $E_n=0$ for any $n< n_0$, then from \eqref{eq:En} and \eqref{eq:minus} we have
\begin{multline}\label{eq:En1}
\left|\int_{0}^{\infty} k(t)\Big| \sum_{n=n_0}^\infty E_{n}e^{i p_{n}t}+\overline{E_{n}}e^{-i \overline{p_n}t}\Big|^2\ d t
-2\pi T\sum_{n=n_0}^\infty\ \frac{1+e^{-2\Im p_{n} T}}{\pi^2+4T^2(\Im p_{n})^2}|E_{n}|^2\right|
\\
\le
\frac{4\pi}{TM^2}\sum_{n=n_0}^\infty\ (1+e^{-2\Im p_{n}T})|E_{n}|^2
+2 \sum_{n=n_0}^\infty\ (1+e^{-2\Im p_{n}T})|E_{n}|^2\sum_{m=n_0}^\infty\ |K( p_{n}+p_{m})|
\,.
\end{multline}
Now, we observe that for $m\ge n_0$ we have
\begin{equation*}
4m^2-1
\ge 4m^{3/2}n_0^{1/2}-1
\ge n^{1/2}_0(4m^{3/2}-1)
\,,
\end{equation*}
whence
\begin{equation*}
\sum_{m=n_0}^\infty\frac{1}{4m^{2}-1}
\le
\frac{1}{n^{1/2}_0}
\sum_{m=1}^\infty\frac{1}{4m^{3/2}-1}\,.
\end{equation*}
Therefore from (\ref{eq:plus}), by using the above inequality, we get
\begin{eqnarray*}
\sum_{m=n_0}^\infty|K( p_n+p_m)|
\le 
\frac{{4\pi}}{TM^2n^{1/2}_0}
\sum_{n=1}^\infty\frac{1}{4n^{3/2}-1}
\,,
\end{eqnarray*}
whence
\begin{equation*}
\sum_{n=n_0}^\infty(1+e^{-2\Im p_{n}T})|E_n|^2\sum_{m=n_0}^\infty|K( p_{n}+p_{m})|
\le
\frac{{4\pi}}{TM^2n^{1/2}_0}\sum_{n=1}^\infty\frac{1}{4n^{3/2}-1}\sum_{n=n_0}^\infty(1+e^{-2\Im p_{n}T})|E_n|^2\,.
\end{equation*}
If we choose $n_0\in\N$ large enough to satisfy the condition
\begin{equation*}
\frac{{2}}{n^{1/2}_0}\sum_{n=1}^\infty\frac{1}{4n^{3/2}-1}<\varepsilon\,,
\end{equation*}
we have
\begin{equation*}
\sum_{n=n_0}^\infty(1+e^{-2\Im p_{n}T})|E_n|^2\sum_{m=n_0}^\infty|K( p_{n}+p_{m})|
\le
\frac{{2\pi}}{TM^2}\varepsilon\sum_{n=n_0}^\infty(1+e^{-2\Im p_{n}T})|E_n|^2\,.
\end{equation*}
Plugging the above inequality into \eqref{eq:En1} we get
\begin{multline*}
\left|\int_{0}^{\infty} k(t)\Big| \sum_{n=n_0}^\infty  E_{n}e^{i p_{n}t}+\overline{E_{n}}e^{-i \overline{p_n}t}\Big|^2\ d t
-2\pi T\sum_{n=n_0}^\infty\ \frac{1+e^{-2\Im p_{n} T}}{\pi^2+4T^2(\Im p_{n})^2}|E_{n}|^2\right|
\\
\le
\frac{4\pi}{TM^2}(1+\varepsilon)\sum_{n=n_0}^\infty\ (1+e^{-2\Im p_{n} T})|E_{n}|^2
\,.
\end{multline*}
Finally, by substituting $M$ with $M\sqrt{1+\varepsilon}$ we obtain \eqref{eq:En2} .
\end{Proof}
\vskip1cm
As for the inverse inequality, to prove  direct estimates we need to introduce an auxiliary function. Let $T>0$ and define  
\begin{equation}\label{eq:kcos}
k^*(t):=\left \{\begin{array}{l}
\cos \frac{\pi t}{2T}\,\qquad\qquad \mbox{if}\ |t|\le T\,,\\
\\
0\,\qquad\qquad\quad\  \ \ \  \mbox{if}\  |t|>T\,.
\end{array}\right .
\end{equation}
For the sake of completeness, we list some easy to check properties of $k^*$ in the following lemma.
\begin{lemma} \label{th:k}
Set
\begin{equation*}\label{}
K^*(u):=\frac{4T\pi}{\pi^2-4T^2u^2}\,,\qquad u\in \C\,,
\end{equation*}
the following properties hold  for any $u\in \C$
\begin{equation}\label{eqn:k1}
\int_{-\infty}^{\infty} k^*(t)e^{iu t}dt=\cos(uT)K^*(u)\,,
\end{equation}
\begin{equation}\label{eqn:k2bis}
\overline{K^*(u)}=K^*(\overline{u})\,,
\end{equation}
\begin{equation}\label{eqn:k2}
\big|K^*(u)\big|=\big|K^*(\overline{u})\big|\,.
\end{equation}
If we set 
$K_{T}(u)=\frac{T\pi}{\pi^2-T^2u^2}$, then we have
\begin{equation}\label{eqn:k3}
K^*(u)=2K_{2T}(u)\,.
\end{equation}

\end{lemma}
From now on $c(T)$ will denote a positive constant depending on $T$. 
\begin{proposition}\label{pr:EnD}
Assume {\rm (\ref{eq:hp1})--(\ref{eq:hp2})}.
Let $T>0$, $\varepsilon\in(0,1)$, $M>\frac{\pi}{T\sqrt{1-\varepsilon}}$ and $\{E_n\}$
a complex number sequence such that $\sum_{n=1}^\infty\ |E_{n}|^2<+\infty$.  
There exists $n_0=n_0(\varepsilon)\in\N$ such that if $E_n=0$ for $n<n_0$,
then we have
\begin{equation}\label{eq:EnD2}
\int_{-\infty}^{\infty} k^*(t)\Big| \sum_{n=n_0}^\infty E_{n}e^{i p_{n}t}+\overline{E_{n}}e^{-i \overline{p_n}t}\Big|^2\ d t
\le
4c(T)\Big(\frac{2T}{\pi}+\frac{\pi}{TM^2}\Big(1+\sum_{n=1}^\infty\frac{1}{4n^{2}-1}\Big)\Big)
\sum_{n=n_0}^\infty |E_{n}|^2  
\,.
\end{equation}
\end{proposition}
\begin{Proof}
Let $k^*(t)$  be the function defined by (\ref{eq:kcos}). If we use (\ref{eqn:k1}) and \eqref{eqn:k2bis}, then we have
\begin{multline*}
\int_{-\infty}^{\infty} k^*(t)\Big| \sum_{n=1}^\infty E_{n}e^{i p_{n}t}+\overline{E_{n}}e^{-i \overline{p_n}t}\Big|^2\ d t=
\\
\int_{-\infty}^{\infty} k^*(t)\sum_{n=1}^\infty\big(E_{n}e^{i p_{n}t}+\overline{E_{n}}e^{-i \overline{p_n}t}\big)
\sum_{m=1}^\infty\big(\overline{E_{m}}e^{-i\overline{p_m}t}+E_{m}e^{i p_{m}t}\big)\ dt
\\
=\sum_{n, m=1}^\infty E_{n}\overline{E_{m}} \cos(( p_{n}-\overline{p_m}) T)K^*( p_{n}-\overline{p_m})
+\sum_{n, m=1}^\infty E_{n} E_{m} \cos(( p_{n}+p_{m}) T)K^*(p_{n}+p_{m})
\\
+\sum_{n, m=1}^\infty \overline{E_{n}E_{m}} \cos((\overline{p_{n}+p_{m}})T)K^*( \overline{p_{n}+p_{m}})
+\sum_{n, m=1}^\infty\overline{E_{n}}E_{m} \cos(( \overline{p_n}-p_{m}) T)K^*(\overline{ p_n}-p_{m})
\\
=2\sum_{n, m=1}^\infty \Re \big(E_{n}\overline{E_{m}}\cos(( p_{n}-\overline{p_m}) T)K^*( p_{n}-\overline{p_m})\big)
\\
+2\sum_{n, m=1}^\infty \Re \big(E_{n}E_{m}\cos(( p_{n}+p_{m}) T)K^*( p_{n}+p_{m})\big)\,.
\end{multline*}
Applying  the elementary estimates $\Re z\le | z|$ and $|\cos z|\le\cosh(\Im z)$, $z\in\C$,
we obtain
\begin{multline*}
\int_{-\infty}^{\infty} k^*(t)\Big| \sum_{n=1}^\infty E_{n}e^{i p_{n}t}+\overline{E_{n}}e^{-i \overline{p_n}t}\Big|^2\ d t
\\
\le
2\sum_{n, m=1}^\infty |E_{n}| |E_{m}|\cosh((\Im p_n+\Im p_m)T)\big[ |K^*( p_{n}-\overline{p_m})|+|K^*( p_{n}+p_{m})|\big]
\,.
\end{multline*}
Since the sequence $\{{\Im}p_n\}$ is bounded, 
for any $n,m\in\N$, we have
$$
\cosh((\Im p_n+\Im p_m)T)\le e^{2T\sup|\Im p_{n}|}\,,
$$
and hence
\begin{multline*}
\int_{-\infty}^{\infty} k^*(t)\Big| \sum_{n=1}^\infty E_{n}e^{i p_{n}t}+\overline{E_{n}}e^{-i \overline{p_n}t}\Big|^2\ d t
\le
2e^{2T\sup|\Im p_{n}|}\sum_{n, m=1}^\infty |E_{n}| |E_{m}|\big[ |K^*( p_{n}-\overline{p_m})|+|K^*( p_{n}+p_{m})|\big]
\,.
\end{multline*}
In virtue of \eqref{eqn:k2} we get
$
|K^*( p_n-\overline{p_m})|=|K^*( p_m-\overline{p_n})|\,,
$
so we have
\begin{multline*}
\int_{-\infty}^{\infty} k^*(t)\Big| \sum_{n=1}^\infty E_{n}e^{i p_{n}t}+\overline{E_{n}}e^{-i \overline{p_n}t}\Big|^2\ d t
\le
2e^{2T\sup|\Im p_{n}|}\sum_{n=1}^\infty |E_{n}|^2\sum_{m=1}^\infty\big[ |K^*( p_{n}-\overline{p_m})|+|K^*( p_{n}+p_{m})|\big]
\,.
\end{multline*}
Taking into account the definition of $K^*$ we have
\begin{multline}\label{eq:dir1}
\int_{-\infty}^{\infty} k^*(t)\Big| \sum_{n=1}^\infty E_{n}e^{i p_{n}t}+\overline{E_{n}}e^{-i\overline{ p_n}t}\Big|^2\ d t
\\
\le
8T\pi e^{2T\sup|\Im p_{n}|}\sum_{n=1}^\infty\frac{1}{\pi^2+8T^2(\Im p_n)^2} |E_{n}|^2 
+2e^{2T\sup|\Im p_{n}|}\sum_{n=1}^\infty |E_{n}|^2\sum_{m=1,m\not=n}^\infty |K^*( p_{n}-\overline{p_m})|
\\
+2e^{2T\sup|\Im p_{n}|}\sum_{n=1}^\infty |E_{n}|^2\sum_{m=1}^\infty|K^*( p_{n}+p_{m})|
\,. 
\end{multline} 
Now, we note that in virtue of \eqref{eqn:k3} we can apply Lemma \ref{le:stimaK}:
for any $\varepsilon\in (0,1)$ and $M>\frac{\pi}{T\sqrt{1-\varepsilon}}=\frac{2\pi}{2T\sqrt{1-\varepsilon}}$
there exists $n_0\in\N$  such that for any $n\ge n_0$
\begin{equation*}
\sum_{m=n_0,m\not=n}^\infty |K^*( p_n-\overline{p_m})|
\le
\frac{4\pi}{2TM^2}=\frac{2\pi}{TM^2}
\,, 
\end{equation*}
\begin{eqnarray*}
\sum_{m=n_0}^\infty|K^*( p_n+p_m)|
\le 
\frac{{4\pi}}{2TM^2}
\sum_{m=n_0}^\infty\frac{1}{4m^{2}-1}
\le 
\frac{{2\pi}}{TM^2}
\sum_{n=1}^\infty\frac{1}{4n^{2}-1}
\,.
\end{eqnarray*}
In conclusion, assuming $E_n=0$ for $n< n_0$ and putting the above formulas into (\ref{eq:dir1}), we get
\begin{equation*}
\int_{-\infty}^{\infty} k^*(t)\Big| \sum_{n=n_0}^\infty E_{n}e^{i p_{n}t}+\overline{E_{n}}e^{-i \overline{p_n}t}\Big|^2\ d t
\le
e^{2T\sup|\Im p_{n}|}\Big(\frac{8T}{\pi}+\frac{4\pi}{TM^2}\Big(1+\sum_{n=1}^\infty\frac{1}{4n^{2}-1}\Big)\Big)
\sum_{n=n_0}^\infty |E_{n}|^2 
\,,
\end{equation*} 
that is \eqref{eq:EnD2}.
\end{Proof}

\subsection{Inverse and direct inequalities excluding a finite number of terms.}

Due to the asymptotic assumptions on data,  some properties hold true for sufficiently large integers. For that reason,  first we will show some inverse and direct inequalities in the special case when our series have a finite number of terms vanishing.

Before proceeding, we state the next result, that  can be proved in the same way as in 
\cite[Theorem $5.3$]{LoretiSforza1}, taking into account that the function $k(t)$ is non negative.

From now on we denote with $c(T,\varepsilon)$ a positive constant depending on $T$ and $\varepsilon$.
\begin{theorem}\label{th:Cn}
Under assumptions {\rm (\ref{eq:hom1})--(\ref{eq:hom3})}, for any $\varepsilon\in (0,1)$ and 
$T>\frac{2\pi}{\ga(1-\varepsilon)}$ there exist $n_0=n_0(\varepsilon)\in\N$ and $c(T,\varepsilon)>0$ such that if $C_n=0$ for any
$n< n_0$, then  we have
\begin{equation}\label{eq:inv-sine}
\int_{0}^{\infty} k(t)\Big|\sum_{n=n_0}^\infty R_{n}e^{r_{n}t}+C_{n}e^{i\omega_{n}t}+\overline{C_n}e^{-i\overline{\omega_n}t}\Big|^2\ dt
\ge
c(T,\varepsilon)\sum_{n= n_0}^\infty(1+e^{-2({\Im}\om_n-\alpha)T})|C_n|^2
\,.
\end{equation}
\end{theorem}
In the following finding we give a lower bound for the first component of the solution of coupled system.
\begin{theorem}\label{th:CnDn}
Under assumptions {\rm (\ref{eq:hp1})--(\ref{eq:hom3})}, for any $\varepsilon\in (0,1)$ and
$T>\frac{2\pi}{\ga(1-\varepsilon)}$,
 there exist $n_0=n_0(\varepsilon)\in\N$ and $c(T,\varepsilon)>0$ such that if $C_n=D_n=0$ for any
$n< n_0$, then  we have
\begin{multline}\label{eq:CnDn}
\int_{0}^{T} \Big|\sum_{n=n_0}^\infty R_{n}e^{r_{n}t}+C_{n}e^{i\omega_{n}t}+\overline{C_n}e^{-i\overline{\omega_n}t}+D_{n}e^{i p_{n}t}+\overline{D_{n}}e^{-i\overline{ p_n}t}\Big|^2\ dt
\\
\ge
c(T,\varepsilon)\sum_{n= n_0}^\infty(1+e^{-2({\Im}\om_n-\alpha)T})|C_n|^2
\\
-2\pi T\sum_{n=n_0}^\infty\ \Big( \frac{1}{\pi^2+4T^2(\Im p_{n})^2}+\frac{2}{ T^2\gamma^2}\Big)(1+e^{-2{\Im}p_n T})|D_{n}|^2
\,.
\end{multline}
\end{theorem}

\begin{Proof}
First of all, we set for any $t\ge 0$
\begin{equation*}
F_1(t)=\sum_{n=1}^\infty\big(R_ne^{r_nt}+C_ne^{i\om_nt}+\overline{C_n}e^{-i\overline{\omega_n}t}\big)\in\R\,,
\end{equation*}
\begin{equation*}
F_2(t)=\sum_{n=1}^\infty \big(D_{n}e^{i p_{n}t}+\overline{D_{n}}e^{-i \overline{p_n}t}\big)\in\R\,,
\end{equation*}
and observe that
if $k(t)$  is the function defined by (\ref{eq:k}), we have to estimate  the term
\begin{equation*}
\int_{0}^{\infty} k(t)| F_1(t)+ F_2(t)|^2\ dt\,.
\end{equation*}
Because of the elementary inequality $2|ab|\le\frac12a^2+2b^2$, we observe that
\begin{multline*}
| F_1(t)+ F_2(t)|^2=| F_1(t)|^2+2F_1(t)F_2(t)+ |F_2(t)|^2
\\
\ge
| F_1(t)|^2-\frac12| F_1(t)|^2-2|F_2(t)|^2+ |F_2(t)|^2
=\frac12| F_1(t)|^2-|F_2(t)|^2\,.
\end{multline*}
Since $k(t)$ is positive, from the above inequality we have
\begin{equation*}
\int_{0}^{\infty} k(t)| F_1(t)+ F_2(t)|^2\ dt
\ge
\frac12\int_{0}^{\infty} k(t)| F_1(t)|^2\ d t -\int_{0}^{\infty} k(t)| F_2(t)|^2\ d t\,.
\end{equation*}
Therefore, in view of Theorem \ref{th:Cn} we can apply \eqref{eq:inv-sine} to get
\begin{multline}\label{eq:u1inv}
\int_{0}^{\infty} k(t)| F_1(t)+ F_2(t)|^2\ dt
\ge
\frac12 c(T,\varepsilon)\sum_{n= n_0}^\infty(1+e^{-2({\Im}\om_n-\alpha)T})|C_n|^2
-\int_{0}^{\infty} k(t)| F_2(t)|^2\ d t\,,
\end{multline}
for $n_0$ sufficiently large.
To complete our proof,
we must give an upper bound for the term 
$\int_{0}^{\infty} k(t)| F_2(t)|^2\ d t$.
Indeed, if we take $E_n=D_n$  and $M=\gamma>\frac{2\pi}{T(1-\varepsilon)}$ in Proposition \ref{pr:En}, then by formula \eqref{eq:En2} we have
\begin{multline*}
\int_{0}^{\infty} k(t)| F_2(t)|^2\ d t
=
\int_{0}^{\infty} k(t)\Big|\sum_{n=n_0}^\infty D_{n}e^{i p_{n}t}+\overline{D_{n}}e^{-i \overline{p_n}t}\Big|^2\ d t
\\
\le
2\pi T\sum_{n=n_0}^\infty\ \frac{1+e^{-2\Im p_{n} T}}{\pi^2+4T^2(\Im p_{n})^2}|D_{n}|^2
+\frac{4\pi}{ T\gamma^2}\sum_{n=n_0}^\infty\ (1+e^{-2{\Im}p_n T})|D_{n}|^2
\\
=
2\pi T\sum_{n=n_0}^\infty\ \Big( \frac{1}{\pi^2+4T^2(\Im p_{n})^2}+\frac{2}{ T^2\gamma^2}\Big)(1+e^{-2{\Im}p_n T})|D_{n}|^2
\,.
\end{multline*}
Therefore, putting the above estimate in \eqref{eq:u1inv}  we have
\begin{multline*}
\int_{0}^{\infty} k(t)| F_1(t)+ F_2(t)|^2\ dt
\\
\ge
\frac12 c(T,\varepsilon)\sum_{n= n_0}^\infty(1+e^{-2({\Im}\om_n-\alpha)T})|C_n|^2
-2\pi T\sum_{n=n_0}^\infty\ \Big( \frac{1}{\pi^2+4T^2(\Im p_{n})^2}+\frac{2}{ T^2\gamma^2}\Big)(1+e^{-2{\Im}p_n T})|D_{n}|^2
\,,
\end{multline*}
whence, in virtue of the definition of $k(t)$, \eqref{eq:CnDn} follows.  
\end{Proof}
\begin{proposition}\label{pr:Dn}
Assume {\rm (\ref{eq:hp1})}, {\rm (\ref{eq:hp2})} and {\rm (\ref{eq:alldn})}.
Let $ T>0$, $\varepsilon\in (0,1)$ and 
$M>\frac{2\pi}{T(1-\varepsilon)}$.  There exist $n_0=n_0(\varepsilon)\in\N$ 
and
$c(T,M,\varepsilon)>0$
such that if $D_n=0$ for any
$n< n_0$, then  we have
\begin{multline}\label{eq:dn}
\int_{0}^{T} \Big|e^{\eta t} \sum_{n=n_0}^{\infty}\big(d_nD_{n}e^{i p_{n}t}
+\overline{d_{n}}\overline{D_{n}}e^{-i \overline{p_n}t}\big)
\Big|^2\ dt
\\
\ge
2\pi Tm_1^2\sum_{n=n_0}^\infty\ \Big(\frac{1}{\pi^2+4T^2(\Im p_{n})^2}-\frac{2}{T^2M^2}\Big)(1+e^{-2\Im p_{n} T})|D_{n}|^2 |p_n|^4
\,,
\end{multline}
and
\begin{equation}\label{eq:cpos}
\frac{1}{\pi^2+4T^2(\Im p_{n})^2}-\frac{2}{T^2M^2}>c(T,M,\varepsilon)\,, \qquad\forall n\ge n_0\,.
\end{equation}

\end{proposition}
\begin{Proof}
We will use Proposition \ref{pr:En} again. Indeed, if
we set
\begin{equation*}
G(t)=
\sum_{n=1}^{\infty}\Big(d_nD_{n}e^{i p_{n}t}
+\overline{d_{n}}\overline{D_{n}}e^{-i \overline{p_n}t}\Big)
\end{equation*}
and take 
$E_n=d_nD_{n}$,  
we can apply formula \eqref{eq:En2} for $n_0$ large enough:
\begin{equation}\label{eq:constn}
\int_{0}^{\infty} k(t)| G(t)|^2\ d t 
\ge
2\pi T\sum_{n=n_0}^\infty\ \Big(\frac{1}{\pi^2+4T^2(\Im p_{n})^2}-\frac{2}{T^2M^2}\Big)
 (1+e^{-2\Im p_{n} T}) |D_{n}|^2 |d_n|^2
\,,
\end{equation}
where $k(t)$  is the function defined by (\ref{eq:k}).
Since $\lim_{n\to\infty}\Im p_n=0$, by taking $n_0$ large enough we have
\begin{equation*}
|\Im p_n|<\frac{M\sqrt\varepsilon}{2\sqrt2}
\qquad\forall n\ge n_0\,,
\end{equation*}
and hence, since $M>\frac{2\pi}{T(1-\varepsilon)}$, we get
\begin{equation*}
\frac{1}{\pi^2+4T^2(\Im p_{n})^2}-\frac{2}{T^2M^2}>\frac{1}{\pi^2+T^2M^2\varepsilon/2}-\frac{2}{T^2M^2}
=2\frac{T^2M^2(1-\varepsilon)-2\pi^2}{T^2M^2(2\pi^2+T^2M^2\varepsilon)}
\,,
\end{equation*}
that is, \eqref{eq:cpos} holds true with
$
c(T,M,\varepsilon)=2\frac{T^2M^2(1-\varepsilon)-2\pi^2}{T^2M^2(2\pi^2+T^2M^2\varepsilon)}>0\,.
$
%
%
Thanks to {\rm (\ref{eq:alldn})}, we get
\begin{equation*}
\int_{0}^{\infty} k(t)| G(t)|^2\ dt
\ge
2\pi Tm^2_1\sum_{n=n_0}^\infty\ \Big(\frac{1}{\pi^2+4T^2(\Im p_{n})^2}-\frac{2}{T^2M^2}\Big) (1+e^{-2\Im p_{n}  T})|D_{n}|^2 |p_n|^4\,,
\end{equation*}
whence, in virtue of the definition of $k(t)$, 
\begin{equation}\label{eq:noexp}
\int_{0}^{T} | G(t)|^2\ dt
\ge
2\pi Tm_1^2\sum_{n=n_0}^\infty\ \Big(\frac{1}{\pi^2+4T^2(\Im p_{n})^2}-\frac{2}{T^2M^2}\Big) (1+e^{-2\Im p_{n}  T})|D_{n}|^2 |p_n|^4\,.
\end{equation}
Finally, in view of
\begin{equation*}
\int_{0}^{T} | e^{\eta t}G(t)|^2\ dt
=\int_{0}^{T}e^{2\eta t} | G(t)|^2\ dt
\ge\int_{0}^{T}| G(t)|^2\ dt\,,
\end{equation*}
by \eqref{eq:noexp} it follows
\eqref {eq:dn}.
\end{Proof}
Now, we anticipate a result concerning direct estimates, because we will use it in the next theorem.
\begin{proposition}\label{pr:dirDn0}
Assume {\rm (\ref{eq:hp1})}, {\rm (\ref{eq:hp2})} and {\rm (\ref{eq:alldn})}.
Let $T>0$, $\varepsilon\in (0,1)$ and 
$M>\frac{\pi}{T\sqrt{1-\varepsilon}}$.  There exist $c(T)>0$ and  $n_0=n_0(\varepsilon)\in\N$  such that if $D_n=0$ for any
$n< n_0$, then  we have
\begin{equation}\label{eq:dirDn0}
\int_{-\infty}^{\infty} k^*(t)\Big| \sum_{n=n_0}^{\infty}d_nD_{n}e^{i p_{n}t}
+\overline{d_{n}}\overline{D_{n}}e^{-i \overline{p_n}t}
\Big|^2\ dt
\le
c(T)\sum_{n=n_0}^\infty\ |D_{n}|^2 |p_n|^4
\,,
\end{equation}
\begin{equation}\label{eq:dirDn1}
\int_{-T}^{T}\Big| \sum_{n=n_0}^{\infty}d_nD_{n}e^{i p_{n}t}
+\overline{d_{n}}\overline{D_{n}}e^{-i \overline{p_n}t}
\Big|^2\ dt
\le
c(T)\sum_{n=n_0}^\infty\ |D_{n}|^2 |p_n|^4
\,.
\end{equation}
\end{proposition}
\begin{Proof}
We evaluate the integral by using Proposition \ref{pr:EnD}: indeed, if we take 
$$
E_n=d_nD_{n}\,,
$$
then from \eqref{eq:EnD2} it follows
\begin{equation}\label{eq:dir2.1}
\int_{-\infty}^{\infty} k^*(t)\Big| \sum_{n=n_0}^{\infty}d_nD_{n}e^{i p_{n}t}
+\overline{d_{n}}\overline{D_{n}}e^{-i \overline{p_n}t}\Big|^2\ dt
\le
c(T)
\sum_{n=n_0}^\infty |D_{n}|^2  |d_n|^2
\,.
\end{equation}
Moreover, from (\ref{eq:alldn}) we get
\begin{equation*}
\int_{-\infty}^{\infty} k^*(t)\Big| \sum_{n=n_0}^{\infty}d_nD_{n}e^{i p_{n}t}
+\overline{d_{n}}\overline{D_{n}}e^{-i \overline{p_n}t}\Big|^2\ dt
\le
c(T)
\sum_{n=n_0}^\infty |D_{n}|^2  |p_n|^4
\,,
\end{equation*}
that is \eqref{eq:dirDn0}.

Now, if we consider the last inequality with the function $k^*$ 
replaced by the analogous one relative to $2T$ instead of $T$,
see (\ref{eq:kcos}),
then 
we get 
\begin{equation*}
\int_{-2T}^{2T} \cos \frac{\pi t}{4T}\Big| \sum_{n=n_0}^{\infty}d_nD_{n}e^{i p_{n}t}
+\overline{d_{n}}\overline{D_{n}}e^{-i \overline{p_n}t}\Big|^2\ dt
\\
\le
c(2T)
\sum_{n=n_0}^\infty\ |D_{n}|^2|p_n|^4 ,
\end{equation*}
whence, thanks to
$
\cos \frac{\pi t}{4T}\ge \frac1{\sqrt2}
$
for $|t|\le T$,
it follows
\begin{equation*}
\int_{-T}^{T} \Big| \sum_{n=n_0}^{\infty}d_nD_{n}e^{i p_{n}t}
+\overline{d_{n}}\overline{D_{n}}e^{-i \overline{p_n}t}\Big|^2\ dt
\le
\sqrt2c(2T)
\sum_{n=n_0}^\infty\ |D_{n}|^2|p_n|^4 \,,
\end{equation*}
that is \eqref{eq:dirDn1}.
\end{Proof}

\begin{theorem}\label{th:Dn+D}
Assume {\rm (\ref{eq:hp1})}, {\rm (\ref{eq:hp2})} and {\rm (\ref{eq:alldn})}.
Let $T> T_0>0$, 
$\varepsilon\in (0,1)$ and 
$M>\frac{2\pi}{T(1-\varepsilon)}$.  
If $\D\in\R$, there exist $n_0=n_0(\varepsilon)\in\N$, $C_0>0$  independent of $T$, $c(T)>0$ and $c(T,M,\varepsilon)>0$  such that if $D_n=0$ for any
$n< n_0$, then  we have
\begin{multline}\label{eq:dn+D}
\int_{0}^{T} \Big| e^{\eta t}\sum_{n=n_0}^{\infty}\big(d_nD_{n}e^{i p_{n}t}
+\overline{d_{n}}\overline{D_{n}}e^{-i \overline{p_n}t}\big)
+\D\Big|^2\ dt
\\
\ge
\pi C_0 T
\sum_{n=n_0}^\infty\ \Big(\frac{1}{\pi^2+4T^2(\Im p_{n})^2}-\frac{2}{T^2M^2}\Big)(1+e^{-2\Im p_{n} T})|D_{n}|^2 |p_n|^4
+c(T) |\D|^2
\,,
\end{multline}
\begin{equation}\label{eq:cpos1}
\frac{1}{\pi^2+4T^2(\Im p_{n})^2}-\frac{2}{T^2M^2}>c(T,M,\varepsilon)\,, \qquad\forall n\ge n_0\,.
\end{equation}
\end{theorem}
\begin{Proof}
We introduce the function
\begin{equation*}
G_1(t)=
e^{\eta t}\sum_{n=1}^{\infty}\big(d_nD_{n}e^{i p_{n}t}
+\overline{d_{n}}\overline{D_{n}}e^{-i \overline{p_n}t}\big)+\D\,.
\end{equation*}
To evaluate the integral of 
$
G_1(t)
$
on the left-hand side of \eqref{eq:dn+D}, we will use the operator introduced by Haraux which annihilates the constant $\D$, see \cite{Ha}.
Indeed, if we take $\delta\in ( T_0/4, T_0/2)$, then we have for any $t\in [0, T-\delta]$
\begin{multline}
\int_0^\delta (G_1(t)-G_1(t+s))\ ds
\\
=
e^{\eta t}\sum_{n=1}^{\infty}d_{n}\Big(\delta-\frac{e^{(\eta+i p_{n})\delta}-1}{\eta+i p_{n}}\Big)
D_{n}e^{i p_{n}t}
+\overline{d_{n}}\Big(\delta-\frac{e^{(\eta-i \overline{p_n})\delta}-1}{\eta-i \overline{p_n}}\Big)
\overline{D_{n}}e^{-i \overline{p_n}t}\,.
\end{multline}
We can apply Proposition \ref{pr:Dn} to the function $t\to\int_0^\delta (G_1(t)-G_1(t+s))\ ds$ in the interval $ [0, T-\delta]$. 
If  $M>\frac{2\pi}{T(1-\varepsilon)}$, we note that $\overline{M}=\frac{T}{T-\delta}M$ verifies 
$\overline{M}>\frac{2\pi}{ (T-\delta)(1-\varepsilon)}$, so by \eqref{eq:dn} and \eqref{eq:cpos} we have
\begin{multline}\label{eq:delta1}
\int_{0}^{ T-\delta} \Big| \int_0^\delta (G_1(t)-G_1(t+s))\ ds\Big|^2\ dt
\\
\ge
2\pi(T-\delta)m_1^2\sum_{n=n_0}^\infty\ \Big(\frac{1}{\pi^2+4(T-\delta)^2(\Im p_{n})^2}-\frac{2}{(T-\delta)^2\overline{M}^2}\Big)
\ (1+e^{-2\Im p_{n} ( T-\delta)})\Big|\delta-\frac{e^{(\eta+i p_{n})\delta}-1}{\eta+i p_{n}}\Big|^2|D_{n}|^2 |p_n|^4
\end{multline}
and
\begin{equation*}
\frac{1}{\pi^2+4(T-\delta)^2(\Im p_{n})^2}-\frac{2}{(T-\delta)^2\overline{M}^2}>0\,, \qquad\forall n\ge n_0\,.
\end{equation*}
We have to estimate $\Big|\delta-\frac{e^{(\eta+i p_{n})\delta}-1}{\eta+i p_{n}}\Big|$. First, we observe that
\begin{equation*}
\Big|\delta-\frac{e^{(\eta+i p_{n})\delta}-1}{\eta+i p_{n}}\Big|
\ge
\delta-\frac{|e^{(\eta+i p_{n})\delta}-1|}{|\eta+ip_{n}|}
\ge
\delta-\frac{e^{(\eta-\Im p_{n})\delta}+1}{|\Re p_{n}|}\,.
\end{equation*}
Since the sequence $\{\Im p_{n}\}$ is bounded and $\delta<T_0$, we have
\begin{equation*}
e^{(\eta-\Im p_{n})\delta}
\le
e^{(\eta+\sup|\Im p_{n}|)\delta}
\le
e^{(\eta+\sup|\Im p_{n}|)T_0}\,,
\end{equation*}
and hence
\begin{equation*}
\Big|\delta-\frac{e^{(\eta+i p_{n})\delta}-1}{\eta+i p_{n}}\Big|
\ge
\delta-\frac{2e^{(\eta+\sup|\Im p_{n}|)T_0}}{|\Re p_{n}|}\,.
\end{equation*}
Taking into account that $\lim_{n\to\infty}\Re p_{n}=+\infty$, for $n_0$ sufficiently large we have for any $n\ge n_0$
\begin{equation*}
\frac{2e^{(\eta+\sup|\Im p_{n}|)T_0}}{\Re p_{n}}
\le\frac{T_0}{4}\,, 
\end{equation*}
whence
\begin{equation*}
\Big|\delta-\frac{e^{(\eta+i p_{n})\delta}-1}{\eta+i p_{n}}\Big|
\ge
\delta-\frac{T_0}{4}>0\,.
\end{equation*}
Plugging the above estimate into \eqref{eq:delta1}, we obtain
\begin{multline}\label{eq:delta2}
\int_{0}^{ T-\delta} \Big| \int_0^\delta (G_1(t)-G_1(t+s))\ ds\Big|^2\ dt
\\
\ge
2\pi(\delta-T_0/{4})^2(T-\delta)m_1^2\sum_{n=n_0}^\infty\ \Big(\frac{1}{\pi^2+4(T-\delta)^2(\Im p_{n})^2}-\frac{2}{(T-\delta)^2\overline{M}^2}\Big)
\ (1+e^{-2\Im p_{n} ( T-\delta)})|D_{n}|^2 |p_n|^4
\,.
\end{multline}
Moreover, since  $2\delta<T_0$ and the sequence $\{\Im p_{n}\}$ is bounded we get
\begin{equation*}
e^{2\Im p_{n} \delta}
\ge
e^{-2 \delta|\Im p_{n}|} 
\ge
e^{-T_0|\Im p_{n}| } 
\ge
e^{-T_0\sup|\Im p_{n}|} \,,
\end{equation*}
whence
\begin{equation*}
1+e^{-2\Im p_{n} ( T-\delta)}
=1+e^{-2\Im p_{n} T}e^{2\Im p_{n}\delta}
\ge
e^{-T_0\sup|\Im p_{n}|} (1+e^{-2\Im p_{n} T})\,.
\end{equation*}
In view of the above inequality,  from \eqref{eq:delta2} it follows
\begin{multline}\label{eq:delta3}
\int_{0}^{ T-\delta} \Big| \int_0^\delta (G_1(t)-G_1(t+s))\ ds\Big|^2\ dt
\ge
2\pi(\delta-T_0/{4})^2e^{-T_0\sup|\Im p_{n}|}
(T-\delta)m_1^2
\\
\cdot\sum_{n=n_0}^\infty\ \Big(\frac{1}{\pi^2+4(T-\delta)^2(\Im p_{n})^2}-\frac{2}{(T-\delta)^2\overline{M}^2}\Big)
\ (1+e^{-2\Im p_{n}  T})|D_{n}|^2 |p_n|^4
\,.
\end{multline}
By $(T-\delta)\overline{M}=TM$ and  \eqref{eq:cpos}, for $n_0$ large enough we get  for all $n\ge n_0$
\begin{multline*}
\frac{1}{\pi^2+4(T-\delta)^2(\Im p_{n})^2}-\frac{2}{(T-\delta)^2\overline{M}^2}
\\
=\frac{1}{\pi^2+4(T-\delta)^2(\Im p_{n})^2}-\frac{2}{T^2M^2}
>\frac{1}{\pi^2+4T^2(\Im p_{n})^2}-\frac{2}{T^2M^2}>c(T,M,\varepsilon)>0 \,.
\end{multline*}
In addition, because of $\delta<T/2$, we have $\frac{T-\delta}{T}>\frac12$, so
\begin{multline*}
(T-\delta)\Big(\frac{1}{\pi^2+4(T-\delta)^2(\Im p_{n})^2}-\frac{2}{(T-\delta)^2\overline{M}^2}\Big)
\\
>T\frac{T-\delta}{T}\Big(\frac{1}{\pi^2+4T^2(\Im p_{n})^2}-\frac{2}{T^2M^2}\Big) 
>\frac{T}{2}\Big(\frac{1}{\pi^2+4T^2(\Im p_{n})^2}-\frac{2}{T^2M^2}\Big)\,.
\end{multline*}
Putting the above estimate into \eqref{eq:delta3}, we have
\begin{multline}\label{eq:delta4}
\int_{0}^{ T-\delta} \Big| \int_0^\delta (G_1(t)-G_1(t+s))\ ds\Big|^2\ dt
\\
\ge
\pi(\delta-T_0/{4})^2e^{-T_0\sup|\Im p_{n}|}
Tm_1^2\sum_{n=n_0}^\infty\ \Big(\frac{1}{\pi^2+4T^2(\Im p_{n})^2}-\frac{2}{T^2M^2}\Big)
\ (1+e^{-2\Im p_{n}  T})|D_{n}|^2 |p_n|^4
\,,
\end{multline}
and, in addition, \eqref{eq:cpos1} holds true.

On the other hand
\begin{multline*}
\int_{0}^{T-\delta} \Big| \int_0^\delta (G_1(t)-G_1(t+s))\ ds\Big|^2\ dt
\le
\delta\int_{0}^{T-\delta}  \int_0^\delta \big|G_1(t)-G_1(t+s)\big|^2\ ds\ dt
\\
\le
2\delta\int_{0}^{T-\delta}  \int_0^\delta \big(|G_1(t)|^2+|G_1(t+s)|^2\big)\ ds\ dt
\\
\le2\delta^2\int_{0}^{T} |G_1(t)|^2\ dt
+2\delta\int_0^\delta \int_{0}^{T-\delta} |G_1(t+s)|^2\ dt\ ds
\\
=
2\delta^2\int_{0}^{T} |G_1(t)|^2\ dt
+2\delta\int_0^\delta \int_{s}^{T-\delta+s} |G_1(y)|^2\ dy\ ds
\\
\le
2\delta^2\int_{0}^{T} |G_1(t)|^2\ dt
+2\delta\int_0^\delta \int_{0}^{T} |G_1(y)|^2\ dy\ ds
=
4\delta^2\int_{0}^{ T} | G_1(t)|^2\ dt\,,
\end{multline*}
whence
\begin{equation*}
\int_{0}^{ T} | G_1(t)|^2\ dt
\ge
\frac1{4\delta^2}
\int_{0}^{T-\delta} \Big| \int_0^\delta (G_1(t)-G_1(t+s))\ ds\Big|^2\ dt\,.
\end{equation*}
From the above estimate and  \eqref{eq:delta4},  it follows
\begin{equation}\label{eq:G1const}
\int_{0}^{ T} | G_1(t)|^2\ dt
\ge
2\pi C_0T
\sum_{n=n_0}^\infty\Big(\frac{1}{\pi^2+4T^2(\Im p_{n})^2}-\frac{2}{T^2M^2}\Big)\ (1+e^{-2\Im p_{n}  T})|D_{n}|^2 |p_n|^4
\,,
\end{equation}
where the constant $C_0=\frac{(\delta-T_0/{4})^2}{8\delta^2}e^{-T_0\sup|\Im p_{n}|}m_1^2$ depends on $T_0$, but not on $T$.
Moreover,
\begin{multline}\label{eq:constD}
|\D|^2
=
\frac1T\int_0^T|\D|^2\ dt
=\frac1T\int_0^T
\left|G_1(t)-
e^{\eta t}\sum_{n=n_0}^{\infty}\Big(d_nD_{n}e^{i p_{n}t}
+\overline{d_{n}}\overline{D_{n}}e^{-i \overline{p_n}t}\Big)\right|^2\ dt
\\
\le\frac2T\bigg(\int_0^T|G_1(t)|^2\ dt
+\int_0^T
\Big|e^{\eta t}\sum_{n=n_0}^{\infty}\big(d_nD_{n}e^{i p_{n}t}
+\overline{d_{n}}\overline{D_{n}}e^{-i \overline{p_n}t}\big)\Big|^2\ dt\bigg)\,.
\end{multline}
By \eqref{eq:dirDn1}, \eqref{eq:G1const} and \eqref{eq:cpos1} we have
\begin{multline*}
\int_{0}^{T}\Big| e^{\eta t}\sum_{n=n_0}^{\infty}\big(d_nD_{n}e^{i p_{n}t}
+\overline{d_{n}}\overline{D_{n}}e^{-i \overline{p_n}t}\big)
\Big|^2\ dt
\\
\le
e^{2\eta T}\int_{0}^{T}\Big|\sum_{n=n_0}^{\infty}\big(d_nD_{n}e^{i p_{n}t}
+\overline{d_{n}}\overline{D_{n}}e^{-i \overline{p_n}t}\big)
\Big|^2\ dt
\\
\le
e^{2\eta T}c(T)\sum_{n=n_0}^\infty\ |D_{n}|^2 |p_n|^4
\le
c(T)\int_{0}^{ T} | G_1(t)|^2\ dt
\,.
\end{multline*}
Plugging the above estimate into \eqref{eq:constD}, we obtain
\begin{equation*}
\int_{0}^{ T} | G_1(t)|^2\ dt\ge c(T) |\D|^2\,.
\end{equation*}
Finally, because of
\eqref{eq:G1const}
we get
\begin{equation*}
\int_{0}^{ T} | G_1(t)|^2\ dt
\ge
\pi C_0T
\sum_{n=n_0}^\infty\Big(\frac{1}{\pi^2+4T^2(\Im p_{n})^2}-\frac{2}{T^2M^2}\Big)\ (1+e^{-2\Im p_{n}  T})|D_{n}|^2 |p_n|^4
+c(T) |\D|^2
\,,
\end{equation*}
that is \eqref{eq:dn+D}.
\end{Proof}



Now, we are able to prove an inverse  inequality in the special case when our series have a finite number of terms vanishing.

\begin{theorem}\label{th:Inversa}
Under assumptions {\rm (\ref{eq:hp1})--(\ref{eq:alldn})}, for any $\varepsilon\in (0,1)$ and 
$T>\frac{2\pi}{\ga(1-\varepsilon)}$ there exist $n_0=n_0(\varepsilon)\in\N$ and $c(T,\varepsilon)>0$ such that if $C_n=D_n=0$ for any
$n< n_0$, then  we have
\begin{multline}\label{eq:Inversa}
\int_{0}^{T} \big(|u_1(t)|^2+|u_2(t)|^2\big)\ dt
\\
\ge
c(T,\varepsilon)\bigg(\sum_{n= n_0}^\infty(1+e^{-2({\Im}\om_n-\alpha)T})|C_n|^2+
\sum_{n=n_0}^\infty\ (1+e^{-2\Im p_{n} T})|D_{n}|^2|p_n|^4+|\D|^2 \bigg)
\,.
\end{multline}
\end{theorem}
\begin{Proof}
As a consequence of \eqref{eq:CnDn}, \eqref{eq:dn+D} and \eqref{eq:cpos1} with  $ T_0=\frac{2\pi}{\gamma}$ and
 $M=\gamma>\frac{2\pi}{T(1-\varepsilon)}$
we have, for $n_0$ large enough, 
\begin{multline}\label{eq:FINAL0}
\int_{0}^{T} (|u_1(t)|^2+|e^{\eta t}u_2(t)|^2)\ dt
\ge
c(T,\varepsilon)\Big(\sum_{n= n_0}^\infty(1+e^{-2({\Im}\om_n-\alpha)T})|C_n|^2+|\D|^2 \Big)
\\
+ C_0\pi T
\sum_{n=n_0}^\infty\ \Big(\frac{1}{\pi^2+4T^2(\Im p_{n})^2}-\frac{2}{T^2\ga^2}\Big)(1+e^{-2\Im p_{n} T})|D_{n}|^2 |p_n|^4
\\
-2 \pi T\sum_{n=n_0}^\infty\ \Big( \frac{1}{\pi^2+4T^2(\Im p_{n})^2}+\frac{2}{ T^2\gamma^2}\Big)(1+e^{-2{\Im}p_n T})|D_{n}|^2\,,
\end{multline}
with
\begin{equation}\label{eq:posterm}
\frac{1}{\pi^2+4T^2(\Im p_{n})^2}-\frac{2}{T^2\ga^2}>c(T,\gamma,\varepsilon)>0\,, \qquad\forall n\ge n_0\,.
\end{equation}
Now, we evaluate the sum
\begin{multline}\label{eq:FINAL}
C_0
\sum_{n=n_0}^\infty\ \Big(\frac{1}{\pi^2+4T^2(\Im p_{n})^2}-\frac{2}{T^2\ga^2}\Big)(1+e^{-2\Im p_{n} T})|D_{n}|^2 |p_n|^4
\\
-2\sum_{n=n_0}^\infty\ \Big( \frac{1}{\pi^2+4T^2(\Im p_{n})^2}+\frac{2}{ T^2\gamma^2}\Big)(1+e^{-2{\Im}p_n T})|D_{n}|^2
\\
=
\sum_{n=n_0}^\infty\ \Big(\frac{1}{\pi^2+4T^2(\Im p_{n})^2}-\frac{2}{T^2\ga^2}\Big)(1+e^{-2\Im p_{n} T})|D_{n}|^2
\left[C_0 |p_n|^4
-2\frac{\frac{1}{\pi^2+4T^2(\Im p_{n})^2}+\frac{2}{T^2\ga^2}}{\frac{1}{\pi^2+4T^2(\Im p_{n})^2}-\frac{2}{T^2\ga^2}}
\right]
\,.
\end{multline}
We note that
\begin{equation}\label{eq:incr}
\frac{\frac{1}{\pi^2+4T^2(\Im p_{n})^2}+\frac{2}{T^2\ga^2}}{\frac{1}{\pi^2+4T^2(\Im p_{n})^2}-\frac{2}{T^2\ga^2}}
=
\frac{1+\frac{2(\pi^2+4T^2(\Im p_{n})^2)}{T^2\ga^2}}{1-\frac{2(\pi^2+4T^2(\Im p_{n})^2)}{T^2\ga^2}}
=
\frac{1+\frac{2\pi^2}{T^2\ga^2}+\frac{8(\Im p_{n})^2}{\ga^2}}{1-\Big(\frac{2\pi^2}{T^2\ga^2}+\frac{8(\Im p_{n})^2}{\ga^2}\Big)}\,.
\end{equation}
Since
\begin{equation*}
\frac{2\pi^2}{T^2\ga^2}<\frac{(1-\varepsilon)^2}{2}\,,
\end{equation*}
and for $n_0$ sufficiently large 
\begin{equation*}
|\Im p_n|<\frac{1-\varepsilon}{4}\ga\,,
\qquad
\forall n\ge n_0\,,
\end{equation*}
we have
\begin{equation*}
\frac{2\pi^2}{T^2\ga^2}+\frac{8(\Im p_{n})^2}{\ga^2}<\frac{(1-\varepsilon)^2}{2}+\frac{(1-\varepsilon)^2}{2}
= (1-\varepsilon)^2<1-\varepsilon\,.
\end{equation*}
Therefore, taking into account that the function $x\to\frac{1+x}{1-x}$ is strictly  increasing on $(-\infty,1)$, 
from 
\eqref{eq:incr} we get 
\begin{equation*}
\frac{\frac{1}{\pi^2+4T^2(\Im p_{n})^2}+\frac{2}{T^2\ga^2}}{\frac{1}{\pi^2+4T^2(\Im p_{n})^2}-\frac{2}{T^2\ga^2}}
<\frac{2-\varepsilon}{\varepsilon}\,,
\end{equation*}
whence
\begin{equation*}
C_0 |p_n|^4
-2\frac{\frac{1}{\pi^2+4T^2(\Im p_{n})^2}+\frac{2}{T^2\ga^2}}{\frac{1}{\pi^2+4T^2(\Im p_{n})^2}-\frac{2}{T^2\ga^2}}
>C_0 |p_n|^4-2\frac{2-\varepsilon}{\varepsilon}\,.
\end{equation*}
Recalling that the constant $C_0$ is independent of $T$,  we can take $n_0\in\N$ (independent of $T$), large enough, so that it holds
\begin{equation*}
\frac{C_0}{2}|p_n|^4>2\frac{2-\varepsilon}{\varepsilon}
\qquad\forall n\ge n_0\,,
\end{equation*}
and hence
\begin{equation*}
C_0 |p_n|^4
-2\frac{\frac{1}{\pi^2+4T^2(\Im p_{n})^2}+\frac{2}{T^2\ga^2}}{\frac{1}{\pi^2+4T^2(\Im p_{n})^2}-\frac{2}{T^2\ga^2}}
>\frac{C_0}{2}|p_n|^4\,.
\end{equation*}
Plugging the above formula into \eqref{eq:FINAL}, we obtain
\begin{multline*}
C_0
\sum_{n=n_0}^\infty\ \Big(\frac{1}{\pi^2+4T^2(\Im p_{n})^2}-\frac{2}{T^2\ga^2}\Big)(1+e^{-2\Im p_{n} T})|D_{n}|^2 |p_n|^4
\\
-2\sum_{n=n_0}^\infty\ \Big( \frac{1}{\pi^2+4T^2(\Im p_{n})^2}+\frac{2}{ T^2\gamma^2}\Big)(1+e^{-2{\Im}p_n T})|D_{n}|^2
\\
>
\frac{C_0}2
\sum_{n=n_0}^\infty\ \Big(\frac{1}{\pi^2+4T^2(\Im p_{n})^2}-\frac{2}{T^2\ga^2}\Big)(1+e^{-2\Im p_{n} T})|D_{n}|^2
|p_n|^4\,.
\end{multline*}
Finally, from \eqref{eq:FINAL0}
it follows
\begin{multline*}
\int_{0}^{T} (|u_1(t)|^2+|e^{\eta t}u_2(t)|^2)\ dt
\\
\ge
c(T,\varepsilon)\bigg(\sum_{n= n_0}^\infty(1+e^{-2({\Im}\om_n-\alpha)T})|C_n|^2+|\D|^2 \bigg)
\\
+T\pi\frac{C_0}2\sum_{n=n_0}^\infty\ \Big(\frac{1}{\pi^2+4T^2(\Im p_{n})^2}-\frac{2}{T^2\ga^2}\Big)(1+e^{-2\Im p_{n} T})|D_{n}|^2
|p_n|^4
\,,
\end{multline*}
and hence, in view of \eqref{eq:posterm},
we obtain
\begin{multline*}
\int_{0}^{T} (|u_1(t)|^2+|e^{\eta t}u_2(t)|^2)\ dt
\\
\ge
c(T,\varepsilon)\bigg(\sum_{n= n_0}^\infty(1+e^{-2({\Im}\om_n-\alpha)T})|C_n|^2
+\sum_{n=n_0}^\infty\ (1+e^{-2\Im p_{n} T})|D_{n}|^2|p_n|^4+|\D|^2 \bigg)
\,.
\end{multline*}
In conclusion,
\begin{multline*}
\int_{0}^{T} (|u_1(t)|^2+|u_2(t)|^2)\ dt
\ge
e^{-2\eta T}\int_{0}^{T} (|u_1(t)|^2+|e^{\eta t}u_2(t)|^2)\ dt
\\
\ge
e^{-2\eta T}c(T,\varepsilon)\bigg(\sum_{n= n_0}^\infty(1+e^{-2({\Im}\om_n-\alpha)T})|C_n|^2
+\sum_{n=n_0}^\infty\ (1+e^{-2\Im p_{n} T})|D_{n}|^2|p_n|^4+|\D|^2 \bigg)
\,,
\end{multline*}
that is
\eqref{eq:Inversa}.
\end{Proof}

As regards the direct inequality, first we recall the following result, see \cite[Theorem $4.2$]{LoretiSforza1}. 
\begin{theorem}\label{th:dir.ingham}
Under assumptions {\rm (\ref{eq:hom1})--(\ref{eq:hom3})}, for any $\varepsilon\in(0,1)$ and 
$T>\frac{\pi}{\ga\sqrt{1-\varepsilon}}$ there exist $c(T)>0$ and $n_0=n_0(\varepsilon)\in\N$ such that if $C_n=0$ for $n< n_0$,
then we have
\begin{equation}\label{eq:dir-cosine}
 \int_{-\infty}^{\infty} k^*(t)\Big|\sum_{n=n_0}^\infty R_{n}e^{r_{n}t}+C_{n}e^{i\omega_{n}t}+\overline{C_n}e^{-i\overline{\omega_n}t}\Big|^2\ dt
 \le c(T)\sum_{n= n_0}^\infty|C_n|^2\,.
\end{equation}
\end{theorem}
\begin{theorem}\label{th:dirDn}
Assume {\rm (\ref{eq:hp1})}, {\rm (\ref{eq:hp2})} and {\rm (\ref{eq:alldn})}.
Let $T>0$, $\varepsilon\in (0,1)$ and 
$M>\frac{\pi}{T\sqrt{1-\varepsilon}}$.  There exist $c(T)>0$ and  $n_0=n_0(\varepsilon)\in\N$  such that if $D_n=0$ for any
$n< n_0$, then  we have
\begin{equation}\label{eq:dirDn}
\int_{-\infty}^{\infty} k^*(t)\bigg| \sum_{n=n_0}^{\infty}d_nD_{n}e^{i p_{n}t}
+\overline{d_{n}}\overline{D_{n}}e^{-i \overline{p_n}t}
+\D e^{-\eta t}\bigg|^2\ dt
\le
c(T)\bigg(\sum_{n=n_0}^\infty\ |D_{n}|^2 |p_n|^4+|\D|^2\bigg)
\,.
\end{equation}
\end{theorem}
\begin{Proof}
Since the function $k^*(t)$ is positive, we have
\begin{multline}\label{eq:dir2}
\int_{-\infty}^{\infty} k^*(t)\bigg| \sum_{n=n_0}^{\infty}d_nD_{n}e^{i p_{n}t}
+\overline{d_{n}}\overline{D_{n}}e^{-i \overline{p_n}t}
+\D e^{-\eta t}\bigg|^2\ dt
\\
\le
2\int_{-\infty}^{\infty} k^*(t)\bigg| \sum_{n=n_0}^{\infty}d_nD_{n}e^{i p_{n}t}
+\overline{d_{n}}\overline{D_{n}}e^{-i \overline{p_n}t}\bigg|^2\ dt
+2\int_{-\infty}^{\infty} k^*(t)e^{-2\eta t}\ dt\
|\D|^2
\,.
\end{multline}
We evaluate the first integral by using Proposition \ref{pr:dirDn0}: indeed, 
from \eqref{eq:dirDn0} we get
\begin{equation*}
\int_{-\infty}^{\infty} k^*(t)\bigg| \sum_{n=n_0}^{\infty}d_nD_{n}e^{i p_{n}t}
+\overline{d_{n}}\overline{D_{n}}e^{-i \overline{p_n}t}\bigg|^2\ dt
\le
c(T)
\sum_{n=n_0}^\infty |D_{n}|^2  |p_n|^4
\,.
\end{equation*}
Putting the previous estimate in \eqref{eq:dir2}, we obtain
\begin{multline}\label{eq:dir2a}
\int_{-\infty}^{\infty} k^*(t)\bigg| \sum_{n=n_0}^{\infty}d_nD_{n}e^{i p_{n}t}
+\overline{d_{n}}\overline{D_{n}}e^{-i \overline{p_n}t}
+\D e^{-\eta t}\bigg|^2\ dt
\\
\le
2c(T)
\sum_{n=n_0}^\infty |D_{n}|^2  |p_n|^4
+2\int_{-\infty}^{\infty} k^*(t)e^{-2\eta t}\ dt\
|\D|^2
\,.
\end{multline}
In addition, formula \eqref{eqn:k1}  yields
\begin{equation*}
\int_{-\infty}^{\infty} k^*(t)e^{-2\eta t}\ dt=
\cosh(2\eta T)\frac{4T\pi}{\pi^2+16T^2\eta^2}\,.
\end{equation*}
Because of the above formula from \eqref{eq:dir2a} it follows
\begin{multline*}
\int_{-\infty}^{\infty} k^*(t)\bigg| \sum_{n=n_0}^{\infty}d_nD_{n}e^{i p_{n}t}
+\overline{d_{n}}\overline{D_{n}}e^{-i \overline{p_n}t}
+\D e^{-\eta t}\bigg|^2\ dt
\\
\le
2c(T)
\sum_{n=n_0}^\infty |D_{n}|^2  |p_n|^4
+\cosh(2\eta T)\frac{8T\pi}{\pi^2+16T^2\eta^2}\
|\D|^2
\,,
\end{multline*}
that is \eqref{eq:dirDn}  .
\end{Proof}

Finally, thanks to Theorems \ref{th:dir.ingham} and \ref{th:dirDn} we are able to prove an Ingham type direct estimate for the solution 
$(u_1,u_2)$
of coupled systems in the special case when our series have a finite number of terms vanishing.
\begin{theorem}\label{th:Diretta}
Under assumptions {\rm (\ref{eq:hp1})--(\ref{eq:alldn})}, for any $\varepsilon\in (0,1)$ and 
$T>\frac{\pi}{\ga\sqrt{1-\varepsilon}}$ there exist $c(T)>0$ and $n_0=n_0(\varepsilon)\in\N$ such that if $C_n=D_n=0$ for any
$n< n_0$, then  we have
\begin{equation}\label{eq:Diretta}
\int_{-T}^{T} \big(|u_1(t)|^2+|u_2(t)|^2\big)\ dt
\le
c(T)\bigg(\sum_{n= n_0}^\infty |C_n|^2+
\sum_{n=n_0}^\infty\ |D_{n}|^2|p_n|^4+|\D|^2 \bigg)
\,.
\end{equation}
\end{theorem}
\begin{Proof}
First of all, since the function $k^*(t)$ is positive, we can write
\begin{multline*}
\int_{-\infty}^{\infty} k^*(t)\bigg|\sum_{n=n_0}^\infty R_{n}e^{r_{n}t}+C_{n}e^{i\omega_{n}t}+\overline{C_n}e^{-i\overline{\omega_n}t}+D_{n}e^{i p_{n}t}+\overline{D_{n}}e^{-i\overline{ p_n}t}\bigg|^2\ dt
\\
\le
2\int_{-\infty}^{\infty} k^*(t)\bigg|\sum_{n=n_0}^\infty R_{n}e^{r_{n}t}+C_{n}e^{i\omega_{n}t}
+\overline{C_n}e^{-i\overline{\omega_n}t}\bigg|^2\ dt
\\
+2\int_{-\infty}^{\infty} k^*(t)\bigg|\sum_{n=n_0}^\infty D_{n}e^{i p_{n}t}+\overline{D_{n}}e^{-i\overline{ p_n}t}\bigg|^2\ dt\,.
\end{multline*}
So, we can apply Theorem \ref{th:dir.ingham}: plugging into the above formula the inequality \eqref{eq:dir-cosine}, we obtain 
\begin{multline*}
\int_{-\infty}^{\infty} k^*(t)\bigg|\sum_{n=n_0}^\infty R_{n}e^{r_{n}t}+C_{n}e^{i\omega_{n}t}+\overline{C_n}e^{-i\overline{\omega_n}t}+D_{n}e^{i p_{n}t}+\overline{D_{n}}e^{-i\overline{ p_n}t}\bigg|^2\ dt
\\
\le
c(T)\sum_{n= n_0}^\infty|C_n|^2+2\int_{-\infty}^{\infty} k^*(t)
\bigg|\sum_{n=n_0}^\infty D_{n}e^{i p_{n}t}+\overline{D_{n}}e^{-i\overline{ p_n}t}\bigg|^2\ dt\,.
\end{multline*}
In Proposition \ref{pr:EnD} we can take $E_n=D_n$ and $M=\gamma$, so by the previous inequality and \eqref{eq:EnD2} we get
\begin{multline*}
\int_{-\infty}^{\infty} k^*(t)\bigg|\sum_{n=n_0}^\infty R_{n}e^{r_{n}t}+C_{n}e^{i\omega_{n}t}+\overline{C_n}e^{-i\overline{\omega_n}t}+D_{n}e^{i p_{n}t}+\overline{D_{n}}e^{-i\overline{ p_n}t}\bigg|^2\ dt
\\
\le
c(T)\bigg(\sum_{n= n_0}^\infty|C_n|^2+\sum_{n= n_0}^\infty|D_n|^2\bigg)\,.
\end{multline*}
Moreover, by the above estimate and \eqref{eq:dirDn} we obtain
\begin{equation*}
\int_{-\infty}^{\infty} k^*(t) \big(|u_1(t)|^2+|u_2(t)|^2\big)\ dt
\le
c(T)\bigg(\sum_{n= n_0}^\infty |C_n|^2+
\sum_{n=n_0}^\infty\ |D_{n}|^2|p_n|^4 +|\D|^2\bigg)
\,.
\end{equation*}
Now, if we consider the last inequality with the function $k^*$ 
replaced by the analogous one relative to $2T$ instead of $T$,
see (\ref{eq:kcos}),
then 
we get 
\begin{equation*}
\int_{-2T}^{2T} \cos \frac{\pi t}{4T}\big(|u_1(t)|^2+|u_2(t)|^2\big)\ dt
\le
c(2T)\bigg(\sum_{n= n_0}^\infty |C_n|^2+
\sum_{n=n_0}^\infty\ |D_{n}|^2|p_n|^4 +|\D|^2\bigg),
\end{equation*}
whence, thanks to
$
\cos \frac{\pi t}{4T}\ge \frac1{\sqrt2}
$
for $|t|\le T$,
it follows
\begin{equation*}
\int_{-T}^{T} \big(|u_1(t)|^2+|u_2(t)|^2\big)\ dt
\le
\sqrt2c(2T)\bigg(\sum_{n= n_0}^\infty |C_n|^2+
\sum_{n=n_0}^\infty\ |D_{n}|^2|p_n|^4 +|\D|^2\bigg)\,.
\end{equation*}
So, the proof of \eqref{eq:Diretta} is complete.
\end{Proof}

\subsection{Haraux type estimates}

To prove our results, we need to introduce 
a suitable family of operators which annihilate a finite number of terms  in the Fourier series.
For the reader's convenience, we proceed to recall the definition of
operators, which was  given in \cite{LoretiSforza1}
and is slightly different from those introduced in
\cite{Ha} and
\cite{KL1}.

Given $\delta >0$ and $z\in\C$ arbitrarily, we define the linear operator $I_{\delta,z}$ as follows: for every continuous function $u:\R\to\C$ the function
$I_{\delta,z}u:\R\to\C$ is given by the formula
\begin{equation}\label{eq:defI}
I_{\delta,z}u(t):=u(t)-\frac 1\delta\int_0^\delta e^{-iz s} u(t+s)\ ds\,,\qquad t\in\R\,.
\end{equation}
A list of properties connected with operators  $I_{\delta,z}$ is now in order.

\begin{lemma} \label{le:opI1}
For any $\delta >0$ and $z\in\C$ the following statements hold true.

\begin{itemize}
\item[(i)] $I_{\delta,z}(e^{iz t})= 0 $\,.
\item[(ii)] For any $z'\in\C$, $z'\not=z$, we have
\begin{equation*}
I_{\delta,z}(e^{iz' t})=
\bigg(1-\frac{e^{i(z'-z) \delta}-1}{i(z'-z) \delta}\bigg)e^{iz' t}
\,.
\end{equation*}
\item[(iii)] The linear operators $I_{\delta,z}$ commute, that is,
for any $\delta'>0$, $z'\in\C$ and continuous function $u:\R\to\C$ we have
\begin{equation*}
I_{\delta,z}I_{\delta',z'}u=I_{\delta',z'}I_{\delta,z}u\,.
\end{equation*}
\item[(iv)]
For any $T>0$ and  continuous function $u:\R\to\C$ we have
\begin{equation}\label{eq:boundedI}
\int_0^T|I_{\delta,z}u(t)|^2\ dt\le
2(1+e^{2|{\Im}z| \delta})\int_0^{T+\delta}|u(t)|^2\ dt\,.
\end{equation}
\end{itemize}
\end{lemma}
We now  define another operator
\begin{equation}\label{eq:defI1}
I_{\delta,r,\omega, p}:= I_{\delta, -ir}\circ I_{\delta,\omega}\circ I_{\delta,-\overline{\omega}} \circ I_{\delta,p}\circ I_{\delta,-\overline{p}}
\qquad \delta >0\,,r\in\R\,,\,\omega\,,p\in\C\,,
\end{equation}
where the symbol $\circ$ denotes the usual composition among operators.

By using Lemma \ref{le:opI1} one can easily prove the following properties concerning operators $I_{\delta,r,\omega, p}$.
\begin{lemma} \label{le:opI2}
For any $\delta >0$ and $r\in\R\,,\omega\,,p\in\C$ the following statements hold true.
\begin{itemize}
\item[ (i)]
$
I_{\delta,r,\omega,p}(e^{r t})=I_{\delta,r,\omega,p}(e^{i\omega t})=I_{\delta,r,\omega,p}(e^{-i\overline{\omega} t})
=I_{\delta,r,\omega,p}(e^{ip t})=I_{\delta,r,\omega,p}(e^{-i\overline{p} t})=  0 \,.
$
\item[ (ii)]
For any $r'\in\R$, $r'\not\in\{r,i\om,-i\overline{\om},ip,-i\overline{p}\}$, we have
\begin{equation*}
I_{\delta,r,\omega, p}(e^{r' t})
=
\prod_{z\in\{r,i\om,-i\overline{\om},ip,-i\overline{p}\}}\bigg(1-\frac{e^{(r'-z) \delta}-1}{(r'-z) \delta}\bigg)e^{r' t}
\,.
\end{equation*}
\item[ (iii)]
For any $z'\in\C$, $z'\not\in\{-ir,\om,-\overline{\om},p,-\overline{p}\}$, we have
\begin{equation*}
I_{\delta,r,\omega, p}(e^{iz' t})
=
\prod_{z\in\{-ir,\om,-\overline{\om},p,-\overline{p}\}}\bigg(1-\frac{e^{i(z'-z) \delta}-1}{i(z'-z) \delta}\bigg)e^{iz' t}
\,.
\end{equation*}
\item[(iv)] The linear operators $I_{\delta,r,\omega, p}$ commute, that is,
for any $\delta'>0$, $r'\in\R\,,\omega'\,,p'\in\C$ and continuous function $u:\R\to\C$ we have
\begin{equation*}
I_{\delta,r,\omega, p}I_{\delta',r',\omega', p'}u=I_{\delta',r',\omega', p'}I_{\delta,r,\omega, p}u\,.
\end{equation*}
\end{itemize}
\end{lemma}

\begin{corollary}\label{co:boundedI}
For any $T>0$, $\delta>0$, $r\in\R\,,\omega\,,p\in\C$ and continuous function $u:\R\to\C$ we have
\begin{equation}\label{eq:boundedI1}
\int_0^T|I_{\delta,r,\omega, p}u(t)|^2\ dt\le
2^5(1+e^{2|r| \delta})(1+e^{2|{\Im}\om| \delta})^2(1+e^{2|{\Im}p| \delta})^2\int_0^{T+5\delta}|u(t)|^2\ dt
\,.
\end{equation}
\end{corollary}
\begin{Proof}
By applying (\ref{eq:boundedI}) repeatedly,  we obtain
\begin{multline*}
\int_0^T|I_{\delta,r,\omega, p}u(t)|^2\ dt
=\int_0^T| I_{\delta, -ir} I_{\delta,\omega} I_{\delta,-\overline{\omega}}  I_{\delta,p} I_{\delta,-\overline{p}}u(t)|^2\ dt
\\
\le
2^2(1+e^{2|{\Im}p| \delta})^2\int_0^{T+2\delta}|I_{\delta, -ir} I_{\delta,\omega} I_{\delta,-\overline{\omega}}u(t)|^2\ dt\\
\le
2^4(1+e^{2|{\Im}p| \delta})^2(1+e^{2|{\Im}\om| \delta})^2\int_0^{T+4\delta}|I_{\delta, -ir} u(t)|^2\ dt\\
\le
2^5(1+e^{2|{\Im}p| \delta})^2(1+e^{2|{\Im}\om| \delta})^2(1+e^{2|r| \delta})\int_0^{T+5\delta}|u(t)|^2\ dt
\,,
\end{multline*}
that is (\ref{eq:boundedI1}).
\end{Proof}

\begin{proposition}\label{pr:haraux-inv}
Let $\{\om_n\}_{n\in\N}$, $\{r_n\}_{n\in\N}$  and
$\{p_n\}_{n\in\N}$ be sequences of pairwise  distinct numbers
such that $\om_n\not= p_m$, $\om_n\not=\overline{p_m}$, 
$r_n\not= i\om_m$, $r_n\not= ip_m$, $r_n\not=-\eta$, $p_n\not=0$, for any $n\,,m\in\N$, 
\begin{equation}\label{eq:ha1}
\lim_{n\to\infty}|\om_n|=\lim_{n\to\infty}|p_n|=+\infty\,,
\end{equation}
and
\begin{equation}\label{eq:alldnbis}
|d_n|\ge m_1|p_n|^2\,\qquad\forall n\in\N
\qquad (m_1>0)
\,.
\end{equation}
Assume that there exists $n_0\in\N$ such that 
for any sequences $\{R_n\}$, $\{C_n\}$ and $\{D_n\}$ verifying
\begin{equation*}
R_n=C_n=D_n=0\qquad\mbox{for any}\quad n< n_0\,,
\end{equation*}
the estimates
\begin{equation}\label{eq:haraux-inv}
\int_{0}^{T} \big(|u_1(t)|^2+|u_2(t)|^2\big)\ dt
\ge
c_1\bigg(\sum_{n= n_0}^\infty|C_n|^2+
\sum_{n=n_0}^\infty\ |D_{n}|^2|p_n|^4 +|\D|^2 \bigg)\,,
\end{equation} 
\begin{equation}\label{eq:haraux-inv1}
\int_{0}^{T} \big(|u_1(t)|^2+|u_2(t)|^2\big)\ dt
\le
c_2\bigg(\sum_{n= n_0}^\infty|C_n|^2+
\sum_{n=n_0}^\infty\ |D_{n}|^2|p_n|^4 +|\D|^2 \bigg)\,,
\end{equation} 
are satisfied for some constants $c_1,c_2>0$.

Then, there exists $C_1>0$ such that for any sequences $\{R_n\}$, $\{C_n\}$ and $\{D_n\}$
the  estimate  
\begin{equation}\label{eq:haraux-inv22}
\int_{0}^{T} \big(|u_1(t)|^2+|u_2(t)|^2\big)\ dt
\ge
C_1\bigg(\sum_{n= 1}^\infty|C_n|^2+
\sum_{n=1}^\infty\ |D_{n}|^2|p_n|^4 +|\D|^2 \bigg)
\end{equation}
holds. 
\end{proposition}

\begin{Proof}
To begin with, we will transform the functions 
\begin{equation*}
u_1(t)=\sum_{n=1}^{\infty}\Big(R_ne^{r_nt}+C_ne^{i\om_nt}+\overline{C_n}e^{-i\overline{\omega_n}t}
+D_{n}e^{i p_{n}t}+\overline{D_{n}}e^{-i \overline{p_n}t}\Big)
\end{equation*}
\begin{equation*}
u_2(t)=\sum_{n=1}^{\infty}
\Big(d_nD_{n}e^{i p_{n}t}
+\overline{d_n}\overline{D_{n}}e^{-i \overline{p_n}t}\Big)
+\D e^{-\eta t}
\end{equation*}
in such a way that the
series have null   terms corresponding to indices $n=1,\cdots,n_0-1$,
because so we can apply our assumptions \eqref{eq:haraux-inv} and \eqref{eq:haraux-inv1}. 

To this end, we fix $\varepsilon>0$ and choose $\delta\in \big(0,\frac{\varepsilon}{5n_0}\big)$.  Let us denote by
$\I$ the composition of all linear operators $I_{\delta,r_j,\omega_j,p_j}$, $j=1,\cdots,n_0-1$. We note that by Lemma \ref{le:opI2}-(iv) the definition of $\I$ does not depend on the
order of the operators $I_{\delta,r_j,\omega_j,p_j}$. 

By using Lemma \ref{le:opI2}, we get
\begin{equation*}
\I u_1(t)=\sum_{n=n_0}^{\infty}\Big(R'_ne^{r_nt}+C'_ne^{i\om_nt}+\overline{C'_n}e^{-i\overline{\omega_n}t}
+D'_{n}e^{i p_{n}t}+\overline{D'_{n}}e^{-i \overline{p_n}t}\Big)
\end{equation*}
\begin{equation*}
\I u_2(t)=\sum_{n=n_0}^{\infty}
\Big(d_nD'_{n}e^{i p_{n}t}
+\overline{d_n}\overline{D'_{n}}e^{-i \overline{p_n}t}\Big)
+\D' e^{-\eta t}
\end{equation*}
where
$$
R'_n:=R_n
\prod_{j=1}^{n_0-1}
\prod_{z\in\{r_j,i\om_j,-i\overline{\om_j},ip_j,-i\overline{p_j}\}}\bigg(1-\frac{e^{(r_n-z) \delta}-1}{(r_n-z) \delta}\bigg)
\,,
$$
$$
C'_n:=C_n
\prod_{j=1}^{n_0-1}
\prod_{z\in\{-ir_j,\om_j,-\overline{\om_j},p_j,-\overline{p_j}\}}\bigg(1-\frac{e^{i(\om_n-z) \delta}-1}{i(\om_n-z) \delta}\bigg)
\,,
$$
$$
D'_n:=D_n
\prod_{j=1}^{n_0-1}
\prod_{z\in\{-ir_j,\om_j,-\overline{\om_j},p_j,-\overline{p_j}\}}\bigg(1-\frac{e^{i(p_n-z) \delta}-1}{i(p_n-z) \delta}\bigg)
\,,
$$
$$
\D':=\D\prod_{j=1}^{n_0-1}
\prod_{z\in\{r_j,i\om_j,-i\overline{\om_j},ip_j,-i\overline{p_j}\}}\bigg(1+\frac{e^{-(\eta+z) \delta}-1}{(\eta+z) \delta}\bigg)
\,.
$$
Therefore, we are in condition to
apply estimate (\ref{eq:haraux-inv}) to functions  $\I u_1(t)$ and $\I u_2(t)$:
\begin{equation}\label{eq:haraux-inv2} 
\int_{0}^{T} (|\I u_1(t)|^2+|\I u_2(t)|^2)\ dt
\ge
c_1\bigg(\sum_{n= n_0}^\infty|C'_n|^2+
\sum_{n=n_0}^\infty\ |D'_{n}|^2|p_n|^4 +|\D'|^2 \bigg)\,.
\end{equation} 
Next, we choose $\delta\in \big(0,\frac{\varepsilon}{5n_0}\big)$ such that for any $n\ge n_0$ none of the products
\begin{equation}\label{eq:products}
\prod_{j=1}^{n_0-1}
\prod_{z\in\{-ir_j,\om_j,-\overline{\om_j},p_j,-\overline{p_j}\}}\bigg(1-\frac{e^{i(\om_n-z) \delta}-1}{i(\om_n-z) \delta}\bigg)
\end{equation}
vanishes. This is possible because the analytic function 
$$w\longmapsto1-\frac{e^w-1}{w}$$ 
does not vanish identically and, since every
number $\om_n-z$, 
with $z\in\{-ir_j,\om_j,-\overline{\om_j},p_j,-\overline{p_j}\}$, is different from zero, we have to exclude only a countable set of values of $\delta$.

Then, we note that there exists a constant $c'>0$ such that for any $n\ge n_0$
\begin{equation}\label{eq:C'}
\Bigg|\prod_{j=1}^{n_0-1}
\prod_{z\in\{-ir_j,\om_j,-\overline{\om_j},p_j,-\overline{p_j}\}}\bigg(1-\frac{e^{i(\om_n-z) \delta}-1}{i(\om_n-z) \delta}\bigg)\Bigg|^2\ge c'
\,.
\end{equation}
Indeed, it is sufficient to observe that for any fixed $j=1,\cdots,n_0-1$ 
and $z\in\{-ir_j,\om_j,-\overline{\om_j},p_j,-\overline{p_j}\}$ we have
$$
\Big|\frac{e^{i(\om_n-z) \delta}-1}{i(\om_n-z) \delta}\Big|\le\frac{e^{-{\Im}(\om_n-z) \delta}+1}{|\om_n-z| \delta}
\to 0\qquad
\mbox{as}\quad n\to \infty\,,
$$
thanks to  \eqref{eq:ha1}. As a result, the product in \eqref{eq:products}
tends to $1$ as $n\to \infty$, so that it is minorized, e.g., by $1/2$ for  $n$ large enough. 
By repeating the same argumentations used to get \eqref{eq:C'}, we also have 
\begin{equation}\label{eq:D'}
\Bigg|\prod_{j=1}^{n_0-1}
\prod_{z\in\{-ir_j,\om_j,-\overline{\om_j},p_j,-\overline{p_j}\}}\bigg(1-\frac{e^{i(p_n-z) \delta}-1}{i(p_n-z) \delta}\bigg)\Bigg|^2\ge c'
\,.
\end{equation}
In addition, we can assume that
$$
\Bigg|\prod_{j=1}^{n_0-1}
\prod_{z\in\{r_j,i\om_j,-i\overline{\om_j},ip_j,-i\overline{p_j}\}}\bigg(1+\frac{e^{-(\eta+z) \delta}-1}{(\eta+z) \delta}\bigg)\Bigg|^2\ge c'\,,
$$
and hence
$$
|\D'|^2\ge c'|\D|^2\,.
$$
Therefore, the above estimate and
\eqref{eq:haraux-inv2}--\eqref{eq:D'} 
 yield
\begin{equation}\label{eq:haraux-inv3}
\int_{0}^{T} \big(|\I u_1(t)|^2+|\I u_2(t)|^2\big)\ dt
\ge
c'c_1\bigg(\sum_{n= n_0}^\infty|C_n|^2+
\sum_{n=n_0}^\infty\ |D_{n}|^2|p_n|^4 +|\D|^2 \bigg)\,.
\end{equation} 
On the other hand, applying  (\ref{eq:boundedI1}) repeatedly with $r=r_j$, $\om=\om_j$ and $p=p_j$, $j=1,\cdots,n_0-1$, we have
\begin{multline*}
\int_{0}^{T} \big(|\I u_1(t)|^2+|\I u_2(t)|^2\big)\ dt
\\
\le
2^{5(n_0-1)}\prod_{j=1}^{n_0-1}(1+e^{2|r_j| \delta})(1+e^{2|{\Im}\om_j| \delta})^2(1+e^{2|{\Im}p_j| \delta})^2
\int_0^{T+5(n_0-1)\delta}(| u_1(t)|^2+| u_2(t)|^2)\ dt 
\,.
\end{multline*}
From the above inequality,
by using \eqref{eq:haraux-inv3} and $5n_0\delta<\varepsilon$, it follows 
\begin{multline*}
\sum_{n= n_0}^\infty |C_n|^2+
\sum_{n=n_0}^\infty\ |D_{n}|^2|p_n|^4 +|\D|^2
\\
\le
\frac{2^{5(n_0-1)}}{c'c_1}\prod_{j=1}^{n_0-1}(1+e^{|r_j| \varepsilon/n_0})(1+e^{|{\Im}\om_j| \varepsilon/n_0})^2
(1+e^{|{\Im}p_j| \varepsilon/n_0})^2
\int_0^{T+\varepsilon}\big(| u_1(t)|^2+| u_2(t)|^2\big)\ dt 
\,,
\end{multline*}
whence, passing to the limit as $\varepsilon\to 0^+$, we have
\begin{equation}\label{eq:C''}
\sum_{n= n_0}^\infty|C_n|^2+
\sum_{n=n_0}^\infty\ |D_{n}|^2|p_n|^4 +|\D|^2
\le
\frac{2^{10(n_0-1)}}{c'c_1}\int_0^{T}\big(| u_1(t)|^2+| u_2(t)|^2\big)\ dt  
\,.
\end{equation}
Moreover, thanks to the triangle inequality, we get
\begin{multline}\label{eq:n0u1}
\int_0^T\bigg|\sum_{n=1}^{n_0-1}
R_ne^{r_nt}+C_ne^{i\om_nt}+\overline{C_n}e^{-i\overline{\omega_n}t}
+D_{n}e^{i p_{n}t}+\overline{D_{n}}e^{-i \overline{p_n}t}
\bigg|^2 \ d t
\\
=
\int_0^T\bigg|u_1(t)-\sum_{n=n_0}^{\infty}\Big(R_ne^{r_nt}+C_ne^{i\om_nt}+\overline{C_n}e^{-i\overline{\omega_n}t}
+D_{n}e^{i p_{n}t}+\overline{D_{n}}e^{-i \overline{p_n}t}\Big)\bigg|^2 d t\\
\le
2\int_0^T |u_1(t)|^2 d t+2\int_0^T\Big|\sum_{n=n_0}^{\infty}R_ne^{r_nt}+C_ne^{i\om_nt}+\overline{C_n}e^{-i\overline{\omega_n}t}
+D_{n}e^{i p_{n}t}+\overline{D_{n}}e^{-i \overline{p_n}t}\bigg|^2 d t
\end{multline}
and
\begin{multline}\label{eq:n0u2}
\int_0^T\bigg|\sum_{n=1}^{n_0-1}
d_nD_{n}e^{i p_{n}t}+\overline{d_n}\overline{D_{n}}e^{-i \overline{p_n}t}\bigg|^2 \ d t
\\
=
\int_0^T\bigg|u_2(t)-\sum_{n=n_0}^{\infty}
\Big(d_nD_{n}e^{i p_{n}t}+\overline{d_n}\overline{D_{n}}e^{-i \overline{p_n}t}\Big)
-\D e^{-\eta t}\bigg|^2 d t\\
\le
2\int_0^T |u_2(t)|^2 d t+2\int_0^T\Big|\sum_{n=n_0}^{\infty}\Big(d_nD_{n}e^{i p_{n}t}+\overline{d_n}\overline{D_{n}}e^{-i \overline{p_n}t}\Big)
+\D e^{-\eta t}\bigg|^2 d t
\,.
\end{multline}
Putting together \eqref{eq:n0u1} and \eqref{eq:n0u2} and using
\eqref{eq:haraux-inv1} and \eqref{eq:C''} we have
\begin{multline}\label{eq:1+C2C''}
\int_0^T\Big|\sum_{n=1}^{n_0-1}
R_ne^{r_nt}+C_ne^{i\om_nt}+\overline{C_n}e^{-i\overline{\omega_n}t}
+D_{n}e^{i p_{n}t}+\overline{D_{n}}e^{-i \overline{p_n}t}
\Big|^2 \ dt
\\
+\int_0^T\Big|\sum_{n=1}^{n_0-1}
d_nD_{n}e^{i p_{n}t}
+\overline{d_n}\overline{D_{n}}e^{-i \overline{p_n}t}\Big|^2 d t
\\
\le
2\int_0^{T}\big(| u_1(t)|^2+| u_2(t)|^2\big)\ dt 
+2c_2\bigg(\sum_{n= n_0}^\infty|C_n|^2+
\sum_{n=n_0}^\infty\ |D_{n}|^2|p_n|^4 +|\D|^2 \bigg)
\\
\le
2\Big(1+c_2\frac{2^{10(n_0-1)}}{c'c_1}\Big)\int_0^{T}\big(| u_1(t)|^2+| u_2(t)|^2\big)\ dt \,.
\end{multline}
Let us note that the expression
\begin{multline*}
\int_0^T\Big|\sum_{n=1}^{n_0-1}
R_ne^{r_nt}+C_ne^{i\om_nt}+\overline{C_n}e^{-i\overline{\omega_n}t}
+D_{n}e^{i p_{n}t}+\overline{D_{n}}e^{-i \overline{p_n}t}
\Big|^2 \ dt
\\
+\int_0^T\Big|\sum_{n=1}^{n_0-1}
d_nD_{n}e^{i p_{n}t}
+\overline{d_n}\overline{D_{n}}e^{-i \overline{p_n}t}\Big|^2 d t
\end{multline*}
is a positive semidefinite quadratic form of the variable 
$$
\big(\{R_n\}_{n<n_0},\{C_n\}_{n<n_0},\{d_nD_n\}_ {n<n_0}\big)\in\R^{n_0-1}\times\C^{n_0-1}\times \C^{n_0-1}\,.
$$
Moreover, it is positive {\it definite}, because the functions $e^{r_nt}$, $e^{i\om_nt}$, $e^{ip_nt}$, ${n<n_0}$, are linearly independent. Hence, there exists a constant
$c''>0$ such that
\begin{multline*}
\int_0^T\Big|\sum_{n=1}^{n_0-1}
R_ne^{r_nt}+C_ne^{i\om_nt}+\overline{C_n}e^{-i\overline{\omega_n}t}
+D_{n}e^{i p_{n}t}+\overline{D_{n}}e^{-i \overline{p_n}t}
\Big|^2 \ dt
\\
+\int_0^T\Big|\sum_{n=1}^{n_0-1}
d_nD_{n}e^{i p_{n}t}
+\overline{d_n}\overline{D_{n}}e^{-i \overline{p_n}t}\Big|^2 d t
\ge c''\sum_{n=1}^{n_0-1}\Big(|R_n|^2+|C_n|^2+|D_{n}|^2|p_n|^4\Big)\,,
\end{multline*}
taking into account \eqref{eq:alldnbis}.
So, from \eqref{eq:1+C2C''} and the above inequality  we deduce that
\begin{eqnarray*}
\sum_{n=1}^{n_0-1}\Big(|C_n|^2+|D_{n}|^2|p_n|^4\Big)\le 
\frac{2}{c''}\Big(1+c_2\frac{2^{10(n_0-1)}}{c'c_1}\Big)\int_0^{T}\big(| u_1(t)|^2+| u_2(t)|^2\big)\ dt
\,.
\end{eqnarray*}
Finally, the above estimate and \eqref{eq:C''} yield the required inequality (\ref{eq:haraux-inv22}). 
\end{Proof}

\subsection{Inverse and direct inequalities}

We recall that
\begin{equation*}
u_1(t)=\sum_{n=1}^{\infty}\Big(R_ne^{r_nt}+C_ne^{i\om_nt}+\overline{C_n}e^{-i\overline{\omega_n}t}
+D_{n}e^{i p_{n}t}+\overline{D_{n}}e^{-i \overline{p_n}t}\Big)
\,,
\end{equation*}
\begin{equation*}
u_2(t)=\sum_{n=1}^{\infty}
\Big(d_nD_{n}e^{i p_{n}t}
+\overline{d_n}\overline{D_{n}}e^{-i \overline{p_n}t}\Big)
+\D e^{-\eta t}
\,,
\end{equation*}
with
\begin{equation}\label{eq:exprD}
\D=-\frac\beta A\sum_{n=1}^{\infty}\Re p_{n}
\Big(\frac{D_{n}}{\eta+ip_n}+\frac{\overline{D_{n}}}{\eta-ip_n}
\Big)\,.
\end{equation}

\begin{theorem}\label{th:inv.ingham1}
Let $\{\om_n\}_{n\in\N}$, $\{r_n\}_{n\in\N}$  and
$\{p_n\}_{n\in\N}$ be sequences of pairwise  distinct numbers
such that $\om_n\not= p_m$, $\om_n\not=\overline{p_m}$, 
$r_n\not= i\om_m$, $r_n\not= ip_m$, $r_n\not=-\eta$, $p_n\not=0$, for any $n\,,m\in\N$.
Assume
\begin{equation*}
\lim_{n\to\infty}(\Re p_{n+1}-\Re p_{n})=+\infty\,,
\end{equation*}
\begin{equation*}
\lim_{n\to\infty}\Im p_{n}=0\,,
\end{equation*}
and for some $\gamma>0$, $\alpha\in\R$, $n'\in\N$, $\mu>0$, $\nu> 1/2$, $m_1\,,m_2>0$
\begin{equation*}
\liminf_{n\to\infty}({\Re}\om_{n+1}-{\Re}\om_{n})=\gamma\,,
\end{equation*}
\begin{equation*}
\lim_{n\to\infty}{\Im}\om_n=\alpha
\,,
\qquad
r_n\le -{\Im}\om_n\,\qquad\forall\ n\ge n'\,,
\end{equation*}
\begin{equation*}
|R_n|\le \frac{\mu}{n^{\nu}}|C_n|\,\quad\forall\ n\ge n'\,,
\qquad
|R_n|\le \mu|C_n|\,\quad\forall\ n\le n'\,,
\end{equation*}
\begin{equation*}
m_1|p_n|^2\le|d_n|\le m_2|p_n|^2\,\qquad\forall n\in\N\,.
\end{equation*}
Then, for any
$T>2\pi/\ga$ we have 
 \begin{equation}\label{eq:inv.ingham}
 \int_{0}^{T} \big(|u_1(t)|^2+|u_2(t)|^2\big)\ dt
 \ge
 c_1(T)\bigg(\sum_{n= 1}^\infty |C_n|^2+
\sum_{n=1}^\infty\ |D_{n}|^2|p_n|^4+|\D|^2 \bigg)
\,,
\end{equation}
where $c_1(T)$ is a positive constant. 
\end{theorem}
\begin{Proof}
Since $T>2\pi/\ga$, there exists $0<\varepsilon<1$ such that
$T>\frac{2\pi}{\ga(1-\varepsilon)}$. By applying Theorems \ref{th:Inversa} and \ref{th:Diretta}, there exist $n_0\in\N$,
$c(T,\varepsilon)>0$  and $c(T)>0$ such that if $R_n=C_n=D_n=0$ for
$n< n_0$, then we have
\begin{multline*}
c(T,\varepsilon)\bigg(\sum_{n= n_0}^\infty |C_n|^2+
\sum_{n=n_0}^\infty\ |D_{n}|^2|p_n|^4+|\D|^2 \bigg)
\\
\le
\int_{0}^{T} \big(|u_1(t)|^2+|u_2(t)|^2\big)\ dt
\\
\le
c(T)\bigg(\sum_{n= n_0}^\infty |C_n|^2+
\sum_{n=n_0}^\infty\ |D_{n}|^2|p_n|^4+|\D|^2 \bigg)
\,.
\end{multline*}
Finally, thanks to Proposition \ref{pr:haraux-inv} we can conclude.
\end{Proof}
Before to prove the direct inequality, we need a result, which allows us to recover a finite number of terms in the Fourier series.
\begin{proposition}\label{pr:haraux-dir}
Assume
\begin{equation}\label{eq:alldntris}
|d_n|\le m_2|p_n|^2\,\qquad\forall n\in\N
\qquad (m_2>0)
\,,
\end{equation}
\begin{equation}\label{eq:a0bis}
|p_n|\ge a_0\qquad\forall n\in\N
\qquad (a_0>0)\,,
\end{equation}
and that there exists $n_0\in\N$ such that 
for any sequences $\{R_n\}$, $\{C_n\}$ and $\{D_n\}$ verifying
\begin{equation*}
R_n=C_n=D_n=0\qquad\mbox{for any}\quad n< n_0\,,
\end{equation*}
the estimate
\begin{equation}\label{eq:haraux-dir}
\int_{-T}^{T}  \big(|u_1(t)|^2+|u_2(t)|^2\big)\ dt
\le
c_2\bigg(\sum_{n= n_0}^\infty|C_n|^2+
\sum_{n=n_0}^\infty\ |D_{n}|^2|p_n|^4 +|\D|^2 \bigg)
\end{equation} 
is satisfied for some  $c_2>0$.
Then, 
for any sequences $\{R_n\}$, $\{C_n\}$ and $\{D_n\}$ such that
\begin{equation}\label{eq:ha4bis}
|R_n|\le \mu|C_n|\qquad\mbox{for any}\quad n<n_0 
\qquad (\mu>0)\,,
\end{equation}
the  estimate  
\begin{equation}\label{eq:haraux-dir0}
\int_{-T}^{T} \big(|u_1(t)|^2+|u_2(t)|^2\big)\ dt
\le
C_2\bigg(\sum_{n=1}^\infty|C_n|^2+
\sum_{n=1}^\infty\ |D_{n}|^2|p_n|^4 +|\D|^2 \bigg)
\end{equation} 
holds for some  $C_2>0$.
\end{proposition}

\begin{Proof}
Let $\{R_n\}$, $\{C_n\}$ and $\{D_n\}$ be  arbitrary sequences and 
assume that  (\ref{eq:ha4bis}) holds.
If we use \eqref{eq:haraux-dir}, then we have  
\begin{multline}\label{eq:haraux-direp2}
\int_{-T}^{T} \Big|\sum_{n=n_0}^{\infty}R_ne^{r_nt}+C_ne^{i\om_nt}+\overline{C_n}e^{-i\overline{\omega_n}t}
+D_{n}e^{i p_{n}t}+\overline{D_{n}}e^{-i \overline{p_n}t}\bigg|^2  d t
\\
+\int_{-T}^{T} \Big|\sum_{n=n_0}^{\infty}\Big(d_nD_{n}e^{i p_{n}t}+\overline{d_n}\overline{D_{n}}e^{-i \overline{p_n}t}\Big)
+\D e^{-\eta t}\bigg|^2 d t
\\
\le c_2\bigg(\sum_{n= n_0}^\infty|C_n|^2+
\sum_{n=n_0}^\infty\ |D_{n}|^2|p_n|^4 +|\D|^2 \bigg)\,.
\end{multline}
Now, we  will prove that
\begin{multline}\label{eq:haraux-direp3}
\int_{-T}^{T}\Big|\sum_{n=1}^{n_0-1}
R_ne^{r_nt}+C_ne^{i\om_nt}+\overline{C_n}e^{-i\overline{\omega_n}t}
+D_{n}e^{i p_{n}t}+\overline{D_{n}}e^{-i \overline{p_n}t}
\Big|^2 \ dt
\\
+\int_{-T}^T\Big|\sum_{n=1}^{n_0-1}
d_nD_{n}e^{i p_{n}t}
+\overline{d_n}\overline{D_{n}}e^{-i \overline{p_n}t}\Big|^2 d t
\le c'_2\sum_{n=1}^{n_0-1}\Big(|C_n|^2+\ |D_{n}|^2|p_n|^4  \Big)\,,
\end{multline}
for some constant $c'_2>0$. Indeed, by applying the Cauchy-Schwarz inequality we get
\begin{multline*}
\Big|\sum_{n=1}^{n_0-1}
R_ne^{r_nt}+C_ne^{i\om_nt}+\overline{C_n}e^{-i\overline{\omega_n}t}
+D_{n}e^{i p_{n}t}+\overline{D_{n}}e^{-i \overline{p_n}t}
\Big|^2
\\
\le
\bigg(\sum_{n=1}^{n_0-1}|R_n|e^{r_nt}+2|C_n|e^{-\Im\om_nt}+2|D_n|e^{-\Im p_nt}\bigg)^2\\
\le
12(n_0-1)\sum_{n=1}^{n_0-1}\Big(|R_n|^2e^{2r_nt}+|C_n|^2e^{-2\Im\om_nt}+|D_n|^2e^{-2\Im p_nt}\Big)
\end{multline*}
and in view also of \eqref{eq:alldntris}
\begin{multline*}
\Big|\sum_{n=1}^{n_0-1}
d_nD_{n}e^{i p_{n}t}
+\overline{d_n}\overline{D_{n}}e^{-i \overline{p_n}t}\Big|^2
\le
4\bigg(\sum_{n=1}^{n_0-1}|d_nD_n|e^{-\Im p_nt}\bigg)^2
\\
\le
4(n_0-1)\sum_{n=1}^{n_0-1}|d_n|^2|D_n|^2e^{-2\Im p_nt}
\le
4(n_0-1)m^2_2\sum_{n=1}^{n_0-1}|p_n|^4|D_n|^2e^{-2\Im p_nt}
\,.
\end{multline*}
If we use the previous inequalities, (\ref{eq:ha4bis}) and \eqref{eq:a0bis}, then we get
\begin{multline*}
\int_{-T}^{T}\Big|\sum_{n=1}^{n_0-1}
R_ne^{r_nt}+C_ne^{i\om_nt}+\overline{C_n}e^{-i\overline{\omega_n}t}
+D_{n}e^{i p_{n}t}+\overline{D_{n}}e^{-i \overline{p_n}t}
\Big|^2 \ dt
\\
+\int_{-T}^T\Big|\sum_{n=1}^{n_0-1}
d_nD_{n}e^{i p_{n}t}
+\overline{d_n}\overline{D_{n}}e^{-i \overline{p_n}t}\Big|^2 d t
\le 12(n_0-1)\sum_{n=1}^{n_0-1}|C_n|^2\int_{-T}^{T}\big(\mu^2e^{2r_nt}+e^{-2\Im\om_nt}\big)\ dt
\\
+4(n_0-1)(3a_0^{-4}+m^2_2)\sum_{n=1}^{n_0-1}|p_n|^4|D_n|^2\int_{-T}^{T}e^{-2\Im p_nt}\ dt
\,,
\end{multline*}
whence (\ref{eq:haraux-direp3}) follows with
$
\displaystyle
c'_2=12(n_0-1)|\max_{n<n_0}\Big\{\int_{-T}^{T}\big(\mu^2e^{2r_nt}+e^{-2\Im\om_nt}+(a_0^{-4}+m^2_2)e^{-2\Im p_nt}\big)\ dt\Big\}\,.
$ 

Finally, from (\ref{eq:haraux-direp2}) and (\ref{eq:haraux-direp3}) we deduce that
\begin{multline*}
\int_{-T}^{T} \big(|u_1(t)|^2+|u_2(t)|^2\big)\ dt
\\
\le
2c_2\bigg(\sum_{n= n_0}^\infty|C_n|^2+
\sum_{n=n_0}^\infty\ |D_{n}|^2|p_n|^4 +|\D|^2 \bigg)
+2c'_2\bigg(\sum_{n=1}^{n_0-1}|C_n|^2+
\sum_{n=1}^{n_0-1}\ |D_{n}|^2|p_n|^4  \bigg)\,,
\end{multline*}
so (\ref{eq:haraux-dir0}) holds with 
$
C_2=2\max\{c_2,c'_2\}\,.
$
\end{Proof}
Finally, we are in a position to prove the direct inequality.
\begin{theorem}\label{th:diringham}
Assume $p_n\not=0$,
\begin{equation*}
\lim_{n\to\infty}(\Re p_{n+1}-\Re p_{n})=+\infty\,,
\end{equation*}
\begin{equation*}
\lim_{n\to\infty}\Im p_{n}=0\,,
\end{equation*}
and for some $\gamma>0$, $\alpha\in\R$, $n'\in\N$, $\mu>0$, $\nu> 1/2$, $m_1\,,m_2>0$
\begin{equation*}
\liminf_{n\to\infty}({\Re}\om_{n+1}-{\Re}\om_{n})=\gamma\,,
\end{equation*}
\begin{equation*}
\lim_{n\to\infty}{\Im}\om_n=\alpha
\,,
\qquad
r_n\le -{\Im}\om_n\,\qquad\forall\ n\ge n'\,,
\end{equation*}
\begin{equation*}
|R_n|\le \frac{\mu}{n^{\nu}}|C_n|\,\quad\forall\ n\ge n'\,,
\qquad
|R_n|\le \mu|C_n|\,\quad\forall\ n\le n'\,,
\end{equation*}
\begin{equation*}
m_1|p_n|^2\le|d_n|\le m_2|p_n|^2\,\qquad\forall n\in\N\,.
\end{equation*}
Then, for any
$T>\pi/\ga$ we have 
\begin{equation}\label{eq:diringham}
 \int_{-T}^{T}  \big(|u_1(t)|^2+|u_2(t)|^2\big)\ dt
 \le
 c_2(T)\bigg(\sum_{n= 1}^\infty |C_n|^2+
\sum_{n=1}^\infty\ |D_{n}|^2|p_n|^4 \bigg)\,,
 \end{equation}
where $c_2(T)$ is a positive constant\,.
\end{theorem}
\begin{Proof}
Since $T>\pi/\ga$, there exists $0<\varepsilon<1$ such that
$T>\frac{\pi}{\ga\sqrt{1-\varepsilon}}$. By applying Theorem \ref{th:Diretta}, 
there exist $c(T)>0$ and $n_0=n_0(\varepsilon)\in\N$ such that if $C_n=D_n=0$ for any
$n< n_0$, then  we have
\begin{equation*}
\int_{-T}^{T} \big(|u_1(t)|^2+|u_2(t)|^2\big)\ dt
\le
c(T)\bigg(\sum_{n= n_0}^\infty |C_n|^2+
\sum_{n=n_0}^\infty\ |D_{n}|^2|p_n|^4+|\D|^2 \bigg)
\,.
\end{equation*}
We can use Proposition \ref{pr:haraux-dir} to obtain,
for any arbitrary sequences $\{R_n\}$, $\{C_n\}$ and $\{D_n\}$, 
\begin{equation}\label{eq:diringham1}
\int_{-T}^{T} \big(|u_1(t)|^2+|u_2(t)|^2\big)\ dt
\le
C_2\bigg(\sum_{n=1}^\infty|C_n|^2+
\sum_{n=1}^\infty\ |D_{n}|^2|p_n|^4 +|\D|^2 \bigg)
\end{equation} 
for some  $C_2>0$.
Moreover, if we take
\begin{equation*}
\D=-\frac\beta A\sum_{n=1}^{\infty}\Re p_{n}
\Big(\frac{D_{n}}{\eta+ip_n}+\frac{\overline{D_{n}}}{\eta-ip_n}
\Big)\,,
\end{equation*}
then, for some $C>0$, we have
\begin{equation}\label{eq:stimaG2}
|\D|^2
\le C
\sum_{n=1}^\infty|D_{n}|^2|p_{n}|^4\,.
\end{equation}
Indeed, we observe that
\begin{equation*}
\Big|\sum_{n=1}^{\infty}\Re p_n
\Big(\frac{D_{n}}{\eta+ip_n}+\frac{\overline{D_{n}}}{\eta-i\overline{p_n}}\Big)\Big|
\le2\sum_{n=1}^{\infty}|p_{n}|\big|\Re\frac{D_{n}}{\eta+ip_n}\big|
\le2\sum_{n=1}^{\infty}|p_{n}|\frac{|D_{n}|}{|\eta+ip_n|}\,,
\end{equation*}
whence, thanks to $\lim_{n\to\infty} \frac{\Re p_n}{n}=+\infty$, we get
\begin{multline*}
\Big|\sum_{n=1}^{\infty}\Re p_n
\Big(\frac{D_{n}}{\eta+ip_n}+\frac{\overline{D_{n}}}{\eta-i\overline{p_n}}\Big)\Big|^2
\le
4\sum_{n,m=1}^{\infty}\frac{|D_{n}||p_{n}|}{|\eta+ip_m|}
\frac{|D_{m}||p_{m}|}{|\eta+ip_n|}
\\
\le2\sum_{n=1}^\infty|D_{n}|^2|p_{n}|^2\sum_{m=1}^{\infty}\frac{1}{(\eta-\Im p_m) ^2+\Re p_m^2}
+2\sum_{m=1}^\infty|D_{m}|^2|p_{m}|^2\sum_{n=1}^{\infty}\frac{1}{(\eta-\Im p_n) ^2+\Re p_n^2}
\\
=4\sum_{n=1}^{\infty}\frac{1}{(\eta-\Im p_n) ^2+\Re p_n^2}\sum_{n=1}^\infty|D_{n}|^2|p_{n}|^2\,.
\end{multline*}
Therefore, taking into account \eqref{eq:a0},  we have that \eqref{eq:stimaG2} holds true
with
$$C=4a_0^2\frac\beta A\sum_{n=1}^{\infty}\frac{1}{(\eta-\Im p_n) ^2+\Re p_n^2}.$$
In conclusion, from \eqref{eq:diringham1} and \eqref{eq:stimaG2} it follows
\eqref{eq:diringham}.
\end{Proof}

\section{A reachability result}

Finally, by applying our abstract results of Sections  \ref{se:specana} and  \ref{se:invdir} we are able to show our reachability result for  wave--Petrovsky coupled systems with a memory term.
\begin{theorem}\label{th:reachres}
Let $\eta>3\beta/2$. For  any $T>2\pi$ and
$$(u_{10},u_{11},u_{20},u_{21})\in  L^{2}(0,\pi)\times H^{-1}(0,\pi)\times H_0^{1}(0,\pi)\times H^{-1}(0,\pi)\,,$$
 there exist $g_i\in L^2(0,T)$, $i=1,2$, such that the weak solution  $(u_1,u_2)$ of system 
\begin{equation}\label{eq:problem-usix}
\begin{cases}
\displaystyle 
u_{1tt}(t,x) -u_{1xx}(t,x)+\beta\int_0^t\ e^{-\eta(t-s)} u_{1xx}(s,x)ds+Au_2(t,x)= 0\,,
\\
\phantom{u_{1tt}(t,x) -u_{1xx}(t,x)+\int_0^t\ k(t-s) u_{1xx}(s,x)ds+Au_2(t,x)= 0}
t\in (0,T)\,,\quad x\in(0,\pi)
\\
\displaystyle
u_{2tt}(t,x) +u_{2xxxx}(t,x)+Bu_1(t,x)= 0
\,,
\end{cases}
\end{equation} 
with null initial conditions 
\begin{equation}
u_1(0,x)=u_{1t}(0,x)=u_2(0,x)=u_{2t}(0,x)=0\qquad  x\in(0,\pi)\,,
\end{equation} 
and boundary conditions
\begin{equation}\label{eq:bound-u1r}
u_1(t,0)=0\,,\qquad u_1(t,\pi)=g_1(t)\qquad t\in (0,T) \,,
\end{equation}
\begin{equation}\label{eq:bound-u2r}
u_2(t,0)=u_{2xx}(t,0)=u_2(t,\pi)=0\,,\qquad u_{2xx}(t,\pi)=g_2(t)\qquad t\in (0,T) \,.
\end{equation}
verifies the final conditions
\begin{equation}\label{eq:findataT}
u_1(T,x)=u_{10}(x)\,,\qquad u_{1t}(T,x)=u_{11}(x)\,,
\qquad x\in(0,\pi)\,,
\end{equation}
\begin{equation}\label{eq:findataT1}
u_2(T,x)=u_{20}(x)\,,\qquad u_{2t}(T,x)=u_{21}(x)\,,
\qquad x\in(0,\pi)\,.
\end{equation}
\end{theorem}
\begin{Proof}
To prove our claim, we will apply the Hilbert Uniqueness Method described in Section \ref{se:HUM}.
Let 
$
X= L^2(0,\pi )
$
be endowed with the usual scalar product and norm 
$$
\|v\|:=\left(\int_0^\pi |v(x)|^{2}\ dx\right)^{1/2}\qquad
v\in L^2(0,\pi)\,.
$$
We consider the operator $L:D(L)\subset X\to X$ defined by
$$
\begin{array}{l}
D(L)=H^2(0,\pi )\cap H_0^1(0,\pi ) \\
\\
Lv=\displaystyle -v_{xx}\qquad v\in D(L)\,.
\end{array}
$$
It is well known that $L$ is a self-adjoint positive
 operator on $X$ with dense domain $D(L)$, $\{n^2\}_{n\ge1}$ is the sequence  of  eigenvalues for  $L$ and 
  $\{\sin(nx)\}_{n\ge1}$ is the sequence of the corresponding eigenvectors.
We can apply our spectral analysis (see Section \ref{se:specana}) to the adjoint system of (\ref{eq:problem-usix}). Indeed, the adjoint system  is given by 
\begin{equation}\label{eq:adjointr}
\begin{cases}
\displaystyle 
z_{1tt}(t,x) -z_{1xx}(t,x)+\beta\int_t^T\ e^{-\eta(s-t)} z_{1xx}(s,x)ds+Bz_2(t,x)= 0\,,\\
\phantom{u_{1tt}(t,x) -u_{1xx}(t,x)+\int_0^t\ k(t-s) u_{1xx}(s,x)ds+Au_2(t,x)= 0\,,\qquad}
t\in (0,T)\,,\quad x\in(0,\pi)
\\
\displaystyle
z_{2tt}(t,x) +z_{2xxxx}(t,x)+Az_1(t,x)= 0\,,
\\
z_1(t,0)=z_1(t,\pi)=z_2(t,0)=z_2(t,\pi)=z_{2xx}(t,0)=z_{2xx}(t,\pi)=0\qquad t\in [0,T]\,,
\end{cases}
\end{equation}
with  final data 
\begin{equation} \label{eq:finalr}
z_1(T,\cdot)=z_{10}\,,\quad z_{1t}(T,\cdot)=z_{11}\,,\quad 
z_2(T,\cdot)=z_{20}\,,\quad z_{2t}(T,\cdot)=z_{21}\,,
\end{equation}
where
\begin{equation}\label{eq:fi1}
z_{10}(x)=\sum_{n=1}^{\infty}\alpha_{1n}\sin(nx)\,,\qquad
              \|z_{10}\|^2_{H^1_0}=\sum_{n=1}^\infty\alpha^2_{1n} n^2
\,,
\end{equation} 
\begin{equation}\label{eq:fi2}
z_{11}(x)=\sum_{n=1}^{\infty}\rho_{1n}\sin(nx)\,,\qquad
\|z_{11}\|^2= \sum_{n=1}^\infty\rho^2_{1n}
              \,,
\end{equation} 
\begin{equation}\label{eq:fi3}
z_{20}(x)=\sum_{n=1}^{\infty}\alpha_{2n}\sin(nx)\,,\qquad
\|z_{20}\|^2_{H^1_0}=\sum_{n=1}^\infty\alpha^2_{2n} n^2
       \,,
\end{equation} 
\begin{equation}\label{eq:fi4}
z_{21}(x)=\sum_{n=1}^{\infty}\rho_{2n}\sin(nx)\,,\qquad
        \|z_{21}\|_{H^{-1}}^2=\sum_{n=1}^\infty\frac{\rho^2_{2n}} {n^2}
\,.
\end{equation} 
The solution $(z_1,z_2)$ of system \eqref{eq:adjointr}-\eqref{eq:finalr} can be written in the following way (see formulas \eqref{eq:vsum}--\eqref{eq:esprD}): for any $(t,x)\in [0,T]\times [0,\pi]$
\begin{equation*}
z_1(t,x)=\sum_{n=1}^{\infty}\Big(R_ne^{r_n(T-t)}+C_ne^{i\om_n(T-t)}
+\overline{C_n}e^{-i\overline{\om_n}(T-t)}+D_ne^{ip_n(T-t)}+\overline{D_n}e^{-i\overline{p_n}(T-t)}\Big)\sin(n x)
\,,
\end{equation*}
\begin{multline*}
z_2(t,x)=\sum_{n=1}^{\infty}
\Big(d_nD_{n}e^{i p_{n}(T-t)}
+\overline{d_n}\overline{D_{n}}e^{-i \overline{p_n}(T-t)}\Big)\sin(n x)
\\
-\frac\beta Ae^{-\eta (T-t)}\sum_{n=1}^{\infty}\Re p_{n}
\Big(\frac{D_{n}}{\eta+ip_n}+\frac{\overline{D_{n}}}{\eta-ip_n}
\Big)\sin(n x)
\,.
\end{multline*}
In particular,  by estimates \eqref{eq:|Cj2|} and \eqref{eq:|Cj4|} we have the following relationships between the coefficients $C_n$, $D_n$ and the Fourier coefficients of the final data:

\begin{equation}\label{eq:rel1}
                   \sum_{n=1}^\infty    n^2\vert C_{n}\vert^2
               \asymp
              \sum_{n=1}^\infty\alpha^2_{1n} n^2+ \sum_{n=1}^\infty\rho^2_{1n}\,,
  \end{equation}
    \begin{equation}\label{eq:rel2}         
           \sum_{n=1}^\infty n^2\vert D_{n}\vert^2\vert p_{n}\vert^4
      \asymp
        \sum_{n=1}^\infty\alpha^2_{2n} n^2+ \sum_{n=1}^\infty\frac{\rho^2_{2n}} {n^2}\,.
       \end{equation}
Moreover, for any $t\in [0,T]$
\begin{equation*}
z_{1x}(t,\pi)=\sum_{n=1}^{\infty}(-1)^n n\Big(R_ne^{r_n(T-t)}+C_ne^{i\om_n(T-t)}
+\overline{C_n}e^{-i\overline{\om_n}(T-t)}+D_ne^{ip_n(T-t)}+\overline{D_n}e^{-i\overline{p_n}(T-t)}\Big)
\,,
\end{equation*}
\begin{multline*}
z_{2x}(t,\pi)=\sum_{n=1}^{\infty}
(-1)^n n\Big(d_nD_{n}e^{i p_{n}(T-t)}
+\overline{d_n}\overline{D_{n}}e^{-i \overline{p_n}(T-t)}\Big)
\\
-\frac\beta Ae^{-\eta (T-t)}\sum_{n=1}^{\infty}(-1)^n n\Re p_{n}
\Big(\frac{D_{n}}{\eta+ip_n}+\frac{\overline{D_{n}}}{\eta-ip_n}
\Big)
\,.
\end{multline*}
We can apply Theorems \ref{th:inv.ingham1} and \ref{th:diringham}  to $(z_{1x}(t,\pi),z_{2x}(t,\pi))$.
Indeed, thanks to inequalities  \eqref{eq:inv.ingham} and \eqref{eq:diringham}  we have   
\begin{equation*}
 \int_{0}^{T}  \big(|z_{1x}(t,\pi)|^2+|z_{2x}(t,\pi)|^2\big)\ dt
  \asymp
 \sum_{n= 1}^\infty n^2\Big( |C_n|^2+
 |D_{n}|^2|p_n|^4 \Big)\,.
 \end{equation*}
Therefore, taking into account \eqref{eq:rel1}-\eqref{eq:rel2} and the expressions of the norms of the final data, see \eqref{eq:fi1}-\eqref{eq:fi4},
we get
\begin{equation}\label{eq:}
 \int_{0}^{T}  \big(|z_{1x}(t,\pi)|^2+|z_{2x}(t,\pi)|^2\big)\ dt
  \asymp
 \|z_{10}\|^2_{H^1_0}+\|z_{11}\|^2
 +\|z_{20}\|^2_{H^1_0}+\|z_{21}\|_{H^{-1}}^2
 \,.
 \end{equation}
Finally, Theorem \ref{th:uniqueness1} holds true and
the space $F$ introduced at the end of Section \ref{se:HUM} is 
$$H^1_0(0,\pi)\times L^2(0,\pi)\times H^1_0(0,\pi)\times H^{-1}(0,\pi)\,,$$ 
so, if we take
\begin{equation*}
g_1(t)=z_{1x}(t,\pi)-\beta\int_t^T\ e^{-\eta(s-t)}z_{1x}(s,\pi)ds\,,
\qquad
g_2(t)=-z_{2x}(t,\pi)\,,
\end{equation*}
then our proof is complete.
\end{Proof}

\end{document}